\documentclass[11pt,reqno]{amsart}
\usepackage{amssymb,amscd,amsbsy,mathrsfs}
\newcommand{\lb}{\linebreak}

\renewcommand{\a}{\alpha}
\renewcommand{\b}{\beta}

\renewcommand{\d}{\delta}
\newcommand{\e}{\varepsilon}
\newcommand{\vk}{\varkappa}
\newcommand{\z}{\zeta}

\renewcommand{\l}{\lambda}

\newcommand{\s}{\sigma}
\renewcommand{\t}{\tau}

\newcommand{\f}{\varphi}
\renewcommand{\o}{\omega}

\newcommand{\G}{\Gamma}
\newcommand{\D}{\Delta}
\renewcommand{\L}{\Lambda}

\renewcommand{\O}{\Omega}

\newcommand{\n}{\nabla}

\newcommand{\A}{{\mathscr A}}

\newcommand{\E}{{\mathscr E}}
\newcommand{\cd}{{\mathscr D}}
\newcommand{\F}{{\mathscr F}}

\newcommand{\h}{{\mathscr H}}

\newcommand{\K}{{\mathscr K}}

\newcommand{\M}{{\mathscr M}}

\newcommand{\X}{{\mathscr X}}

\newcommand{\1}{{\bf 1}}

\newcommand{\C}{{\Bbb C}}
\newcommand{\T}{{\Bbb T}}
\newcommand{\pp}{{\Bbb P}}
\newcommand{\dd}{{\Bbb D}}
\newcommand{\R}{{\Bbb R}}
\newcommand{\Z}{{\Bbb Z}}

\newcommand{\0}{{\boldsymbol{0}}}

\newcommand{\bs}{\boldsymbol}

\newcommand{\m}{{\boldsymbol m}}
\newcommand{\bS}{{\boldsymbol S}}

\newcommand{\rf}[1]{(\ref{#1})}

\newcommand{\df}{\stackrel{\mathrm{def}}{=}}
\newcommand{\dist}{\operatorname{dist}}

\newcommand{\re}{\operatorname{Re}}
\newcommand{\spn}{\operatorname{span}}
\newcommand{\supp}{\operatorname{supp}}
\newcommand{\clos}{\operatorname{clos}}

\newcommand{\trace}{\operatorname{trace}}
\newcommand{\rank}{\operatorname{rank}}
\newcommand{\const}{\operatorname{const}}

\newcommand{\eeq}{\end{equation}}
\newcommand{\beq}{\begin{equation}}
\newcommand{\bay}{\begin{eqnarray}}
\newcommand{\ba}{\begin{align*}}
\newcommand{\ea}{\end{align*}}
\newcommand{\ey}{\end{eqnarray}}
\newcommand{\bey}{\begin{eqnarray*}}
\newcommand{\eey}{\end{eqnarray*}}

\newcommand{\eq}{\Leftrightarrow}
\newcommand{\imp}{\Rightarrow}
\newcommand{\be}{\infty}

\newcommand{\bl}{\blacksquare}
\newcommand{\ess}{\operatorname{ess}}

\newcommand{\Pf}{{\bf Proof. }}
\newcommand{\im}{\operatorname{Im}}
\renewcommand{\re}{\operatorname{Re}}
\newcommand{\ov}{\overline}

\newtheorem{thm}{\hspace{\parindent}Theorem}[section]

\newtheorem{cor}[thm]{\hspace{\parindent}Corollary}
\newtheorem{lem}[thm]{\hspace{\parindent}Lemma}

\usepackage{amsmath}

\newcommand\Li{{\rm Lip}}
\newcommand\OL{{\rm OL}}
\newcommand\CM{{\rm CM}}
\newcommand\fM{\frak M}
\newcommand\fF{\frak F}

\newcommand\dg{\frak D}

\newcommand\CA{{\rm C}_{\rm A}}
\newcommand\CL{{\rm CL}}
\newcommand\OD{{\rm OD}}
\newcommand\sgn{\operatorname{sgn}}
\newcommand{\nn}{{\Bbb N}}
\newcommand\mB{\mathcal{B}}
\newcommand\mD{\mathcal{D}}
\newcommand\mM{\mathcal{M}}
\newcommand\mS{\mathcal{S}}
\newcommand\mT{\mathcal{T}}
\newcommand\mE{\mathcal{E}}
\newcommand\mP{\mathcal{P}}

\newcommand{\dlR}{{\operatorname{Aut}}(\widehat{\Bbb R})}
\newcommand{\lC}{{\operatorname{Aut}}(\Bbb C)}
\newcommand{\lR}{{\operatorname{Aut}}(\Bbb R)}
\newcommand{\dlC}{{\operatorname{Aut}}(\widehat{\Bbb C})}

\newcommand\zh{|\!\!|\!|}

\begin{document}

\

\newcommand{\vse}{\vspace{.2in}}
\numberwithin{equation}{section}

\title{Operator Lipschitz functions}
\author{A.B. Aleksandrov and V.V. Peller}
\thanks{Research of the first author is partially supported by grant RFFI 14-01-00198;
Research of the second author is partially supported by NSF grant DMS 130092}

\begin{abstract}
The purpose of this survey is a comprehensive study of operator Lip\-schitz functions. A continuous function $f$ on the real line $\R$ os called operator Lipschitz if $\|f(A)-f(B)\|\le\const\|A-B\|$ for arbitrary self-adjoint operators $A$ and $B$. We give sufficient conditions and necessary conditions for operator Lipschitzness. We also study the class of operator differentiable functions on $\R$. Next, we consider operator Lipschitz functions on closed subsets of the plane and introduce the class of commutator Lipschitz functions on such subsets. An important role for the study of such classes of functions is played by double operator integrals and Schur multipliers.
\end{abstract}

\maketitle

{\bf
\footnotesize
\tableofcontents
\normalsize
}

\vspace*{-1cm}

\setcounter{section}{0}
\section{\bf Introduction}
\setcounter{equation}{0}
\label{intr}

\medskip

One of the most important problems in perturbation theory is the study of the question to which extent a function  $f(A)$ of an operator $A$ can change under small perturbations of the operator. In particular, in a natural way we arrive at the problem to describe the class of continuous functions $f$ on the real line $\R$ such that the following inequality holds:
\bay
\label{OLf}
\|f(A)-f(B)\|\le\const\|A-B\|
\ey
for arbitrary (bounded) self-adjoint operators $A$ and $B$ on Hilbert space. Such functions are called {\it operator Lipschitz}. Recall that functions of self-adjoint (normal) operators are defined as the integrals of theses functions with respect to the spectral measures of the operators, see \cite{R}.

We will denote the class of operator Lipschitz functions on $\R$ by $\OL(\R)$. Note that if $f$ is an operator Lipschitz function, then inequality \rf{OLf} also holds for unbounded self-adjoint operators $A$ and $B$ with bounded difference, see Theorem \ref{anbnrn} below; moreover, the constant on the right remains the same. The minimal value of this constant is, by definition, the norm  $\|f\|_{\OL}=\|f\|_{\OL(\R)}$ of the function $f$ 
in the space $\OL(\R)$
(strictly speaking, it is a seminorm that that becomes a norm after the identification of functions that differ from each other be a constant function). 

Clearly, if $f$ is an operator Lipschitz function, then it is {\it Lipschitz}, i.e.,
$$
|f(x)-f(y)|\le\const|x-y|,
$$
for arbitrary real $x$ and $y$ (we use the notation $\Li(\R)$ for the class of Lipschitz functions on $\R$). The converse is wrong. Farforovskaya constructed in \cite{F1} an example of a Lipschitz function that is not operator Lipschitz. Later it was shown in \cite{Mc} and \cite{Kat} that the Lipschitz function $x\mapsto|x|$ is not operator Lipschitz.

Operator Lipschitz functions play an important role in operator theory and mathematical physics. In particular, they appear when studying the applicability of the Lifshits--Krein trace formula:
\bay
\label{LiKr}
\trace\big((f(A)-f(B)\big)=\int_\R f'(t)\xi(t)\,dt
\ey
(see \cite{Kr}).
Here $A$ and $B$ are self-adjoint operators on Hilbert space such that $A-B$ is a trace class operator (i.e., $A-B\in\bS_1$) and $\xi$ is a function of class $L^1(\R)$ ({\it the spectral shift function}), which is determined only by $A$ and $B$. Clearly, the right-hand side of \rf{LiKr} makes sense for an arbitrary Lipschitz function $f$. As for the left-hand side, as seen from the example of Farforovskaya in \cite{F2}, the conditions $A-B\in\bS_1$ and $f\in\Li(\R)$ do not guarantee that $f(A)-f(B)\in\bS_1$. Thus, for the applicability of trace formula \rf{LiKr} for all self-adjoint operators with trace class difference, one has to impose a stronger assumption on $f$. At least $f$ has to possess the following property:
\bay
\label{yad}
A-B\in\bS_1\quad\Longrightarrow\quad f(A)-f(B)\in\bS_1
\ey
for self-adjoint operators $A$ and $B$. For a function $f$ 
on $\R$, property \rf{yad}
holds for arbitrary (not necessarily bounded) self-adjoint operators if and only if $f$ is operator Lipschitz (see Theorem \ref{tsentrez} below).
It turns out (see the recent paper \cite{Pe7}) that the operator Lipschitzness of $f$ 
is not only necessary for the validity of trace formula \rf{LiKr} for arbitrary not necessarily bounded) self-adjoint operators $A$ and $B$ with trace class difference, but also sufficient.

The class of operator Lipschitz functions possesses certain peculiar properties. For example, operator Lipschitz functions must be differentiable everywhere, but not necessarily continuously differentiable (see Theorem \ref{pro} and Example 7 in \S\:\ref{prim}).

It turns out that operator Lipschitzness can be characterized in terms of Schur multipliers (see \S\:\ref{razdrazn}).
We will see that a continuous function $f$ on $\R$ is operator Lipschitz if and only if it is differentiable everywhere and the divided difference $\dg f$,
$$
(\dg f)(x,y)\df\frac{f(x)-f(y)}{x-y},\quad x,\:y\in\R,
$$
is a Schur multiplier.

In a similar way one can consider the same problem for functions on the circle and for unitary operators. A continuous function $f$ on the unit circle $\T$ is called {\it operator Lipschitz} if
$
\|f(U)-f(V)\|\le\const\|U-V\|
$
for arbitrary unitary operators $U$ and $V$.

In Chapter I of this survey we discuss necessary conditions and sufficient conditions for functions on the line $\R$ and on the circle $\T$ to be operator Lipschitz.
Note that in the case of self-adjoint operators a key role is played by the inequality
\bay
\label{sig}
\|f(A)-f(B)\|\le\const\s\|f\|_{L^\be}\|A-B\|
\ey
for arbitrary self-adjoint operators $A$ and $B$ with bounded difference and an arbitrary  bounded function $f$ on $\R$ whose Fourier transform is supported in  $[-\s,\s]$, $\s>0$. This inequality was obtained in \cite{Pe1} and \cite{Pe3}.
Later, it was shown in \cite{AP4} that inequality \rf{sig} holds with constant 1. 

By analogy with operator Lipschitz functions, it would be natural to consider operator 
H\"older functions. Let $0<\a<1$. We say that a function $f$ on $\R$ {\it is operator 
H\"older of order} $\a$ if the inequality
$$
\|f(A)-f(B)\|\le\const\|A-B\|^\a
$$
holds
for arbitrary self-adjoint operators $A$ and $B$ on Hilbert space. However (see \S\:\ref{oHold}), the situation here is tremendously different from the operator Lipschitz estimates: a  $f$ function is operator 
H\"older of order $\a$ if and only if it belongs to the class $\L_\a(\R)$ of {\it 
H\"older functions of order} $\a$, i.e.,
$|f(x)-f(y)|\le\const|x-y|^\a,\quad x,\;y\in\R$.

In Chapter II we discuss double operator integrals, i.e., expressions of the form
$$
\iint\Phi(x,y)\,dE_1(x)T\,dE_2(y).
$$
Here $\Phi$ is a bounded measurable function, $T$ is a bounded linear operator on Hilbert space, and $E_1$ and $E_2$ are spectral measures.
Double operator integrals appeared in the paper by Yu.L. Daletskii and S.G. Krein \cite{DK} and were studied systematically by M.S. Birman and M.Z. Solomyak in \cite{BS1}--\cite{BS3}. It is in those papers it became clear what important role double operator integrals play in perturbation theory. Double operator integrals for arbitrary bounded linear operators $T$ are defined in the case when the function $\Phi$ is a {\it Schur multiplier} with respect to $E_1$ and $E_2$. In Chapter II we study the space of such Schur multipliers. First, we study so-called {\it discrete Schur multipliers} and then we use them to study Schur multipliers with respect to spectral measures.

Next, in Chapter III we consider the class $\OL(\fF)$ of {\it operator Lipschitz functions on an arbitrary closed subset} $\fF$ of the complex plant $\C$,
which consists of continuous functions $f$ on $\fF$ such that
\bay
\label{neognor}
\|f(N_1)-f(N_2)\|\le\const\|N_1-N_2\|
\ey
for arbitrary normal operators $N_1$ and $N_2$ whose spectra are contained in
$\fF$. We also study in detail the class of {\it commutator Lipschitz functions} on
$\fF$, i.e., the class of continuous functions $f$ on $\fF$ such that
$$
\|f(N_1)R-Rf(N_2)\|\le\const\|N_1R-RN_2\|
$$
for every bounded linear operator $R$ and for arbitrary normal operators $N_1$ and $N_2$ with spectra in $\fF$. To study these classes of functions, we use results of Chapter II.

For the study of the class of operator Lipschitz functions on the whole plane, as in the case of self-adjoint operators, a key role is played by the following generalization of inequality \rf{sig}:
\bay
\label{nors}
\|f(N_1)-f(N_2)\|\le\const\s\|f\|_{L^\be}\|N_1-N_2\|
\ey
for arbitrary normal operators $N_1$ and $N_2$ with bounded difference and for an arbitrary bounded function $f$ on $\R^2$ whose Fourier transform is supported in 
$[-\s,\s]\times[-\s,\s]$. Note that the proof of inequality \rf{sig} obtained in \cite{Pe1} and \cite{Pe3} cannot be generalized to the case of normal operators. A new method of obtaining such estimates was found in \cite{APPS}. 

We also obtain a sufficient condition for the commutator Lipschitzness of functions on proper closed subsets of the plane in terms of Cauchy integrals of measures on the compliment of the set; it was found in \cite{A2}. We use this condition to deduce the sufficient condition by Arazy--Bartman--Friedman \cite{ABF} for the commutator Lipschitzness of functions analytic in the disk as well as its analog for the half-plane.

Finally, we study in Chapter III properties of commutator Lipschitz functions on the unit circle $\T$ that admit analytic extensions to the unit disk $\dd$; these results are grouped around the results of Kissin and Shulman of \cite{KS3}.

In the final section ``Concluding Remarks'' we mention briefly certain results that were not included in the survey.

The authors would like to express a sincere gratitude to V.S. Shulman for helpful remarks.

\medskip

\section{\bf Preliminaries and notation}
\setcounter{equation}{0}
\label{Prel}

\medskip

{\bf1. Besov classes.} Let $w$ be an infinitely differentiable function on  $\R$ such that
\bay
\label{w}
w\ge0,\quad\supp w\subset\left[1/2,2\right],\quad\mbox{and} \quad w(s)=1-w\left(\frac s2\right)\quad\mbox{for}\quad s\in[1,2].
\ey

We define the functions $W_n$, $n\in\Z$, on $\R^d$ by 
$$
\big(\F W_n\big)(x)=w\left(\frac{\|x\|}{2^n}\right),\quad n\in\Z, \quad x=(x_1,\cdots,x_d),
\quad\|x\|\df\Big(\sum_{j=1}^dx_j^2\Big)^{1/2},
$$
where $\F$ is the {\it Fourier transform} defined on $L^1\big(\R^d\big)$ by
$$
\big(\F f\big)(t)=\!\int\limits_{\R^d} f(x)e^{-{\rm i}(x,t)}\,dx,\!\quad 
x=(x_1,\cdots,x_d),
\quad t=(t_1,\cdots,t_d), \!\quad(x,t)\df \sum_{j=1}^dx_jt_j.
$$
Clearly,
$$
\sum_{n\in\Z}(\F W_n)(t)=1,\quad t\in\R^d\setminus\{0\}.
$$

With each tempered distribution $f$ in ${\mathscr S}^\prime\big(\R^d\big)$ we associate the sequence $\{f_n\}_{n\in\Z}$,
\bay
\label{fn}
f_n\df f*W_n.
\ey
The formal series $\sum_{n\in\Z}f_n$, being a Paley--Wiener type expansion of $f$, does not necessarily converge to $f$. 
First, we define the (homogeneous) Besov class $\dot B^s_{p,q}\big(\R^d\big)$,
$s\in\R$, $0<p,\,q\le\be$, as the space of distributions 
$f$ 
such that
\bay
\label{Wn}
\{2^{ns}\|f_n\|_{L^p}\}_{n\in\Z}\in\ell^q(\Z),\quad
\|f\|_{B^s_{p,q}}\df\big\|\{2^{ns}\|f_n\|_{L^p}\}_{n\in\Z}\big\|_{\ell^q(\Z)}.
\ey
In accordance with this definition, $\dot B^s_{p,q}(\R^d)$ contains all polynomials and 
$\|f\|_{B^s_{p,q}}=0$ for every polynomial $f$. Moreover, a distribution $f$ is uniquely determined by the sequence $\{f_n\}_{n\in\Z}$
modulo polynomials. It is easy to see that the series
$\sum_{n\ge0}f_n$ converges in ${\mathscr S}^\prime(\R^d)$\footnote{Here and in what follows we assume that the space  ${\mathscr S}^\prime(\R^d)$ is equipped with the weak topology
$\s\big({\mathscr S}^\prime(\R^d),{\mathscr S}(\R^d)\big)$.}. 
However, the series $\sum_{n<0}f_n$ can diverge in general. Nevertheless, it can be proved that the series
\bay
\label{ryad}
\sum_{n<0}\frac{\partial^r f_n}{\partial x_1^{r_1}\cdots\partial x_d^{r_d}}\qquad \mbox{for}\quad r_j\ge0,\quad
1\le j\le d,\quad\sum_{j=1}^dr_j=r,
\ey
converges uniformly on $\R^d$, whenever $r\in\Z_+$ and 
$r>s-d/p$. Note that for $q\le1$, the series \rf{ryad}
converges uniformly under the weaker assumption $r\ge s-d/p$.

We can define now the modified (homogeneous) Besov space $B^s_{p,q}\big(\R^d\big)$. We say that  $f\in B^s_{p,q}(\R^d)$ if \rf{Wn} holds and
$$
\frac{\partial^r f}{\partial x_1^{r_1}\cdots\partial x_d^{r_d}}
=\sum_{n\in\Z}\frac{\partial^r f_n}{\partial x_1^{r_1}\cdots\partial x_d^{r_d}}\quad
\mbox{for}\quad 
r_j\ge0,\quad\mbox
1\le j\le d,\quad\sum_{j=1}^dr_j=r,
$$
in the space ${\mathscr S}^\prime\big(\R^d\big)$, where $r$ is the minimal nonnegative number such that
$r>s-d/p$ ($r\ge s-d/p$  if $q\le1$). Now $f$ is determined uniquely by the sequence $\{f_n\}_{n\in\Z}$ modulo a polynomial of degree less than $r$. Also, a polynomial $g$ belongs to $B^s_{p,q}\big(\R^d\big)$ if and only if $\deg g<r$. 

In the case $p=q$ we use the notation $B_p^s(\R^d)$ for $B_{p,p}^s(\R^d)$.

Consider now the scale $\L_\a(\R^d)$, $\a>0$,  of {\it H\"older--Zygmund classes}. They can be defined by $\L_\a(\R^d)\df B_\be^\a(\R^d)$.

Besov classes admit many other descriptions. We give the one in terms of finite differences.
For $h$ in $\R^d$, we define the difference operator $\D_h$ by
$(\D_hf)(x)=\lb f(x+h)-f(x)$, $x\in\R^d$.

Let $s>0$,  $m\in\Z$, and $m-1\le s<m$.  
Suppose that $p,q\in[1,+\be]$.
The Besov class $B_{p,q}^s\big(\R^d\big)$ can be defined as the set of function
 $f$ in $L^1_{\rm loc}\big(\R^d\big)$ such that
$$
\int_{\R^d}|h|^{-d-sq}\|\D^m_h f\|_{L^p}^q\,dh<\be,~\; q<\be;\qquad
\sup_{h\not=0}\frac{\|\D^m_h f\|_{L^p}}{|h|^s}<\be,~\; q=\be.
$$
However, with this definition, Besov classes can contain polynomials of degree higher than in the case of the definition in terms of convolutions with the functions $W_n$.

The space $B_{pq}^s$ can be defined in terms of Poisson integral.
Let $P_d(x,t)$ be the Poisson kernel on
$\R^{d+1}_+\df\{(x,t): x\in\R^d, t>0\}$, i.e., $P_d(x,t)= c_dt(|x|^2+t^2)^{-\frac{d+1}2}$,
$c_d=\pi^{-\frac{d+1}2}\G(\frac{d+1}2)$. With each function   
$f$ in $L^1\big(\R^d,(\|x\|+1)^{-(d+1)}\,dx\big)$,
we can associate the Poisson integral $\mP f$,
$$
(\mP f )(x,t)=\int_{\R^d}P_d(x-y,t)f(y)\,dy.
$$
Then for every positive integer $m$, the following equality holds:
$$
\frac{\partial^m (\mP f)}{\partial^m t}(x,t)=\int_{\R^d}\frac{\partial^m P_d(x-y,t)}{\partial^m t}f(y)\,dy.
$$
Note that the second integral makes sense for all
$f\in L^1\big(\R^d,(\|x\|+1)^{-(d+m+1)}\,dx\big)$ which
allows us to define $\frac{\partial^m}{\partial t^m}\mP f$.

Let $m\in\Z$, \mbox{$m-1\le s<m$}, $1\le p,\:q\le+\be$. We can define 
$B_{pq}^s$ as the set of functions 
$f\in L^1\big(\R^d,(\|x\|+1)^{-(d+m+1)}\,dx\big)$
such that 
$$
\left(\int_0^{\be}t^{(m-s)q-1}\left\|\Big(\frac{\partial^m}{\partial t^m}\mP f\Big)(\cdot,t)\right\|_{L^p(\R^d)}^q\,dt\right)^{\frac1q}
<+\be,\quad q<+\be,
$$
$$
\sup_{t>0}t^{m-s}\left\|\Big(\frac{\partial^m}{\partial t^m}\mP f\Big)(\cdot,t)\right\|_{L^p(\R^d)}<+\be,\quad
q=+\be.
$$
It is also true that with this definition Besov classes can contain polynomials of degree  higher than in the case of the definition in terms of the convolutions with $W_n$.
Note also that this definition in terms of Poisson integral can also be used under certain provisions in the case when $p<1$  or $q<1$.

\medskip

We proceed now to {\it Besov classes of functions on the unit circle $\T$}. 
Let $w$ be a function satisfying \rf{w}.
We define the trigonometric polynomials $W_n$, $n\ge0$, by
$$
W_n(\z)\df\sum_{j\in\Z}w\left(\frac{|j|}{2^n}\right)\z^j,\quad n\ge1,
\quad W_0(\z)\df\sum_{\{j:\,|j|\le1\}}\z^j,\quad\z\in\T.
$$
If $f$ is a distribution on $\T$, put
$f_n=f*W_n$, $n\ge0$,
and say that $f$ belongs to the Besov class $B_{p,q}^s(\T)$, $s\in\R$, 
$0<p,\,q\le\be$, if
\bay
\label{Bperf}
\big\{2^{ns}\|f_n\|_{L^p}\big\}_{n\ge0}\in\ell^q.
\ey

Let $s\in\R$, $s>\max\{0,1/p-1\}$ and let $m$ be a positive integer such that $m>\max\{s,s+1/p-1\}$. Then a distribution $f$ on $\T$ belongs to
$B_{p,q}^s(\T)$  if and only if
$$
\int_0^{1}r(1-r^2)^{(m-s)q-1}\left\|\frac{\partial^m}{\partial r^m}\Big((\mP f)(r\z)\Big)\right\|_{L^p(\T)}^q\,dr
<+\be,\quad q<+\be,
$$
$$
\sup_{r\in[0,1)}(1-r^2)^{m-s}\left\|\frac{\partial^m}{\partial r^m}\Big((\mP f)(r\z)\Big)\right\|_{L^p(\T)}<+\be,\quad q=+\be,
$$
where $\mP f$  denotes the Poisson integral of the distribution $f$.

In the definitions of Besov classes in terms of Poisson integral we have considered the  
$m$th derivative in the variable $t$ in the first case and in the variable $r$ in the second case.
It is well known that in both cases we would get an equivalent definition if we required 
that similar expressions that involve all partial derivatives (including mixed ones) of order $m$ must be finite.

We refer the reader to \cite{Pee} and \cite{T} for more detailed information about Besov classes.

\medskip

{\bf2. Schatten--von Neumann classes.} For a bounded linear operator $T$ on Hilbert space, its {\it singular values} $s_j(T)$, $j\ge0$, are defined by
$$
s_j(T)\df\inf\big\{\|T-R\|:~\rank R\le j\big\}.
$$
The {\it Schatten--von Neumann class} $\bS_p$, $0<p<\be$, consists, by definition,
from operators $T$, for which
$$
\|T\|_{\bS_p}\df\Big(\sum_{j\ge0}\big(s_j(T)\big)^p\Big)^{1/p}<\be.
$$
For $p\ge1$,  this is a normed ideal of operators on Hilbert space. 
The class $\bS_1$ is called {\it trace class}. If $T$ is a trace class operator on a Hilbert space $\h$, its {\it trace} is defined by
$\trace T\df\sum_{j\ge0}(Te_j,e_j)$,
where $\{e_j\}_{j\ge0}$ is an orthonormal basis in $\h$. The right-hand side does not depend on the choice of a basis.

The class $\bS_2$ is called {\it Hilbert Schmidt class}. It forms a Hilbert space with inner product
$(T,R)_{\bS_2}\df\trace(TR^*)$.

For $p\in(1,\be)$, the dual space $(\bS_p)^*$ can be identified isometrically with the space $\bS_{p'}$, $1/p+1/p'=1$, with the help of the bilinear form
$\langle T,R\rangle\df\trace(TR)$.
The space dual to $\bS_1$ can be identified with the space of bounded linear operators with the help of the same bilinear form, while the space dual to the space of compact operators can be identified with $\bS_1$.

We refer the reader to \cite{GK} for more detailed information.

\medskip

{\bf3. Hankel operators.} For a function $\f$ of class $L^2$ on the unit circle 
$\T$, the {\it Hankel operator} $H_\f$ is defined on the dense subset of polynomials in the Hardy class $H^2$ by $H_\f f\df\pp_-\f f$, where $\pp_-$ is the orthogonal projection from $L^2$ onto $H^2_-\df L^2\ominus H^2$. By Nehari's theorem, $H_\f$ extends to a bounded operator from $H^2$ to $H^2_-$ if and only if there exists a function $\psi$ of class $L^\be$ on $\T$ whose Fourier coefficients $\hat\psi(n)$ satisfy the equality  
$\hat\psi(n)=\hat\f(n)$ for $n<0$. The last property, in turn, is equivalent, by Ch. Fefferman's theorem, to the fact that 
$\pp_-\f$ belongs to the class ${\rm BMO}$. 

A Hankel operator $H_\f$ belongs to the Schatten--von Neumann class $\bS_p$ if and only if the function $\pp_-\f$ belongs to the Besov class $B_p^{1/p}(\T)$. For 
$p\ge1$, this was proved in \cite{Pe0}, while for $p\in(0,1)$, in \cite{Pe<1}; see also \cite{Pek} and \cite{Se}, where other proofs are given for $p<1$.

It is easy to see that the operator $H_\f$ has matrix $\{\hat\f(-j-k-1)\}_{j\ge0,\,k>1}$
in the bases $\{z^j\}_{j\ge0}$ and $\{\ov z^k\}_{k\ge1}$. Such matrices, i.e., matrices of the form 
$\{\a_{j+k}\}_{j,k\ge0}$ are called {\it Hankel matrices}. The criterion for Hankel operators to belong to $\bS_p$ can be reformulated in the following way: {\it the operator on $\ell^2$ with Hankel matrix $\{\a_{j+k}\}_{j,k\ge0}$ belongs to $\bS_p$, $p>0$, if and only if the function $\sum_{j\ge0}\a_jz^j$ belongs to
$B_p^{1/p}(\T)$}.

We refer the reader to the monograph \cite{Pe} for proofs of the above results and for more detailed information on Hankel operators.

\medskip

{\bf4. Notation.}
We give here a list of some notation used in this survey.

\medskip

$\OL(\fF)$ is the space of operator Lipschitz functions on a closed subset $\fF$ of the complex plane $\C$;

$\CL(\fF)$ is the space of commutator Lipschitz functions on a closed subset $\fF$ of the complex plane $\C$;

$\OD(\R)$ is the space of operator differentiable functions on $\R$;

$\mB(\h_1,\h_2)$ is the space of bounded linear operators from a Hilbert space  $\h_1$ to a Hilbert space $\h_2$, $\mB(\h)\df\mB(\h,\h)$; 

$\mB_{\rm sa}(\h)$ is the space of bounded self-adjoint operators  on a Hilbert space 
$\h$;

$\m$ is normalized Lebesgue measure on the circle $\T$;

$\m_2$ is Lebesgue measure on the plane.

\

\begin{center}
\bf\Large Chapter I
\end{center}


\begin{center}
\bf\Large 
Operator Lipschitz functions on the line and on the circle. The first round
\end{center}
\label{1per}
\renewcommand{\thesection}{1.\arabic{section}}
\setcounter{section}{0}

\addtocontents{toc}{\vspace*{.33cm}\hspace*{-.45cm}\textbf{{\bf Chapter I}. Operator Lipschitz functions on the line and on the circle.\\ \hspace*{1.38cm}The first round}\hfill\pageref{1per}}

\

In this introductory chapter we consider operator Lipschitz functions on the real line  
$\R$ and on the unit circle $\T$. 
Later, in Chapter III, the class of operator Lipschitz functions will be subjected
to a more detailed study. We will also proceed to the study of operator Lipschitz functions on closed subsets of the complex plane $\C$.

We use the notation $\OL(\R)$ for the class of operator Lipschitz functions on $\R$ and for  $f\in\OL(\R)$ we put 
$$
\|f\|_{\OL(\R)}\df\sup\left\{\frac{\|f(A)-f(B)\|}{\|A-B\|}:~
A~\mbox{and}~B~\mbox{are self-adjoint operators},~A\ne B\right\}.
$$
In a similar way we introduce the space $\OL(\T)$ of operator Lipschitz functions on 
$\T$ be replacing self-adjoint operators with unitary operators.

It turns out that the class $\OL(\R)$ has somewhat unusual properties. In particular, functions of this class must be differentiable everywhere on $\R$ and also must have derivative at infinity, i.e., the limit
$\lim_{|t|\to\be}\frac{f(t)}{t}$
must exist (see Theorem \ref{pro} below). Note that this implies the result by McIntosh--Kato that has been mentioned in the introduction: the function $x\mapsto|x|$ is not operator Lipschitz. On the other hand, functions of class $\OL(\R)$ do not have to be continuously differentiable. In particular, the function $x\mapsto x^2\sin(1/x)$, not being continuously differentiable, is operator Lipschitz, see Theorem \ref{teorox2f1x} below. 

We begin this chapter with elementary examples of operator Lipschitz functions (see \S\:\ref{prim}).

We introduce in \S\:\ref{odif} the class of operator differentiable functions and the class of locally operator differentiable functions. It turns out that for the definition of these classes, it does not matter whether we consider differentiability in the sense of G\^ateaux or in the sense of Fr\'echet. We will see that (locally) operator differentiable functions must be continuously differentiable and that operator differentiable functions must be operator Lipschitz. However, not every operator Lipschitz function is operator differentiable. 

Besides operator Lipschitz functions, we consider in \S\:\ref{CoL} {\it commutator  Lipschitz functions}, i.e., functions $f$ on $\R$ such that
$$
\|f(A)R-Rf(B)\|\le\const\|AR-RB\|
$$
for arbitrary self-adjoint operators $A$ and $B$ (again, no matter bounded or not necessarily bounded) and for every bounded linear operator $R$. The
{\it commutator  Lipschitz norm}
$\|f\|_{\CL(\R)}$ of $f$ is defined as the minimal constant, for which the inequality holds. 
In a similar way we can define commutator  Lipschitz functions on the unit circle if instead of self-adjoint operators we consider unitary operators .

It turns out that for functions {\it on the line (as well as for functions on the circle) the class of commutator Lipschitz functions coincides with the class of operator Lipschitz functions}. In Chapter III we will see that for functions on an arbitrary closed subset of the plane $\R^2$ this is no longer true.

We obtain in this chapter a sufficient condition for operator Lipschitzness  on the line and on the circle (see \S\:\ref{Dost}) as well as a necessary condition (see \S\:\ref{Neob}) and compare them with each other.

It would be also natural to consider the class of {\it operateor H\"older functions of order} $\a$, $0<\a<1$, i.e., the class of functions $f$ such that
$$
\|f(A)-f(B)\|\le\const\|A-B\|^\a,
$$
for self-adjoint operators $A$ and $B$ on Hilbert space. However, this term turns out to be short-lived because {\it every function $f$ on $\R$ of H\"older class of order $\a$ is necessarily operator H\"older or order $\a$,} see \S\:\ref{oHold}.

\medskip

\section{\bf Elementary examples of operator Lipschitz functions}
\setcounter{equation}{0}
\label{prim}

\medskip

In this section we give examples of operator Lipschitz functions on the line and on the circle and obtain simple sufficient conditions for operator Lipschitzness.

\medskip

{\bf Example 1.}  {\it For every $\l$ in $\C\setminus\R$, the function $(\l-x)^{-1}$ is operator Lipschitz on $\R$
and $\|(\l-x)^{-1}\|_{\OL(\R)}=|\im\l|^{-2}$}.

\medskip

\Pf The Hilbert resolvent identity
$$
(\l I-A)^{-1}-(\l I-B)^{-1}=(\l I-A)^{-1}(A-B)(\l I-B)^{-1}
$$
immediately implies that $\|(\l-x)^{-1}\|_{\OL(\R)}\le|\im\l|^{-2}$. It remains to observe that
\lb$\|(\l-x)^{-1}\|_{\OL(\R)}\ge\|(\l-x)^{-1}\|_{\Li(\R)}=|\im\l|^{-2}$. $\bl$

\medskip

{\bf Example 1$'$.}  {\it For every $\l$ in $\C\setminus\T$, the function $(\l-z)^{-1}$ is operator Lipschitz on $\T$
and $\|(\l-z)^{-1}\|_{\OL(\T)}=(|\l|-1)^{-2}$}.

\medskip

{\bf Example 2.}  {\it The function $x\mapsto\log(1+{\rm i}x)$ is operator Lipschitz on $\R$ and
\lb$\|\log(1+{\rm i}x)\|_{\OL(\R)}=1$. Here $\log$ means the principal branch of logarithm.}

\medskip

\Pf
Clearly, 
$$
\log(1+{\rm i}x)=\int_0^{+\be}\left(\frac1{1+t}-\frac1{1+t+{\rm i}x}\right)\,dt.
$$
It follows that
\begin{align*}
\|\log(1+{\rm i}x)\|_{\OL(\R)}
&\le\int_0^{+\be}\left\|\frac1{1+t}-\frac1{1+t+{\rm i}x}\right\|_{\OL(\R)}dt\\[.2cm]
&=\int_0^{+\be}\left\|\frac1{1+t+{\rm i}x}\right\|_{(\OL)(\R)}dt
=\int_0^{+\be}\frac{dt}{(1+t)^2}=1.
\end{align*}
On the other hand, the inequality $\|\log(1+{\rm i}x)\|_{\OL(\R)}\ge1$ is obvious because 
$$
\|\log(1+{\rm i}x)\|_{\OL(\R)}\ge\|\log(1+{\rm i}x)\|_{\Li(\R)}=1.\quad\bl
$$

\medskip

In a similar way one can prove that for every $\l$ in $\C\setminus\R$, we have the equality \lb$\|\log(\l-x)\|_{\OL(\R)}=|\im\l|^{-1}$, where $\log(\l-x)$ denotes any
regular branch of the function $\log(\l-z)$ on $\R$.

\medskip

{\bf Example 3.}  {\it The function $\arctan$ is operator Lipschitz and $\|\arctan\|_{\OL(\R)}=1$}.

\medskip

\Pf It suffices to verify that $\|\arctan\|_{\OL(\R)}\le1$. To this end, we observe that $\arctan x=\im\log(1+{\rm i}x)$, $x\in\R$. $\bl$

\medskip

{\bf Example 4.} {\it For every positive integer $n$, the following equality holds:
$$
\|(\l -x)^{-n}\|_{\OL(\R)}=n|\im\l|^{-n-1}\quad\mbox{for every}\quad\l\in\C\setminus\R.
$$}

\medskip

\Pf Substituting $X=(\l I-A)^{-1}$ and $Y=(\l I-B)^{-1}$, in the elementary identity
\bay
\label{xny}
X^n-Y^n=\sum_{k=1}^nX^{n-k}(X-Y)Y^{k-1},
\ey
we obtain
$$
(\l I-A)^{-n}-(\l I-B)^{-n}=\sum_{k=1}^n(\l I-A)^{k-n}\big((\l I-A)^{-1}-(\l I-B)^{-1}\big)(\l I-B)^{1-k}.
$$
Therefore, for arbitrary self-adjoint operators $A$ and $B$
we have
\begin{align*}
\|(\l I-A)^{-n}&-(\l I-B)^{-n}\|\\[.2cm]
&\le\sum_{k=1}^n\|(\l I-A)^{k-n}\|\cdot\big\|(\l I-A)^{-1}-(\l I-B)^{-1}
\big\|\cdot\|(\l I-B)^{1-k}\|\\[.2cm]
&\le\sum_{k=1}^n|\im\l|^{k-n}|\im\l|^{-2}\|A-B\|\cdot|\im\l|^{1-k}=n|\im\l|^{-n-1}\|A-B\|.
\end{align*}
Thus, we have proved that
$\|(\l I-x)^{-n}\|_{\OL(\R)}\le n|\im\l|^{-n-1}$.
It remains to observe that $\|(\l -x)^{-n}\|_{\OL(\R)}\ge\|(\l -x)^{-n}\|_{\Li(\R)}=n|\im\l|^{-n-1}$. $\bl$

\medskip

{\bf Example 5.} {\it The function 
$x\mapsto e^{{\rm i}ax}$, $a\in\R$, is operator Lipschitz and 
$\|e^{{\rm i}ax}\|_{\OL(\R)}=|a|$}.

\medskip
 
\Pf Again, it suffices to establish only the upper estimate. We may
assume that $a=1$. Let $A$ and $B$ be self-adjoint operators. Then
$$
\Big(e^{{\rm i}tA}e^{-{\rm i}tB}\Big)'={\rm i}Ae^{{\rm i}tA}e^{-{\rm i}tB}-{\rm i}e^{{\rm i}tA}e^{-{\rm i}tB}B
={\rm i}e^{{\rm i}tA}(A-B)e^{-{\rm i}tB},
$$
whence
\begin{align*}
\|e^{{\rm i}A}-e^{{\rm i}B}\|&=\|e^{{\rm i}A}e^{-{\rm i}B}-I\|
=\left\|{\rm i}\int_0^1e^{{\rm i}tA}(A-B)e^{-{\rm i}tB}\,dt\right\|\\[.2cm]
&\le\int_0^1\big\|e^{{\rm i}tA}(A-B)e^{-{\rm i}tB}\big\|\,dt
=\int_0^1\|A-B\|\,dt=\|A-B\|.\quad\bl
\end{align*}

In all the above examples we have the equality
$\|f\|_{\OL(\R)}=\|f'\|_{L^\be(\R)}$ which is rather an exception than a rule.

Example 5 immediately implies the following fact:

\begin{thm}
\label{char}
Let $f$ be a primitive of the Fourier transform $\F\mu$
of a complex Borel measure $\mu$ on $\R$. Then $f\in\OL(\R)$
and $\|f\|_{\OL(\R)}\le\|\mu\|$.
\end{thm}

\Pf We may assume that $f(0)=0$. Then
\begin{align*}
f(x)&=\int_0^x(\F\mu)(t)\,dt=
\int_0^x\left(\int_\R e^{{-\rm i}st}\,d\mu(s)\right)\,dt\\[.2cm]
&=\int_0^1\left(\int_\R xe^{{-\rm i}stx}\,d\mu(s)\right)\,dt
={\rm i}\int_\R\frac{e^{{-\rm i}sx}-1}{s}\,d\mu(s).
\end{align*}
Hence,
$$
\|f\|_{\OL(\R)}\le\int_\R\left\|\frac{e^{{-\rm i}sx}-1}{s}\right\|_{\OL(\R)}\,d|\mu|(s)
\le\int_\R d|\mu|(s)=\|\mu\|.\quad\bl
$$

\begin{cor}
\label{corchar}
Let $f\in C^1(\R)$. Suppose that the function $f'$ is positive definite.
Then $\|f\|_{\OL(\R)}=\|f\|_{\Li(\R)}=f'(0)$.
\end{cor}

\Pf By the classical theorem of Bochner, see, for example, \cite{I}, the positive definite function $f'$ can be represented in the form  $f'=\F\mu$, where $\mu$ is a finite Borel positive measure on $\R$. It remains to observe that
$$
\|\mu\|=f'(0)=|f'(0)|\le\|f\|_{\Li(\R)}\le\|f\|_{\OL(\R)}\le \|\mu\|,
$$
where the last inequality follows from Theorem \ref{char}. $\bl$

\medskip

In this section practically all the above examples of an explicit evaluation 
of the seminorm in $\OL(\R)$ are based more or less on Corollary \ref{corchar}.
Nevertheless, one can construct an example of a function $f$ in $\OL(\R)$ such that $\|f\|_{\OL(\R)}=\|f\|_{\Li(\R)}=f'(0)$ and $f$ does not satisfy the hypotheses of Corollary  \ref{corchar}.

On the other hand, if $\|f\|_{\OL(\R)}=\|f\|_{\Li(\R)}=\frac{f(a)-f(0)}{a}=1$ for  $a\in\R$, $a\ne0$,
then $f(x)=x+f(0)$ for every $x\in\R$.

Example 5 admits one more generalization, the so-called operator Bernstein inequality. This will be discussed in \S\:\ref{Bern}.
In particular, it will be shown in \S\:\ref{Bern} that $L^\be(\R)\cap \E_\s\subset\OL(\R)$, where the symbol $\E_\s$ denotes the space of entire functions of exponential type at most $\s$. 

Consider now examples of operator Lipschitz functions on the unit circle $\T$.

\medskip

{\bf Example 6.}  {\it Let $n\in\Z$.
Then $\|z^n\|_{\OL(\T)}=|n|$ for all $n\in\Z$}.

\medskip

\Pf  It suffices to consider the case when $n>0$; then everything reduces to the verification of the inequality:
$
\|U^n-V^n\|\le n\|U-V\|
$
for arbitrary unitary operators $U$ and $V$.  To prove it,
it suffices to substitute  $X=U$ and $Y=V$ in \rf{xny}. $\bl$

\medskip

This example immediately leads to an analog of Theorem \ref{char} for the circle.

\begin{thm} 
Let $f$ be a continuous function on the unit circle $\T$
such that \lb$\sum\limits_{n\in\Z}|n|\cdot|\hat f(n)|<\be$. Then $f\in\OL(\T)$
and $\|f\|_{\OL(\T)}\le\sum\limits_{n\in\Z}|n|\cdot|\hat f(n)|$.
\end{thm}

Note that stronger results will be given soon in \S\:\ref{Dost}.

\medskip 

{\bf Example 7.} The function $x\mapsto x^2\sin\frac1x$ is operator Lipschitz.
To get convinced in this, we prove the following theorem:

\begin{thm}
\label{teorox2f1x}
let $f\in\OL(\R)$  and $f(0)=0$. 
Put $g(x)=x^{2}f(x^{-1})$ for $x\ne0$ and $g(0)=0$.
Then $g\in\OL(\R)$ and
\bay
\label{nervox21x}
\frac13\|f\|_{\OL(\R)}\le\|g\|_{\OL(\R)}\le3\|f\|_{\OL(\R)}.
\ey
\end{thm}

\medskip

\Pf It suffices to prove the second inequality because the first inequality reduces to the second one. We may assume that $\|f\|_{\OL(\R)}=1$. 
As we have mentioned in the introduction to this chapter, for functions on the line, operator Lipschitzness is equivalent to commutator Lipschitzness and the corresponding norms coincide: 
$\|f\|_{\OL(\R)}=\|f\|_{\CL(\R)}$ (see \S\:\ref{CoL}
and \S\:\ref{oplip}). Therefore, it suffices to prove that the inequality
\bay
\label{71}
\|f(A)R-Rf(A)\|\le\|AR-RA\|
\ey
for every bounded operator $R$ and every bounded self-adjoint operator $A$
implies that
$$
\|g(A)R-Rg(A)\|\le3\|AR-RA\|
$$
for every bounded operator $R$ and every self-adjoint operator $A$.
Suppose that, $A$ is invertible. This case reduces to the following assertion:
\bay
\label{72}
\|A^2f(A^{-1})R-RA^2f(A^{-1})\|\le3\|AR-RA\|
\ey
for every bounded operator $R$ and every invertible self-adjoint operator $A$.
We have
\begin{align*}
f(A^{-1})A^2R&-RA^2f(A^{-1})=f(A^{-1})A(AR-RA)\\[.2cm]
&+f(A^{-1})ARA-ARAf(A^{-1})
+(AR-RA)Af(A^{-1}). 
\end{align*}
Clearly,
$$
\|Af(A^{-1})\|\le\sup_{t\ne0}|t^{-1} f(t)|\le\|f\|_{\Li(\R)}\le\|f\|_{\OL(\R)}=1.
$$
Consequently, 
$$
\|f(A^{-1})A(AR-RA)\|\le\|AR-RA\|\quad\mbox{and}\quad\|(AR-RA)Af(A^{-1})\|\le\|AR-RA\|.
$$
Substituting in \rf{71} the operators $ARA$ and $A^{-1}$, we obtain
$$
\|f(A^{-1})ARA-ARAf(A^{-1})\|\le\|A^{-1}ARA-ARAA^{-1}\|=\|AR-RA\|
$$
which immediately implies \rf{72}.
To consider the general case, it is sufficient to observe that for every positive number $\d$, there exists an invertible self-adjoint operator $A_\d$ such that
$AA_\d=A_\d A$ and $\|A-A_\d\|<\d$. Then for all $\d>0$,
\begin{align*}
\|g(A)R&-Rg(A)\|\le\|g(A)-g(A_\d)\|\cdot\|R\|+\|g(A_\d)R-Rg(A_\d)\|\\[.2cm]
&+\|g(A_\d)-g(A)\|\cdot\|R\|
\le2\d\|R\|\cdot\|g\|_{\Li(\R)}+3\|A_\d R-RA_\d\|\\[.2cm]
&\le6\d\|R\|\cdot\|f\|_{\Li(\R)}+3\|AR-RA\|+6\d\|R\|
\le3\|AR-RA\|+12\d\|R\|.\quad\bl
\end{align*}

\medskip

{\bf Remark.} Now it is clear that in view of Example 5, the function $g$ defined by
$g(x)=x^2\sin\frac1x$, is operator Lipschitz.  {\it The function $g$ gives an example of an operator Lipschitz function that is not continuously differntiable}. 
The problem of the existence of such functions was posed in \cite{Wi} and was solved in \cite{KS1}. The fact that $g$ is operator Lipschitz on every finite interval was established in \cite{KS2}.
Recall (see Theorem \ref{pro} below) that every operator Lipschitz function on $\R$ 
must be differentiable everywhere.

Note also that it was proved in \cite{A2} that a set on the line is the set of discontinuities of the derivative of an operator Lipschitz function if and only if it is an $F_\s$ set of the first category. In other words, the sets of discontinuity points of the derivatives of operator Lipschitz functions are the same as the sets of discontinuity
points of functions of the first Baire class.

\medskip

In \S\:\ref{dlp} Theorem \ref{teorox2f1x} will be generalized to the case of arbitrary linear-fractional transformations.

\medskip

\section{\bf Operator Lipschitzness in comparison with operator differentiability}
\setcounter{equation}{0}
\label{odif}

\medskip

Let $H$ be a function on a subset $\L$ of the real line $\R$ with values in a Banach space $X$. It is called {\it Lipschitz} if there is a nonnegative number $c$ such that
\bay
\label{Xlip+}
\|H(s)-H(t)\|_X\le c|s-t|,\quad s,\:t\in\L.
\ey
We denote the set of all such functions by $\Li(\L,X)$. Let $\|H\|_{\Li(\L,X)}$ denote the least constant $c$ satisfying \rf{Xlip+}. 
As usual, we put $\|H\|_{\Li(\L,X)}\df\be$ if $H\not\in\Li(\L,X)$.
 
Let $f$ be a continuous function on $\R$.
With each self-adjoint operator $A$ and each bounded self-adjoint operator
$K$, we associate the function $H_{f,A,K}$, $H_{f,A,K}(t)=\lb f(A+tK)-f(A)$,
defined for those $t$ in $\R$, for which $f(A+tK)-f(A)\in\mB(\h)$.

Note that if $f\in\OL(\R)$, then $H_{f,A,K}\in\Li(\R,\mB(\h))$ and
$\|H_{f,A,K}\|_{\Li(\R,\mB(\h))}\le\|K\|\cdot\|f\|_{\OL(\R)}$.
It is easy to see that the following result holds.

\begin{lem} 
\label{HfAK+}
Let $f$ be a continuous function on $\R$. Then
\begin{align*}
\|f\|_{\OL(\R)}&=\sup\big\{\|H_{f,A,K}\|_{\Li(\R,\mB(\h))}:~A,~K\in\mB_{\rm sa}(\h),\|K\|=1\big\}\nonumber\\[.2cm]
&=\sup\big\{\|H_{f,A,K}\|_{\Li(\R,\mB(\h))}:~K\in\mB_{\rm sa}(\h),~
\|K\|=1,~A^*=A \big\}.\quad\bl
\end{align*} 
\end{lem}

We need the following well-known elementary fact.
For the sake of reader's convenience, we give one of its existing proofs.
 
\begin{lem} 
\label{verh+}
Let $H$ be a function with values in a Banach space $X$ that is defined on a nondegenerate interval $\L$, $\L\subset\R$. Then
$$
\|H\|_{\Li(\L,X)}=\sup_{t\in\L}\varlimsup_{h\to0}\frac{\|H(t+h)-H(t)\|_X}{|h|}.
$$
\end{lem}

\Pf The inequality $\ge$ is evident.
To prove the opposite inequality, it suffices to show that inequality \rf{Xlip+} holds whenever $c$ satisfies the condition
\bay
\label{cXlip+}
c>\sup_{t\in\L}\varlimsup_{h\to0}\frac{\|H(t+h)-H(t)\|_X}{|h|}.
\ey
We fix such a number $c$ and an arbitrary point $t$ in $\L$.  Let $\L_t$ be the set
of points $s$ in $\L$ that satisfy inequality \rf{Xlip+}. It follows immediately from \rf{cXlip+} that the set $\L_t$
 is at the same time open and closed in $\L$. 
Moreover, $\L_t\ne\varnothing$
for $t\in\L$. Consequently, $\L_t=\L$ because $\L$ is connected. $\bl$

\begin{thm} 
\label{lippoint+}
Let $f$ be a continuous function on $\R$.
Suppose that
$$
\varlimsup_{t\to0}\frac{\|f(A+tK)-f(A)\|}{|t|}<+\be
$$
for every (not necessarily bounded) self-adjoint operator $A$ and for every bounded self-adjoint operator $K$. Then $f\in\OL(\R)$.
\end{thm}

\Pf It follows easily from Lemmata \ref{HfAK+} and \ref{verh+} that
$$
\|f\|_{\OL(\R)}=\sup\left\{\varlimsup_{t\to0}\frac{\|f(A+tK)-f(A)\|}{|t|}:~A,~K\in\mB_{\rm sa}(\h),\|K\|=1\right\}
$$
Thus, if we assume that $\|f\|_{\OL(\R)}=\be$, then for each 
$n$ in $\Z_+$, there exist operators 
$A_n,\:K_n\in\mB_{\rm sa}(\h)$ such that $\|K_n\|=1$ and 
$$
\varlimsup_{t\to0}\frac{\|f(A_n+tK_n)-f(A_n)\|}{|t|}>n.
$$
Consider the self-adjoint operators $\A$ and ${\mathcal K}$ on the Hilbert space 
$\ell^2(\h)$ defined by
\bay
\label{oprA+}
\A(x_0,x_1,x_2,\cdots)=(A_0x_0,A_1x_1,A_2x_2,\cdots),\quad
(x_0,x_1,x_2,\cdots)\in\ell^2(\h)
\ey
and
\bay
\label{Knal2}
{\mathcal K}(x_0,x_1,x_2,\cdots)=(K_0x_0,K_1x_1,K_2x_2,\cdots),\quad
(x_0,x_1,x_2,\cdots)\in\ell^2(\h).
\ey
Then 
$$
\varlimsup_{t\to0}\frac{\|f(\A+t\mathcal K)-f(\A)\|}{|t|}\ge
\varlimsup_{t\to0}\frac{\|f(A_n+tK_n)-f(A_n)\|}{|t|}>n
$$
for every nonnegative integer $n$ and we arrive at a contradiction. $\bl$

\medskip

{\bf Remark.} It can be seen from the proof of Theorem \ref{lippoint+} that the following equalities hold:
\begin{align*}
\|f\|_{\OL(\R)}&=\sup
\left\{\varlimsup_{t\to0}\frac{\|f(A+tK)-f(A)\|}{|t|}:~
A,~K\in\mB_{\rm sa}(\h),~\|K\|_{\mB_{\rm sa}(\h)}=1\right\}\\[.2cm]
&=\sup\big\{\|H_{f,A,K}\|_{\Li(\R)}:
~A,~K\in\mB_{\rm sa}(\h), \|K\|_{\mB_{\rm sa}(\h)}=1\big\}.
\end{align*}

\medskip
 
To state the next theorem, we observe that the function $H_{f,A,K}$ is differentiable 
for arbitrary self-adjoint operators $A$ and $K$ if and only if
it is differentiable at $\0$ for arbitrary self-adjoint operators $A$ and $K$ (as usual, the operator $K$ is assumed to be bounded).

The proof of the following theorem uses Theorem \ref{sildif}, which will be proved in Chapter III. 

\begin{thm} 
\label{4i+}
Let $f$ be a continuous function on $\R$. Then the condition

{\rm (a)} $f\in\OL(\R)$

\noindent
is equivalent to each of the following statements for every self-adjoint operator $A$ and every bounded self-adjoint operator $K$:

{\rm (b)} $H_{f,A,K}\in\Li(\R,\mB(\h))$ ;

{\rm (c)} the function $H_{f,A,K}$ is differentiable as a function from $\R$
to the space $\mB(\h)$ equipped with the weak operator topology;

{\rm (d)}  the function $H_{f,A,K}$ is differentiable as a function from $\R$
to the space $\mB(\h)$ equipped with the strong operator topology.
\end{thm}

\Pf The implications (a)$\Longrightarrow$(b) and (d)$\Longrightarrow$(c) are obvious.
The implication (a)$\Longrightarrow$(d) follows from Theorem \ref{sildif} below.
Finally, the implications (c)$\Longrightarrow$(a) and (b)$\Longrightarrow$(a) 
follow immediately from Theorem \ref{lippoint+}. $\bl$

\medskip

We denote by $\OL_{\rm loc}(\R)$ the space of continuous functions $f$ on $\R$ such that $f\big|[-a,a]\in\OL([-a,a])$ for every $a>0$ and by
$\Li_{\rm loc}(\R,\mB(\h))$ the space of continuous functions $f$ on $\R$
such that $f\big|[-a,a]\in\Li([-a,a],\mB(\h))$ for every $a>0$.
All the results of this section also have natural analogs for these spaces.

\begin{thm} 
\label{loclippoint+}
Let $f$ be a continuous function on $\R$.
Suppose that
$$
\varlimsup_{t\to0}\frac{\|f(A+tK)-f(A)\|}{|t|}<\be
$$
for all $A,~K\in\mB_{\rm sa}(\h)$. Then   $f\in\OL_{\rm loc}(\R)$.
\end{thm}

\Pf Suppose that $f\not\in\OL_{\rm loc}(\R)$. Then $f\not\in\OL([-a,a])$
for some $a>0$. Thus, for each $c\ge0$,
there exist operators $A$ and $K$ in $\mB_{\rm sa}(\h)$ such that $\|A\|\le a$, $\|A+K\|\le a$
and $\|f(A+K)-f(A)\|>c\|K\|$. Repeating the reasoning of the proof of Theorem \ref{lippoint+},
we arrive at a contradiction by constructing self-adjoint operators $\A$ and $\A+\mathcal K$  such that
$$
\|\A\|\le a,\quad\|\A+\mathcal K\|\le a\quad\mbox{and}\quad
\varlimsup_{t\to0}\frac{\|f(\A+t\mathcal K)-f(\A)\|}{|t|}=\be.\quad\bl
$$

\begin{thm} 
\label{loc4i+}
Let $f$ be a continuous function on $\R$. The following statements are equivalent:

{\rm (a)} $f\in\OL_{\rm loc}(\R)$;

{\rm (b)} $H_{f,A,K}\in\Li_{\rm loc}(\R,\mB(\h))$ for all $A,~K\in\mB_{\rm sa}(\h)$;

{\rm (c)} for all $A,\:K$ in $\mB_{\rm sa}(\h)$, the function $H_{f,A,K}$ is differentiable as a function from $\R$
to the space $\mB(\h)$ equipped with the weak operator topology;

{\rm (d)} for all $A,\:K$ in $\mB_{\rm sa}(\h)$ the function $H_{f,A,K}$ is differntiable as a function from $\R$
to the space $\mB(\h)$ equipped with the strong operator topology.
\end{thm}

This theorem can be proved by analogy with Theorem \ref{4i+};
however, instead of Theorem \ref{lippoint+} one has to use Theorem \ref{loclippoint+}.

Note that in \cite{KS} it was shown that (a) in Theorem 
\ref{loc4i+} is equivalent to differentiability in the norm in all compact directions for all bounded self-adjoint operators.

It follows from Theorem \ref{4i+} that if $f$ is a continuous function on $\R$, then  
$f\in\OL(\R)$ if and only if  for every self-adjoint operator $A$ and every bounded self-adjoint operator $K$ the limit
\bay
\label{silndif+}
\lim_{t\to0}\frac1t(f(A+tK)-f(A))\df \bs{d}_{f,A}K
\ey
exists in the strong operator topology.
It will also follow from Theorem \ref{sildif} given below that $\bs{d}_{f,A}$
is a bounded linear operator from $\mB_{\rm sa}(\h)$ to $\mB(\h)$.

Similar results also hold for functions $f$ in $\OL_{\rm loc}(\R)$
with the only difference that the operator $A$ must be bounded.

It follows from Theorem \ref{loc4i+} that if $f$ is a continuous function on $\R$, then 
$f\in\OL_{\rm loc}(\R)$ if and only if  
for arbitrary operators
$A$ and $K$ in $\mB_{\rm sa}(\h)$, the limit in \rf{silndif+} exists in the strong operator topology; herewith  $\bs{d}_{f,A}$ is a bounded linear operator from 
$\mB_{\rm sa}(\h)$ to $\mB(\h)$.

\begin{thm} 
Let $f\in\OL_{\rm loc}(\R)$. Then 
\begin{align*}
\|f\|_{\OL(\R)}&=\sup_{A\in\mB_{\rm sa}(\h)}\|\bs{d}_{f,A}\|\nonumber\\[.2cm]
&=\sup\big\{\|\bs{d}_{f,A}\|:~A~\mbox{ is a self-adjoint operator}\big\}.
\end{align*} 
\end{thm}

As usual, the equality $\|f\|_{\OL(\R)}=\be$ means that  
$f\not\in\OL(\R)$.

\medskip

\Pf  It suffices to use Lemma \ref{HfAK+}. $\bl$

\begin{thm} 
\label{operdif+}
Let $f$ be a continuous function on $\R$.
Suppose that for every self-adjoint operator $A$ and for every bounded self-adjoint operator $K$,
the limit in {\em\rf{silndif+}} exists in the operator norm.
Then $f\in\OL(\R)\cap C^1(\R)$,
the map $K\mapsto f(A+K)-f(A)$ ($K\in\mB_{\rm sa}(\h)$) is differentiable in the sense of Fr\'echet at $\0$ for every self-adjoint operator $A$
and its differential at $\0$ is equal to $\bs{d}_{f,A}$.
\end{thm}

\Pf The membership $f\in\OL(\R)$ follows from Theorem \ref{4i+}.
It follows from Theorem \ref{sildif}
that $\bs{d}_{f,A}$ is a bounded linear operator from $\mB_{\rm sa}(\h)$ to $\mB(\h)$.
Let us verify {\it differentiability in the sense of Fr\'echet}, i.e., that $\bs{d}_{f,A}$ is a bounded linear operator (already proved) and
$$
\lim_{t\to0}\frac1t\|f(A+tK)-f(A)-t\bs{d}_{f,A}K\|=0
$$
uniformly in $K$ in the unit sphere of $\mB_{\rm sa}(\h)$.  

It suffices to see that 
\bay
\label{lodF+}
\lim_{n\to\be}\frac1{t_n}\|f(A+t_nK_n)-f(A)-t_n\bs{d}_{f,A}K_n\|=0
\ey
for an arbitrary sequence $\{t_n\}_{n\ge0}$ of nonzero real numbers that tends to zero and for an arbitrary sequence of self-adjoint operators  $\{K_n\}_{n\ge0}$ such that
$\|K_n\|=1$ for every $n$.

Consider the self-adjoint operator $\A$ and the bounded self-adjoint operator ${\mathcal K}$ on $\ell^2(\h)$ defined by \rf{oprA+} with $A_n=A$ and \rf{Knal2}.
Applying the assumptions of the theorem to the operators $\A$ and $\mathcal K$,  we obtain
\bay
\label{lodG+}
\lim_{n\to\be}\frac1{t_n}\|f(\A+t_n{\mathcal K})-f(\A)-t_n\bs{d}_{f,\A}{\mathcal K}\|=0.
\ey
Clearly, $\bs{d}_{f,\A}{\mathcal K}$ is the orthogonal sum of the operators 
$\bs{d}_{f,A}K_n$, $n\ge0$, and so \rf{lodF+} is a consequence of \rf{lodG+}. 

Finally, let us prove that $f\in C^1(\R)$.
We are going to verify the continuity of the derivative $f'$ at an arbitrary point $t_0$.
Let $A$ be the operator of multiplication by $x$ on $L^2([x_0-1,x_0+1])$.
Put $K\df I$.
Then by the assumptions of the theorem, the limit  $\lim_{t\to0} t^{-1}(f(A+tI)-f(A))$ exists in the norm.
Thus, the limit $\lim_{t\to0} t^{-1}(f(x+t)-f(x))=f'(x)$ exists in 
$L^\be([x_0-1,x_0+1])$,
whence $f\in C^1(t_0-1,t_0+1)$.
$\bl$

\medskip

{\bf Definition.}
A function  $f$ satisfying the hypotheses of Theorem \ref{operdif+} is called
{\it operator differentiable}. We denote by $\OD(\R)$ the set of all operator differentiable functions on $\R$.

\medskip 

Recall that for functions defined on Banach spaces, there are different notions of differentiability: the existence of a weak derivative in the sense of G\^ateaux; the existence of a G\^ateaux differential; differentiability in the sense of Fr\'echet. However, as one can see from Theorem \ref{operdif+}, in the case of operator differentiability of functions on the line, all these definitions are equivalent. Note that the equivalence of operator differentiability  in the sense of Fr\'echet and the existence of a G\^ateaux differential that is a bounded linear operator is proved in \cite{KS}. 

The following result can be proved in about the same way as Theorem \ref{operdif+}.

\begin{thm} 
\label{locoperdif+}
Let $f$ be a continuous function on $\R$.
Suppose that for arbitrary $A$ and $K$ in $\mB_{\rm sa}$, 
the limit in  {\em\rf{silndif+}} exists in the operator norm.
Then $f\in\OL_{\rm loc}(\R)\cap C^1(\R)$,
the map $K\mapsto f(A+K)-f(A)$, $K\in\mB_{\rm sa}$, is differentiable 
in the sense of Fr\'echet at $\0$ for every $A$ in $\mB_{\rm sa}$
and its differential at $\0$ is equal to $\bs{d}_{f,A}$.
\end{thm}

If a function $f$ satisfies the hypotheses of Theorem \ref{locoperdif+}, we say that it is {\it locally operator differentiable} and write 
$f\in{\rm OD}_{\rm loc}(\R)$. 

Observe that Theorems \ref{operdif+} and \ref{locoperdif+} affirm, in particular, that 
{\it if $f\in{\rm OD}_{\rm loc}(\R)$, then $f$ is continuously differentiable and belongs to the class $\OL_{\rm loc}(\R)$ and if $f\in{\rm OD}(\R)$, then
$f\in\OL(\R)$}.

\medskip

{\bf Remark.} Note that the function $g$, $g(x)=x^2\sin\frac1x$, not being continuously differentiable, cannot be operator differentiable. Thus, it is impossible to replace in Theorem \ref{teorox2f1x} the class of operator Lipschitz functions with the class of operator differentiable functions.
Indeed, it is easy to verify that the function  $x\mapsto\sin x=\im e^{{\rm i}x}$ is operator differentiable.

\medskip

Our immediate purpose is to prove the continuous dependence (in the operator norm) of the differential $\bs{d}_{f,A}$ on the operator $A$ for (localy) operator differentiable function $f$. The following result was obtained in \cite{KS}.

\begin{thm}
\label{ndlod+}
Let $f$ be a locally operator differentiable function and let $c>0$. Then
for every $\e>0$, there exists $\d>0$ such that
$
\|\bs{d}_{f,A}-\bs{d}_{f,B}\|\le\e,
$
whenever
$A$ and $B$ are self-adjoint operators such that
$\|A\|\le c$, $\|B\|\le c$ and $\|A-B\|\le\d$.
\end{thm}

First, we prove the following lemma obtained in \cite{KS}.

\begin{lem}
\label{01+}
Let $f$ be a locally operator differentiable function.  
Then for arbitrary positive numbers $c$ and $\e$ 
there exists $\d>0$ such that 
$$
\|f(A+K)-f(A)-\bs{d}_{f,A}K\|\le\e\|K\|,
$$
whenever $A$ and $K$ are self-adjoint operators such that $\|K\|\le\d$ and $\|A\|\le c$.
\end{lem}

\Pf Assume the contrary. Then for some positive numbers $c$ and $\e$, there are
sequences of self-adjoint operators $\{A_n\}_{n\ge0}$ and $\{K_n\}_{n\ge0}$
such that $\|K_n\|\to0$, $\|A_n\|\le c$ and
\bay
\label{>e+}
\|f(A_n+K_n)-f(A_n)-\bs{d}_{f,A_n}K_n\|>\e\|K_n\|,\quad n\ge0.
\ey
Let $\A$ be the bounded self-adjoint operator on $\ell^2(\h)$ defined by \rf{oprA+}.
Then $\|\A\|\le c$. Since $f$ is differentiable in the sense of Fr\'echet at the point $A$, there exists $\d>0$ such that
\bay
\label{A_Ke+}
\|f(\A+K)-f(\A)-\bs{d}_{f,\A}K\|\le\e\|K\|
\ey
for every self-adjoint operator $K$ satisfying the condition $\|K\|\le\d$. 
Let us define now the operator ${\mathcal K}_n$ on $\ell^2(\h)$ by
\bay
\label{oprKn+}
{\mathcal K}_n(x_0,x_1,x_2,\cdots)=(\0,\,\cdots,\0,K_nx_n,\0,\0,\cdots),\quad
 (x_0,x_1,x_2,\cdots)\in\ell^2(\h).
\ey

Applying inequality \rf{A_Ke+} for sufficiently large $n$, we obtain
$$
\|f(A_n+K_n)-f(A_n)-\bs{d}_{f,A_n}K_n\|=
\|f(\A+{\mathcal K}_n)-f(\A)-\bs{d}_{f,\A}{\mathcal K}_n\|\le\e\|{\mathcal K}_n\|
=\e\|K_n\|
$$
which contradicts inequality \rf{>e+}. $\bl$

\medskip

{\bf The proof of Theorem \ref{ndlod+}.} Let $c$, $\e$ and $\d$ mean the same as in Lemma \ref{01+}. Consider self-adjoint operators $A$ and $B$ such that $\|A\|\le c$, $\|B\|\le c/2$ and $\|B-A\|\le\min\{\d/2,c/2\}$. Let $K$ be a self-adjoint operator
such that $\|K\|=\d/2$. Then $\|B+K\|\le c$, $\|B-A\|\le\|K\|$ and $\|B-A+K\|\le2\|K\|$. Therefore,
$$
\|f(B+K)-f(B)-\bs{d}_{f,B}K\|\le\e\|K\|,
$$
$$
\|f(B)-f(A)-\bs{d}_{f,A}(B-A)\|\le\e\|B-A\|\le\e\|K\|
$$
and
$$
\|f(B+K)-f(A)-\bs{d}_{f,A}(B-A+K)\|\le\e\|B-A+K\|\le2\e\|K\|.
$$
Using the equality 
$$
\bs{d}_{f,A}(B-A+K)=\bs{d}_{f,A}(B-A)+\bs{d}_{f,A}K,
$$
we obtain
\begin{align*}
\|\bs{d}_{f,B}K-\bs{d}_{f,A}K\|&\le\|\bs{d}_{f,B}K-f(B+K)+f(B)\|\\[.2cm]
&+\|f(B+K)-f(A)-\bs{d}_{f,A}(B-A+K)\|\\[.2cm]
&+\|\bs{d}_{f,A}(B-A)-f(B)+f(A)\|\le4\e\|K\|,
\end{align*}
whence it follows that $\|\bs{d}_{f,B}-\bs{d}_{f,A}\|\le4\e$. $\bl$

\medskip

We proceed now to the case of operator differentiable functions. We say that {\it not necessarily bounded self-adjoint operators $A$ and $B$ are equivalent} if there exists an operator $K$ in $\mB_{\rm sa}(\h)$ such that $B=A+K$. For operators in the same equivalence class we can introduce the metric $\dist(A,B)\df\|B-A\|$.

\begin{thm} 
\label{odineo+}
Let $f$ be an operator differentiable function on $\R$. Then
the map
$
A\mapsto\bs{d}_{f,A}
$
on each equivalence class is continuous in the operator norm.
\end{thm}

\begin{lem} 
\label{ravnots+}
Let $f$ be an operator differentiable function on $\R$. Then for each $\e>0$  
there exists $\d>0$ such that 
$$
\|f(A+K)-f(A)-\bs{d}_{f,A}K\|\le\e\|K\|
$$
for every (not necessarily bounded) self-adjoint operator $A$ and for every self-adjoint operator $K$ whose norm is at most $\d$.
\end{lem}

\Pf Assume the contrary. Then for some $\e>0$, there exist 
sequences of self-adjoint operators $\{A_n\}_{n=1}^\be$ and $\{K_n\}_{n=1}^\be$
such that $\|K_n\|\to0$ and
\bay
\label{verotsep+}
\|f(A_n+K_n)-f(A_n)-\bs{d}_{f,A_n}K_n\|>\e\|K_n\|
\ey
for every $n\ge1$. Let $\A$ and ${\mathcal K}_n$ be the operators on $\ell^2(\h)$,
defined by \rf{oprA+} and \rf{oprKn+}. Since
$f$ is differentiable in the sense of Fr\'echet at the point $\A$, there exists
$\d>0$ such that
$$
\|f(\A+K)-f(\A)-\bs{d}_{f,\A}K\|\le\e\|K\|
$$
for all self-adjoint operators $K$ of norm at most $\d$. 
Applying this inequality to the operator ${\mathcal K}_n$
for sufficiently large $n$, we obtain
$$
\|f(A_n+K_n)-f(A_n)-\bs{d}_{f,A_n}K_n\|=
\|f(\A+{\mathcal K}_n)-f(\A)-\bs{d}_{f,\A}{\mathcal K}_n\|\le\e\|{\mathcal K}_n\|=\e\|K_n\|
$$
which contradicts inequality \rf{verotsep+}. $\bl$

\medskip

{\bf The proof of Theorem \ref{odineo+}.}
Let $\e$ and $\d$ mean the same as in Lemma \ref{ravnots+}. Consider self-adjoint operators $A$ and $B$ such that $\|B-A\|\le\d/2$. Let $K$ be a self-adjoint operator
such that $\|K\|=\d/2$. Then  $\|B-A\|\le\|K\|$ and $\|B-A+K\|\le2\|K\|$. Therefore,
$$
\|f(B+K)-f(B)-\bs{d}_{f,B}K\|\le\e\|K\|,
$$
$$
\|f(B)-f(A)-\bs{d}_{f,A}(B-A)\|\le\e\|B-A\|\le\e\|K\|
$$
and
$$
\|f(B+K)-f(A)-\bs{d}_{f,A}(B-A+K)\|\le\e\|B-A+K\|\le2\e\|K\|.
$$
Using the equality 
$
\bs{d}_{f,A}(B-A+K)=\bs{d}_{f,A}(B-A)+\bs{d}_{f,A}K$,
we obtain
\begin{align*}
\|\bs{d}_{f,B}K-\bs{d}_{f,B}K\|&\le\|\bs{d}_{f,B}K-f(B+K)+f(B)\|\\[.2cm]
&+\|f(B+K)-f(A)-\bs{d}_{f,A}(B-A+K)\|\\[.2cm]
&+\|\bs{d}_{f,A}(B-A)-f(B)+f(A)\|\le4\e\|K\|
\end{align*}
for all self-adjoint operators $K$ such that $\|K\|=\d/2$,
whence it follows that \lb$\|\bs{d}_{f,B}-\bs{d}_{f,A}\|\le4\e$,
whenever $\|B-A\|\le\d/2$. $\bl$

\begin{thm}
Let  $f\in\OL_{\rm loc}(\R)$. Then $f$ is locally operator differentiable if and only if the map $A\mapsto \bs{d}_{f,A}$ is continuous as a map
from the Banach space $\mB_{\rm sa}(\h)$
to the Banach space of bounded operators from $\mB_{\rm sa}(\h)$ to $\mB(\h)$.
\end{thm}

\Pf It follows from Theorem \ref{ndlod+} that it suffices to verify that the continuity of the map $A\mapsto \bs{d}_{f,A}$ implies operator differentiability.
Note that $H_{f,A,K}'(s)= \bs{d}_{f,A+sK}K$ (the derivative is taken in the strong operator topology). Therefore,
\bay
\label{prirost+}
f(A+K)-f(A)=\int_0^1(\bs{d}_{f,A+sK}K)\,ds,
\ey
where the integral is understood in the following sense:
$$
(f(A+K)-f(A))u=\int_0^1\big((\bs{d}_{f,A+sK}K)u\big)\,ds
$$
for every $u\in\h$. Applying \rf{prirost+} to the operator $tK$ instead of $K$, we obtain
$$
t^{-1}(f(A+K)-f(A))-\bs{d}_{f,A}K=\int_0^1\big((\bs{d}_{f,A+stK}-\bs{d}_{f,A})K\big)\,ds.
$$
Assume that $\|K\|=1$. Then it follows from the last identity that 
$$
\|t^{-1}(f(A+K)-f(A))-\bs{d}_{f,A}K\|\le\int_0^1\|\bs{d}_{f,A+stK}-\bs{d}_{f,A}\|\,ds
$$
It remains to observe that
$
\lim_{t\to0}\int_0^1\|\bs{d}_{f,A+stK}-\bs{d}_{f,A}\|\,ds=0
$
uniformly in all self-adjoint operators $K$ of norm 1
because of the continuity of the map $A\mapsto \bs{d}_{f,A}$. $\bl$

\medskip

The following result can be proved in a similar way.

\begin{thm}
\label{neprclass+}
Let  $f\in\OL_{\rm loc}(\R)$. Then $f$ is operator differentiable if and only if
the map $A\mapsto \bs{d}_{f,A}$ is continuous in the operator norm on every equivalence class.
\end{thm}

\begin{thm} The set ${\rm OD}(\R)$ is a closed subspace of $\OL(\R)$. 
\end{thm}

\Pf We have to prove that if $\lim\limits_{n\to\be}f_n=f$  in  $\OL(\R)$
and $f_n\in{\rm OD}(\R)$ for every $n$, then $f\in{\rm OD}(\R)$.
It follows from Theorems  \ref{sildif} and \ref{muld} that
$$
\|\bs{d}_{f_n,A}-\bs{d}_{f,A}\|=\|\bs{d}_{f_n-f,A}\|
\le\|\dg(f_n-f)\|_{\fM(\R\times\R)}=\|f_n-f\|_{\OL(\R)}\to0
$$
as $n\to\be$. Thus, $\lim\limits_{n\to\be}\bs{d}_{f_n,A}=\bs{d}_{f,A}$ in the norm uniformly in all self-adjoint operators $A$.
It remains to apply \rf{neprclass+} because continuity is preserved under uniform convergence. $\bl$

Note here that in the case of functions on finite intervals, the fact that the set of operator differentiable functions is closed in the space of operator Lipschitz functions was established in \cite{KS}. Moreover, it was also shown there that in this case the space of operator differentiable functions coincides with the closure of the set of polynomials in the space of operator Lipschitz functions. Note also that the question of operator differentiability of differentiable functions was posed in \cite{Wid}.

\medskip

\section{\bf Commutator Lipschitzness}
\setcounter{equation}{0}
\label{CoL}

\medskip

Recall that a continuous function $f$ on $\R$ is called {\it commutator Lipschitz} if 
\bay
\label{fAR-RfA}
\|f(A)R-Rf(A)\|\le\const\|AR-RA\|
\ey
for every bounded self-adjoint operator $A$ and for every bounded linear operator $R$. As in the definition of operator Lipschitz functions, if $f$ is commutator Lipschitz, then inequality \rf{fAR-RfA} holds for an arbitrary (not necessarily bounded) self-adjoint operator $A$ and for every bounded linear operator $R$, see Theorem \ref{anbnrn}.

Later we will see that that the following result holds.

\begin{thm}
\label{komlipf}
Let $f$ be a continuous function on $\R$. Then the following statements are equivalent:

{\em(a)} $\|f(A)-f(B)\|\le\|A-B\|$ for arbitrary self-adjoint operators $A$ and $B$;

{\em(b)} $\|f(A)R-Rf(A)\|\le\|AR-RA\|$ for every self-adjoint operator $A$ and for every  bounded linear operator $R$;

{\em(c)} $\|f(A)R-Rf(B)\|\le\|AR-RB\|$ for arbitrary self-adjoint operators $A$ and $B$ and for every  bounded linear operator $R$.
\end{thm}

Operators of the form $f(A)R-Rf(B)$ are called {\it quasicommutators}.

We will deduce Theorem \ref{komlipf} in \S\:\ref{oplip} from a more general result for normal operators in \S\:\ref{oplip}. Note, however, that in the case of functions of normal operators, commutator Lipschitzness is by no means equivalent to operator Lipschitzness.

\medskip

\section{\bf Operator Bernstein's inequalities}
\setcounter{equation}{0}
\label{Bern}

\medskip

In this section we give an elementary proof of the result of
\cite{Pe3} that functions in $L^\be(\R)$ whose Fourier transforms have compact support must be operator operator Lipschitz. Moreover, we obtain the so-called operator Bernstein's inequality with constant 1. We follow the approach of \cite{AP4}. 
We also obtain similar results for functions on the circle.

In \S\:\ref{Dost} we will deduce  from these results that the membership in the Besov class $B_{\be,1}^1(\R)$ is a sufficient condition for operator Lipschitzness.

Let $\s>0$. Recall that an entire function $f$ has {\it exponential type at most} $\s$, if for each $\e>0$,
there is $c>0$ such that $|f(z)|\le c e^{(\s+\e)|z|}$ for every $z\in\C$.

We denote by $\E_\s$ the set of entire functions of exponential type at most $\s$. It is well known that
$\E_\s\cap L^\be(\R)=\{f\in L^\be(\R):\supp\F f\subset[-\s,\s]\}$.

Note also that the space $\E_\s\cap L^\be(\R)$
coincides with the set of entire functions $f$ such that  $f\in L^\be(\R)$ and
\bay
\label{fEs}
|f(z)|\le e^{\s|\im z|}\|f\|_{L^\infty(\R)},\quad z\in\C,
\ey
see, for example, \cite{L}, page 97.

Bernstein's inequality (see \cite{Be}) says that
$$
\sup_{x\in \R}|f^\prime(x)|\le\s\sup_{x\in \R}|f(x)|
$$
for every $f$ in $\mathscr E_\s\cap L^\be(\R)$. It implies that
\bay
\label{corber}
|f(x)-f(y)|\le\s\|f\|_{L^\infty(\R)}|x-y|,\quad f\in\mathscr E_\s\cap L^\be(\R),\quad x,~y\in\R,
\ey
where $\|f\|_{L^\infty(\R)}\df\sup\limits_{x\in \R}|f(x)|$.

Bernstein also proved in \cite{Be} the following improvement of \rf{corber}:
\bay
\label{bern2}
|f(x)-f(y)|\le\b(\s(|x-y|))\|f\|_{L^\infty(\R)},\quad f\in\mathscr E_\s\cap L^\be(\R),\quad x,~y\in\R,
\ey
where
$$
\b(t)\df\left\{\begin{array}{ll}2\sin(t/2),&\text {if}\,\,\,\,0\le t\le\pi,\\[.1cm]
2,&\text {if}\,\,\,\,t>\pi.
\end{array}\right.
$$
Note that $\b(t)\le\min(t,2)$ for every $t\ge0$.

Let $X$ be a complex Banach space. We denote by $\E_\s(X)$ the space of
entire $X$-valued functions $f$ of exponential type at most $\s$, i.e., satisfying the following condition:
for each positive $\e$,
there exists $c>0$ such that $\|f(z)\|_X\le c e^{(\s+\e)|z|}$ for every $z\in\C$.

\medskip

{\bf Bernstein's inequality for vector-valued functions.} {\em Let $f$ be a function in 
$\E_\s(X)\cap L^\be(\R,X)$, where $\s>0$. 
Then 
\bay
\label{bn}
\|f(x)-f(y)\|_X\le\b(\s(|x-y|))\|f\|_{L^\infty(\R,X)}\le \s\|f\|_{L^\be(\R,X)}|x-y|
\ey
for all $x,y\in\R$.}

 \medskip
 
 The vector version of Bernstein's inequality reduces to the scalar version with the help of the Hahn--Banach theorem.

\medskip

\begin{thm}
\label{opbern}
Let $f\in\mathscr E_\s\cap L^\be(\R)$.  Then
\bay
\label{opnerBer}
\|f(A)-f(B)\|\le\b(\s(\|A-B\|))\|f\|_{L^\infty}\le\s\|f\|_{L^\infty}\|A-B\|
\ey
for arbitrary (bounded) self-adjoint operators $A$ and $B$. In particular,
$\|f\|_{\OL(\R)}\lb\le\s\|f\|_{L^\be(\R)}$.
\end{thm}

{\bf The proof of Theorem \ref{opbern}.} Let $A$ and $B$ be self-adjoint
operators on a Hilbert space $\h$. We have to show that
$$
\|f(A)-f(B)\|\le\b(\s\|A-B\|)\|f\|_{L^\infty}.
$$
Put $F(z)=f(A+z(B-A))$. Clearly, $F$ is an entire function with values in the space of operators
$\mB(\h)$ and $\|F(t)\|\le\|f\|_{L^\be(\R)}$ for every $t\in\R$. It follows from von Neumann's inequality (see \cite{SNF}) that 
$F\in\E_{\s\|B-A\|}(\mB(\h))$. To complete the proof, it remains to apply Bernstein's inequality \rf{bn}  to the vector function  $F$ for $x=0$ and $y=1$.  $\bl$

\medskip

Earlier, it was shown in \cite{Pe3} that
\bay
\label{PeR}
\|f\|_{\OL(\R)}\le\const\s\|f\|_{L^\be(\R)},\quad f\in\mathscr E_\s\cap L^\be(\R).
\ey
In particular, $\E_\s\cap L^\be(\R)\subset\OL(\R)$. It follows that 
for every $f\in\E_\s\cap \Li(\R)$, the function $f'$ is operator Lipschitz.

The following example shows that $\E_\s\cap \Li(\R)\not\subset\OL(\R)$.

\medskip

{\bf Example.} {\it Consider the function $f(x)\df\int_0^x{\rm Si}(t)\,dt$, where $\rm Si$  
is the integral sine,
$$
{\rm Si}(x)\df\int_0^x\frac{\sin t}t dt.
$$
Clearly, $f\in\E_1\cap \Li(\R)$,  but $f$ cannot be operator Lipschitz (see theorems {\rm\ref{olm}} and {\rm\ref{pro}} below) for the limit 
$\lim\limits_{|x|\to\be}x^{-1}f(x)$ does not exist (actually, 
$\lim\limits_{x\to\be}x^{-1}f(x)=\lim\limits_{x\to\be}{\rm Si}(x)=\frac\pi2=-\lim\limits_{x\to-\be}x^{-1}f(x)$)}.

\medskip

It is interesting to observe that if we slightly ``corrupt'' the function $f$ in this example by replacing it with the function
$g(x)\df\int_0^x{\rm Si}(|t|)\,dt$, it becomes operator Lipschitz.
It suffices to see that the function $g(x)-\frac\pi2x$ is operator Lipschitz.
This follows (see Proposition 7.8 of \cite{BS3}) from the fact that the derivative of this function belongs to the space $L^2(\R)\cap\Li(\R)$
(this can also be deduced from Theorem \ref{Besdiffer} below).

\medskip

Let us obtain now analogs of Bernstein's inequality for unitary operators. 

\begin{lem} 
\label{IUV}
Let $U$ and $V$ be unitary operators. Then there exists a self-adjoint operator
$A$ such that $V=e^{{\rm i}A}U$, $\|A\|\le\pi$ and
$\b(\|A\|)=\|U-V\|$.
\end{lem}

\Pf We define the operator $A$ by $A=\arg(VU^{-1})$, where the function $\arg$
is defined on $\T$  by $\arg(e^{{\rm i}s})=s, s\in[-\pi,\pi)$.
Obviously, 
$\b(\|A\|)=\|I-e^{{\rm i}A}\|=\|U-V\|$. $\bl$

\begin{thm}  
\label{uniber}
Let $f$ be a trigonometric polynomial of degree at most $n$.
Then
$$
\|f(U)-f(V)\|\le n\|f\|_{L^\be(\T)}\|U-V\|
$$
for arbitrary unitary operators $U$ and $V$.
\end{thm}

\Pf Let $A$ be a self-adjoint operator whose existence is guaranteed by
Lemma \ref{IUV}. Put $\Phi(z)\df f(e^{{\rm i}zA}U)$, $z\in\C$, where the same symbol $f$ stands for the analytic extension of $f$ to $\C\setminus\{0\}$.
Clearly, $\Phi$ is an entire function with values in $\mB(\h)$ and $\|\Phi(t)\|\le\|f\|_{L^\be(\T)}$ for every $t\in\R$. It follows from von Neumann's inequality (see \cite{SNF}) that 
$\Phi\in\E_{n\|A\|}(\mB(\h))$. Applying Bernstein's inequality for vector functions, we obtain
$$
\|f(U)-f(V)\|=\|\Phi(1)-\Phi(0)\|\le\b(n\|A\|)\|f\|_{L^\be(\T)}.
$$
It remains to observe that
$\b(n\|A\|)\le n \b(\|A\|)=n\|U-V\|$. $\bl$

\medskip

Note that it was shown in \cite{Pe1} that
$\|f(U)-f(V)\|\le\const n\|f\|_{L^\be(\T)}\|U-V\|$
for every trigonometric polynomial $f$ of degree $n$ and for arbitrary unitary operators $U$ and $V$.

\medskip

{\bf Remark.} 
It can be seen from the proof of Theorem \ref{uniber} that 
$$
\|f(U)-f(V)\|\le\b(n\|A\|)\|f\|_{L^\be(\T)}=\b\left(2n\arcsin\frac{\|U-V\|}2\right)\|f\|_{L^\be(\T)}.
$$
This estimate is best possible for the function $f(z)=z^n$ because
$$
\sup\big\{|z_1^n-z_2^n|:~z_1\in\T,~z_2\in\T,~
|z_1-z_2|<n\big\}=\b\big(2n\arcsin(\d/2)\big),\quad\d\in(0,2].
$$

\medskip

\section{\bf Necessary conditions for operator Lipschitzness}
\setcounter{equation}{0}
\label{Neob}

\medskip

In this section we obtain necessary conditions for operator Lipschitzness for functions on the line and on the circle. These necessary conditions were obtained substantially in the papers\cite{Pe1} and \cite{Pe3}, in which other methods were used. Here to achieve the purpose, besides the trace class criterion for Hankel operators
(see \S\:\ref{Prel}, Subsection 3), which was also used \cite{Pe1} and \cite{Pe3}, we use the results of Section \ref{olriolt+} of this survey on the behavior of the derivatives of operator Lipschitz functions under linear-fractional transformations.

To prove the next result, we are going to use results of Section 
\ref{yadiyadcom} on the behavior of functions of operators under trace class perturbations. 

\begin{thm}
\label{B11}
Let $f$ be an operator Lipschitz function on $\T$. 
Then $f\in B_1^1(\T)$.
\end{thm}

\Pf By the remark to Theorem \ref{tsentrez}, the function $f$ has the property:
$f(U)-f(V)\in\bS_1$,
whenever $U$ and $V$ are unitary operators such that $U-V\in\bS_1$.

We define the operators $U$ and $V$ on $L^2(\T)$ by
$$
Uf=\bar z f\quad\mbox{and}\quad Vf=\bar zf-2(f,\1)\bar z,\quad f\in L^2.
$$
It is easy to see that $U$ and $V$ are unitary operators and
$\rank(V-U)=1.$
It is also easy to verify that for $n\ge0$, the following equality holds:
$$
V^nz^j=\left\{
\begin{array}{ll}z^{j-n},&j\ge n,\\[.1cm]
-z^{j-n},&0\le j<n,\\[.1cm]
z^{j-n},&j<0.
\end{array}\right.
$$
Hence, for every continuous function $f$ on $\T$,
the following identity holds:
\begin{align*}
\big((f(V)-f(U))z^j,z^k\big)&=\sum_{n>0}\hat f(n)\big((V^nz^j,z^k)-(z^{j-n},z^k)\big)\\
&+\sum_{n<0}\hat f(n)\big((V^nz^j,z^k)-(z^{j-n},z^k)\big)\\
&=-2\left\{\begin{array}{ll}\hat f(j-k),&j\ge0,~k<0,\\[.1cm]
\hat f(j-k), &j<0,~k\ge0,\\[.1cm]
0,&\mbox{otherwise}.
\end{array}\right.
\end{align*}

Thus, if $f(U)-f(V)\in\bS_1$, then the operators on $\ell^2$
with Hankel matrices
\lb$
\{\hat f(j+k)\}_{j\ge0,k\ge1}$ and $\{\hat f(-j-k)\}_{j\ge0,k\ge1}
$
belong to $\bS_1$. Now we can use the trace class criterion for Hankel operators (see \S\:\ref{Prel}, Subsection 3) and conclude that $f\in B_{1}^{1}(\T)$. $\bl$

\medskip

Note here that the construction in the proof of Theorem \ref{B11} is taken from the paper \cite{AP3}.

It is convenient to introduce in this section the notation $\zh M\zh$ for the norm of a matrix $M$.

To state a corollary to Theorem \ref{B11}, we need the Banach space
$(\OL)_{\rm loc}'(\T)$, which will be studied in detail in \S\:\ref{olriolt+}.
Here we only mention that $(\OL)_{\rm loc}'(\T)=\{f'+c\ov z: f\in\OL(\T),\:c\in\C\}$ (see Corollary \ref{rt4}). Moreover,
as always in this paper, the derivative is understood in the complex sense, i.e.,
$f'(\z)\df\lim_{\t\to\z}(\t-\z)^{-1}(f(\t)-f(\z))$.

\begin{cor} 
\label{sumgrad}
Let $u$ be the Poisson integral of a function $f$ in $(\OL)_{\rm loc}'(\T)$. Then 
$\zh\n u\,\zh\in L^1(\dd)$ and 
$\big\|\zh\n u\,\zh\big\|_{L^{1}(\dd)}\le\const\|f\|_{(\OL)_{\rm loc}'(\T)}$. 
\end{cor}

\Pf  Let  $f=g'$, where $g\in\OL(\T)$. Then $f\in B_1^0(\T)$ for $g\in B_1^1(\T)$. It suffices to use the characterization of the Besov class $B_1^0(\T)$ in terms of harmonic extension,  see \S\:\ref{Prel}. It remains to observe that the conclusion of the corollary is obvious for the function $f(z)=z^{-1}=\ov z$. $\bl$

\medskip

To state a stronger necessary condition for operator Lipschitzness, we need the notion of Carleson measures. Let $\mu$ be a positive Borel measure in the open unit disk 
 $\dd$. The well known Carleson theorem says that
the Hardy class $H^p$ is contained in $L^p(\mu)$ ($0<p<+\be$) if and only if
for every point $\z$ of the unit circle $\T$ and for every positive $r$ 
the following inequality holds:
$$
\mu\{z\in\dd:\,|z-\z|<r\}\le\const r.
$$
Such measures $\mu$ are called {\it Carleson measures} in the disk $\dd$.
Note that the Carleson condition does not depend on $p\in(0,+\be)$.
More detailed information about Carleson measures can be found, for example, in \cite{Ni1} and \cite{N}. We need the following equivalent 
reformulation of the Carleson condition:
\bay
\label{Vin}
\sup_{a\in\dd}\int_\dd\frac{1-|a|^2}{|1-\ov a z|^2}\,d\mu(z)<+\be,
\ey
see, for example,  \cite{Ni1}, Lecture VII. 
Note that \rf{Vin} means that $\|k_a\|_{L^2(\mu)}\le\const\|k_a\|_{H^2}$
for every $a\in\dd$, where $k_a(z)\df(1-z\ov a)^{-1}$ is the reproducing kernel of the  Hilbert space $H^2$.

{\it We denote by ${\rm CM}(\dd)$ the space of complex Radon measures $\mu$
in $\dd$ such that $|\mu|$  is a Carleson measure and by $\|\mu\|_{{\rm CM}(\dd)}$ the norm of the identical imbedding from $H^1$ to  $L^1(|\mu|)$}. It is well known that
the (quasi)norm of the operator of identical imbedding from $H^p$ to $L^p(|\mu|)$ is equal to
$\|\mu\|_{{\rm CM}(\dd)}^{1/p}$  for all $p\in(0,+\infty)$.

Everything said about Carleson measures in $\dd$ has natural analogs in the half-plane 
$\C_+$.
In this case the Carleson condition for a positive Borel $\mu$ in $\C_+$
can be rewritten in the following way:
$$
\mu\{z\in\C_+:\,|z-t|<r\}\le\const r
$$
for every $t\in\R$ and every $r>0$. An analogue of \rf{Vin} is the following inequality:
$$
\sup_{a\in\C_+}\int_{\C_+}\frac{\im a}{|z-\ov a |^2}\,d\mu(z)<\be.
$$
In particular, in the same way we can introduce the space ${\rm CM}(\C_+)$ as well as the norm in this space.

Let $f$ be a function  (distribution) on the unit circle $\T$.
{\it We denote by $\mP f$  the Poisson integral of $f$}.

\begin{thm} 
\label{PfCM}
Let  $f\in(\OL)'_{\rm loc}(\T)$.
Then $\zh\n (\mP f)\,\zh\,d\m_2\in\CM(\dd)$.
\end{thm}

\Pf Let $f\in(\OL)_{\rm loc}'(\T)$. Then it follows from Theorem \ref{kp} and Corollary \ref{sumgrad}
that
\bay
\label{afia}
\int_\dd|\!\!|\!|((\nabla u)\circ\f)(z)|\!\!|\!|\cdot|\f'(z)|\,d\m_2
\le\const\|f\|_{(\OL)_{\rm loc}'(\T)}
\ey
for every linear-fractional automorphism $\f$ of the unit disk $\dd$,
where $u=\mP f$. Put now $\f(z)\df(1-\ov a z)^{-1}(a-z)$, where $a\in\dd$.  Making the change of variable $z=\f(w)$ in the integral in \rf{afia}, we obtain
$$
\sup_{a\in\dd}
\int_\dd|\!\!|\!|(\nabla u)(w)|\!\!|\!|\frac{1-|a|^2}{|1-\ov a w|^2}\,d\m_2(w)\le\const\|f\|_{(\OL)'(\T)}.
$$
Therefore, the measure $|\!\!|\!|\nabla u\,|\!\!|\!|\,d\m_2=\zh\n (\mP f)\,\zh\,d\m_2$ satisfies condition \rf{Vin}. $\bl$

\medskip

The following result is a reformulation of Theorem \ref{PfCM}.

\begin{thm} 
Let $f\in\OL(\T)$.
Then\footnote{Here and in what follows ${\rm Hess}$ denotes Hessian, i.e.,
the matrix of second order partial derivatives.} 
$\zh{\rm Hess}\, \mP f\,\zh\,d\m_2\in\CM(\dd)$.
\end{thm}
  
We proceed now to Poisson integrals of functions on $\R$.
If $f\in L^1(\R,(1+x^2)^{-1}\,dx)$, then the Poisson integral can be defined in the standard way.   
We need the Poisson integral of functions $f$ such that $f'\in L^1(\R,(1+x^2)^{-1}\,dx)$.
Clearly, it suffices to consider the a real function $f$. Let $u$ be the Poisson integral of  $f'$. We denote by $v$ the harmonic conjugate of $u$.
The function $u+{\rm i}v$ has a primitive $F$ such that the boundary value function of $\re F$
coincides with $f$ everywhere on $\R$. The function $F$ is not determined uniquely because the harmonic conjugate $v$ is not determined uniquely. The family $\{v+c\}_{c\in\R}$ consists of all functions harmonically conjugate to. We need a primitive of
 $u+{\rm i}(v+c)$ in the form $F+c{\rm i}z+{\rm i}\a$, где $\a\in\R$.
Note that $\re(F+c{\rm i}z+{\rm i}\a)=\re F-cy$. Thus, it is natural to define the poisson 
integral of $f$ as the class of functions $\{\re F-cy\}_{c\in\R}$. Since ${\rm Hess}\,y=0$,
the Hessian of the Poisson integral ${\rm Hess}\,\mP f$ of $f$ is determined uniquely.
In the next statement $(\OL)'(\R)\df\{f': f\in\OL(\R)\}$.

\begin{thm}
Let  $f\in(\OL)'(\R)$.
Then $|\!\!|\!|\n \mP f\,|\!\!|\!|\,d\m_2\in\CM(\C_+)$.
\end{thm}

\Pf 
Let $f\in(\OL)'(\R)$. Then it follows from Theorem \ref{pere} and Corollary \ref{sumgrad} that
\bay
\label{afia1}
\int_\dd|\!\!|\!|((\nabla u)\circ\f)(z)|\!\!|\!|\cdot|\f'(z)|\,d\m_2\le\const\|f\|_{(\OL)'(\R)},
\ey
for every automorphism $\f$ in $\dlC$ such that $\f(\dd)=\C_+$,
where $u=\mP f$.
Take now $\f(z)\df(1-z)^{-1}(a-\ov a z)$, where $a\in\C_+$.  Making the substitution
$z=(w-\ov a)^{-1}(w-a)$ in the integral in
\rf{afia}, we obtain
$$
\sup_{a\in\C_+}
\int_{\C_+}|\!\!|\!|(\nabla u)(w)|\!\!|\!|\frac{2\im a}{|w-\ov a|^2}\,d\m_2(w)\le\const\|f\|_{(\OL)_{\rm loc}'(\T)}.
$$
The last condition is equivalent to the fact that the measure $|\!\!|\!|\nabla u\,|\!\!|\!|\,d\m_2$ is Carleson. $\bl$

\begin{thm} 
Let  $f\in\OL(\R)$.
Then $|\!\!|\!|{\rm Hess}\, \mP f\,|\!\!|\!|\,d\m_2\in\CM(\C_+)$.
\end{thm}

The necessary conditions for operator Lipschitzness given above were obtained originally in \cite{Pe1} and \cite{Pe3}. Namely, it was shown in \cite{Pe1} that if $f\in\OL(\T)$, then both Hankel operators $H_f$ and $H_{\ov f}$ map the Hardy class $H^1$ to the Besov class $B_1^1(\T)$
(the class of such functions $f$ is denoted in \cite{Pe1} by $L$). S. Semmes observed that $f\in L$ if and only if the measure  
$\zh{\rm Hess}\,\mP f\zh\,d\m_2$ is Carleson; see \cite{Pe4}, where the proof of this equivalence is given. A similar result also holds for functions on $\R$, see
\cite{Pe3}. It was also shown in \cite{Pe1} that the necessary condition for operator Lipschitzness discussed above is not sufficient. Moreover, it is not even sufficient for Lipschitzness.

Consider now the spaces $\pp_+(b_\infty^{-1}(\T))$ and 
$\pp_+(b^{-1}_{1,\be}(\T))$, the closures of the set of analytic polynomials in the Besov spaces $B_\infty^{-1}(\T)$ and $B^{-1}_{1,\be}(\T)$.
It is well known that these spaces admit the following descriptions in terms of analytic extension to the unit disk:
$$
\pp_+(b_\infty^{-1}(\T))=
\big\{h:~\lim_{r\to1^-}(1-r)\|h(rz)\|_{L^\be(\T)}=0\big\};
$$
$$
\pp_+(b^{-1}_{1,\be}(\T))=
\big\{h:~\lim_{r\to1^-}(1-r)\|h(rz)\|_{L^1(\T)}=0\big\}.
$$
It was also observed in \cite{Pe1} that the space $\pp_+L$ is dual to the space of analytic functions of functions $g$ in $\dd$ that admit a representation
\bay
\label{predprobes}
\hspace*{-.4cm}g\!=\!\sum_n\f_n\psi_n,~\mbox{where}~\f_n\!\in\!H^1,~
\psi_n\!\in\!\pp_+(b_\infty^{-1}(\T)),~
\sum_n\|\f_n\|_{H^1}\|\psi_n\|_{\pp_+(b_\infty^{-1}(\T))}\!<\be.
\ey
Obviously, such functions $g$ belong to the space $\pp_+(b^{-1}_{1,\be}(\T))$ whose dual space can be identified naturally with the Besov space $\pp_+B^1_{\be,1}(\T)$ of functions analytic in $\dd$. Nevertheless, not every function in $\pp_+(b^{-1}_{1,\be}(\T))$ can be represented in the form \rf{predprobes}. Otherwise,  we would have the equality $L=B_{\be,1}^1(\T)$ which is impossible because the condition $f\in L$, being necessary for operator Lipschitzness, is not sufficient.

\medskip

{\bf Remark.}
Note here that the space
$$
L=\{f\in {\rm BMO}(\R):~\zh{\rm Hess}\,\mP f\zh\,d\m_2\in{\rm CM}(\C_+)\}
$$
is a limit space of the scale of the Triebel--Lizorkin spaces  
$F_{p,q}^s(\R)$; it is denoted by $F_{\be,1}^1(\R)$ 
(see \cite {FJ}, \S\:5). In a similar way we can define the Triebel--Lizorkin space $F_{\be,1}^1(\T)$ of functions on $\T$. The necessary condition for operator Lipschitzness obtained above can be reformulated in the following way:
$\OL(\R)\subset F_{\be,1}^1(\R)$ and $\OL(\T)\subset F_{\be,1}^1(\T)$.

\medskip

\section{\bf A sufficient condition for operator Lipschitzness in terms of Besov classes}
\setcounter{equation}{0}
\label{Dost}

\medskip

In this section we show that the functions in the Besov class $B_{\be,1}^1(\R)$ (see \S\:\ref{Prel}) are operator Lipschitz. We also obtain a similar result for functions on the unit circle. The proofs given here differ from the original proofs of \cite{Pe1} and \cite{Pe3} and are based on operator Bernstein's inequalities (see \S\:\ref{Bern}).

\begin{thm}
\label{Besdost}
Let $f\in B_{\be,1}^1(\R)$. Then $f$ is operator Lipschitz  and
\bay
\label{BeOLsa}
\|f(A)-f(B)\|\le\const\|f\|_{B_{\be,1}^1}\|A-B\|
\ey
for arbitrary self-adjoint operators $A$ and $B$ with bounded difference $A-B$.
\end{thm}

\Pf As we have observed in the introduction (see Theorem \ref{anbnrn} below), it suffices to prove \rf{BeOLsa}
for bounded self-adjoint operators $A$ and $B$.

Without loss of generality we may assume that $f(0)=0$.
Consider the functions $f_n=f*W_n$ defined by \rf{fn}.  
Put $g_n\df f_n-f_n(0)$. It follows from the definition of $B_{\be,1}^1(\R)$ 
(see \S\:\ref{Prel}) that
$\sum_{n=-\be}^\be g'_n=f'$ 
and the series converges uniformly on $\R$. Hence, the series 
$\sum\limits_{n=-\be}^\be g_n$ converges uniformly on each compact subset of $\R$. Thus, 
$$
\sum_{n=-\be}^\be g_n(A)=f(A)\quad\mbox{and}\quad\sum_{n=-\be}^\be g_n(B)=f(B),
$$
the series being absolutely convergent in the operator norm. 

Since obviously, $f_n\in \E_{2^n+1}\cap L^\be(\R)$, operator Bernstein's inequality \rf{opnerBer} allows us to conclude that
\begin{align*}
\|f(A)-f(B)\|&\le\left\|\sum_{n=-\be}^\be\big((g_n(A)-g_n(B)\big)\right\|
=\left\|\sum_{n=-\be}^\be\big((f_n(A)-f_n(B)\big)\right\|\\[.2cm]
&\le\sum_{n=-\be}^\be2^{n+1}\|f_n\|_{L^\be}\|A-B\|
\le\const\|f\|_{B_{\be,1}^1}\|A-B\|.\quad\bl
\end{align*}

In a similar way one can prove the following analog of Theorem \ref{Besdost} for functions on the unit circle.

\begin{thm}
\label{Besokr}
Let $f\in B_{\be,1}^1(\T)$. Then $f$ is operator Lipschitz and
$$
\|f(U)-f(V)\|\le\const\|f\|_{B_{\be,1}^1}\|U-V\|,\quad f\in B_{\be,1}^1(\T),
$$
for arbitrary unitary operators $U$ and $V$.
\end{thm}

The following statement allows us to to combine the necessary conditions obtained in \S\:\ref{Neob} with the sufficient conditions of this section.

\begin{thm}
$B_{\be,1}^1(\T)\subset\OL(\T)\subset F_{\be,1}^1(\T)$ and 
$B_{\be,1}^1(\R)\subset\OL(\T)\subset F_{\be,1}^1(\R)$.
\end{thm}

It turned out that the functions in $B_{\be,1}^1(\R)$ are not only operator Lipschitz, but also operator differentiable.

\begin{thm}
\label{Besdiffer}
Let $f\in B_{\be,1}^1(\R)$. Then $f$ is operator differentiable function.
\end{thm}

We refer the reader to \cite{Pe3} and \cite{Pe5}
for the proof of Theorem \ref{Besdiffer}.

\medskip

\section{\bf Operator H\"older functions}
\setcounter{equation}{0}
\label{oHold}

\medskip

We discuss here other applications of operator Bernstein's inequalities obtained in \S\:\ref{Bern}. We show that the class of operator H\"older functions of order $\a$, $0<\a<1$, coincides with the class of H\"older functions of order $\a$.
We also dwell briefly on the case of arbitrary moduli of continuity.
The results of this section were obtained in \cite{AP1} and \cite{AP2}. Another approach to these problems was found in \cite{FN} where the authors obtained similar results for functions in H\"older classes with somewhat worse constants as well as a somewhat weaker  result for arbitrary moduli of continuity.

It is well known that the (scalar) Bernstein inequality plays a key role in perturbation theory, see, for example, \cite{Ah}, \cite{Dz}, \cite{NI}, \cite{Ti}. We are talking about the description of smoothness type properties in terms of approximation by nice functions.
Direct theorems of approximation theory give us estimates of the rate of approximation of functions in a given function space $X$ (usually, of smooth functions in a sense) by nice functions. 
Inverse theorems allow one to conclude that a given function $f$ belongs to a certain function space if $f$ admits certain estimates of the rate of approximation by nice functions. In the case when for a functions space $X$ the direct theorems ``match'' the inverse theorems, we obtain a complete description of $X$ in terms of approximation by such nice functions.

In this section we consider functions spaces on the unit circle $\T$ and on the real line $\R$. In the first case the role of nice functions is played by the classes $\mathcal P_n$ of trigonometric polynomials of degree at most $n$, while in the second space by the classes $\E_\s$ of functions of exponential type at most $\s$. We are going to consider only uniform approximation.

The classical Bernstein inequalities plays a decisive role in the proof of inverse theorems of approximation theory. It can be seen that using the operator version of Bernstein's inequality in such a proof, we obtain the corresponding smoothness of a function $f$ on the set of unitary operators if we are dealing with functions on the circle and on the set of self-adjoint operators if we are dealing with functions on the line.

Let us illustrate this by an example. 
The classical Jackson theorem says that if $f$ is in the H\"older class
$\L_\a(\T)$, $0<\a<1$, then
\bay
\label{JackBern}
\dist(f,\mathcal P_n)\le\const(n+1)^{-\a}\|f\|_{\L_\a}.
\ey
Bernstein proved that the converse is also true, i.e., if for a function $f$ in $C(\T)$  inequalities \rf{JackBern} hold with $\a\in(0,1)$, then  
$f\in\L_\a(\T)$.

Let us give the standard prove of this Bernstein result.
Without loss of generality we may assume that $c=1$. 
For $n\ge0$, there exists a trigonometric polynomial $f_n$  
such that $\deg f_n<2^n$ and $\|f-f_n\|_{C(\T)}\le 2^{-\a n}$. Clearly, 
$$
\|f_n-f_{n-1}\|_{C(\T)}\le\|f-f_n\|_{C(\T)}+\|f-f_{n-1}\|_{C(\T)}\le2^{-\a n}(1+2^\a)\le3\cdot2^{-\a n}.
$$
Consequently, 
$$
\|f_n-f_{n-1}\|_{\Li(\T)}\le2^n\|f_n-f_{n-1}\|_{C(\T)}\le3\cdot2^{(1-\a)n}
$$
by Bernstein's inequality. Taking into account the obvious equality
$\|f_0\|_{\Li(\T)}=0$, we obtain
$$
\|f_N\|_{\Li(\T)}\le\sum_{n=1}^N\|f_n-f_{n-1}\|_{\Li(\T)}\le3\sum_{n=1}^N2^{(1-\a)n}\le\frac3{1-2^{\a-1}}2^{(1-\a)N},\quad N\in\Z_+.
$$
Let $\z,\,\tau\in\T$. We select $N\in\Z_+$
such that $2^{-N}<|\z-\tau|\le2^{1-N}$. 
Then
\begin{align*}
|f(\z)-f(\tau)|&\le|f(\z)-f_N(\z)|+|f_N(\z)-f_N(\tau)|+|f_N(\xi)-f(\tau)|\\[.2cm]
&\le2\|f-f_N\|_{L^\be}+\|f_N\|_{\Li}|\z-\tau|\le2\cdot 2^{-\a N}+\frac3{1-2^{\a-1}}2^{(1-\a)N}|\z-\tau|\\[.2cm]
&\le2|\z-\tau|^\a+
\frac{3\cdot2^{1-\a}}{1-2^{\a-1}}|\z-\tau|^\a\le\frac8{1-2^{\a-1}}|\z-\tau|^{\a}.\quad\bl
\end{align*}

The following theorem says that every function in $\L_\a(\T)$, $0<\a<1$, is operator H\"older of order $\a$ which is a sharp contrast with the case of Lipschitz functions.

\begin{thm} 
\label{43}
Let $f\in\L_\a(\T)$, where $\a\in(0,1)$. 
Then there exists a constant $c$ such that
$$
\|f(U)-f(V)\|\le c(1-\a)^{-1}\|f\|_{\L_\a}\|U-V\|^\a
$$
for arbitrary unitary operators $U$ and $V$.
\end{thm}

\Pf Let $f\in\L_\a(\T)$. First, we apply the direct theorem of approximation theory (the Jackson theorem in our case).  By this theorem, 
$$
\dist(f,\mathcal P_n)\le\const(n+1)^{-\a}\|f\|_{\L_\a},\quad n\in\Z_+.
$$ 
Repeating practically word for word the proof of the corresponding inverse theorem, replacing $\z$ and $\tau$ with unitary operators $U$ and $V$, and applying the operator Bernstein inequality instead of the scalar one, we arrive at the desired result. $\bl$

\medskip

A similar result also holds in the case of the line.

\begin{thm} 
\label{44}
Let $f\in\L_\a(\R)$, $\a\in(0,1)$. 
Then there exists a constant $c$ such that
$$
\|f(A)-f(B)\|\le c(1-\a)^{-1}\|f\|_{\L_\a}\|A-B\|^\a
$$
for arbitrary self-adjoint operators $A$ and $B$.
\end{thm}

The proof is based on a similar description of the function class $\L_\a(\R)$
in terms of approximation by entire functions of exponential type.

\medskip

{\bf The direct theorem for the space $\L_\a(\R)$.} {\em Let $f\in\L_\a(\R)$ and let $0<\a<1$. Then there exists a positive number $c$ such that 
\bay
\label{LaEs}
\inf\{h\in\E_\s:\|f-h\|_{L^\be(\R)}\}\le c\s^{-\a}\|f\|_{\L_\a(\T)}
\ey
for every $\s>0$}.

\medskip

{\bf The inverse theorem for $\L_\a(\R)$.} {\em Let $0<\a<1$ and let $f$ be
a continuous function on $\R$ such that $\lim\limits_{|x|\to\be} x^{-1}f(x)=0$. Suppose that {\em\rf{LaEs}} holds for a positive number $c$ and for all $\s>0$. Then 
$f\in\L_\a(\R)$ and
$\|f\|_{\L_\a(\R)}\le\frac{5c}{1-2^{\a-1}}.$
}

\medskip

The proofs of Theorems \ref{43} and \ref{44} were given in more detail in \cite{AP2}.
In the same paper a series of other results were obtained; they are based as a matter of fact on certain results of approximation theory.

In particular, analogs of Theorems \ref{43} and \ref{44} were obtained there for all $\a>0$. 
We state here some more results obtained in \cite{AP2} that also can be proved by using methods of approximation theory.

A function $\o:[0,+\be)\to\R$ is called a {\it modulus of continuity} if it is a nonnegative nondecreasing continuous function such that $\o(0)=0$, $\o(x)>0$ for $x>0$ and
$\o(x+y)\le\o(x)+\o(y)$ for all $x,y\in[0,+\be)$.

Denote by $\L_\o(\R)$ the space of continuous functions $f$ on $\R$
such that 
$$
\|f\|_{\L_\o(\R)}\df\sup_{x\ne y}\frac{|f(x)-f(y)|}{\o(|x-y|)}<\be.
$$
Similarly, we can define the space $\L_\o(\T)$.

Put 
\bay
\label{omzv}
\o_*(x)\df x\int_x^\be \frac{\o(t)}{t^2}\,dt.
\ey

\begin{thm} 
\label{prmone}
Let $f\in\L_\o(\R)$, where $\o$ is a modulus of continuity. Then
$$
\|f(A)-f(B)\|\le c\|f\|_{\L_\o(\R)}\o_*(\|A-B\|)
$$  
for arbitrary self-adjoint operators $A$ and $B$
$A$ and  $B$, where $c$ is a numerical constant.
\end{thm}

A similar result also holds for functions $f$ in $\L_\o(\T)$.

\medskip

\section{\bf H\"older functions under perturbations by operators of Schatten--von Neumann classes}
\setcounter{equation}{0}
\label{SchavN}

\medskip

We consider on this section one more application of operator Bertstein inequalities given in \S\:\ref{Bern}. Let $f$ be a function of H\"older class
$\L_\a(\R)$, $0<\a<1$, and let $p\ge1$. Suppose that $A$ and $B$ are self-adjoint operators and $B-A\in\bS_p$. What can we say about the operator $f(A)-f(B)$? This question was studied in detail in \cite{AP3}.
We state here the result of \cite{AP3} in the case $p>1$.

\begin{thm}
\label{aHolSp}
Let $p>1$ and $0<\a<1$. Then 
$$
\|f(A)-f(B)\|_{\bS_{p/\a}}\le const\|f\|_{\L_\a}\|A-B\|_{\bS_p}^\a
$$
for arbitrary self-adjoint operators $A$ and $B$ whose difference belongs to $\bS_p$.
\end{thm}

We omit here the proof of Theorem \ref{aHolSp} and refer the reader to \cite{AP3} for the proof. The case $p=1$ is also considered in detail in \cite{AP3}. Note here that the conclusion of Theorem \ref{aHolSp} is false for $p=1$. Also, in 
\cite{AP3} an analog of Theorem \ref{aHolSp} is obtained for all positive $\a$ and 
more general problems of perturbations by operators in symmetrically normed ideals are considered.

\

\begin{center}
\bf\Large Chapter II
\end{center}


\begin{center}
\bf\Large 
Schur multipliers and double operator integrals
\end{center}
\label{muSch}

\addtocontents{toc}{\vspace*{.3cm}\hspace*{-.35cm}\textbf{{\bf Chapter II}. Schur multipliers and double operator integrals}\hfill\pageref{muSch}}

\renewcommand{\thesection}{2.\arabic{section}}
\setcounter{section}{0}

\medskip

In this section we study Schur multipliers, both discrete ones and Schur multipliers with respect to spectral measures. We use a description of the discrete Schur multipliers that is based on Grothendick's theorem (see Pisier's book
\cite{Pi} and the paper \cite{Pi2}). We refine this result in the case when the initial function is defined on the product of topological spaces and is continuous in each variable. We also obtain a refinement of the general result for Borel functions on the product of topological spaces.

Then we define double operator integrals and Schur multipliers with respect to spectral measures. The study of such Schur multipliers in the case of Borel functions on the product of topological spaces reduces to discrete Schur multipliers. 

\medskip

\section{\bf Discrete Schur multipliers}
\setcounter{equation}{0}
\label{dmsh}

\medskip

We denote by $\ell^p({\mT})$ the space of complex function
$\a: t\mapsto\a_t$ defined on a not necessarily countable or finite set
 $\mT$ such that $\sum\limits_{t\in\mT}|\a_t|^p <\be$
with norm
$\|\a\|_p =\lb\Big(\sum_{t\in\mT}|\a_t|^p\Big)^{1/p}$, where $p\in[1,+\be)$.
For $p=\be$, the space $\ell^p(\mT)$ consists of all bounded complex functions 
 $\a: t\mapsto\a_t$ on $\mT$ and $\|\a\|_\be = \sup_{t\in\mT}|\a_t|$. 
In those cases when we have to specify the set  $\mT$, on which the family $\a$ is defined, we will write $\|\a\|_{\ell^p(\mT)}$
in place of $\|\a\|_p$. We denote by $c_0(\mT)$ the subspace of
$\ell^\be(\mT)$ that consists of the functions that tend to zero at infinity.

Let ${\mS}$ and $\mT$ be arbitrary nonempty sets. With each bounded operator
$A:\ell^2(\mT)\to\ell^2(\mS)$ one can associate a unique matrix $\{a(s,t)\}_{(s,t)\in\mS\times\mT}$ such that
$(Ax)_s=\sum\limits_{t\in T}a(s,t)x_t$ for every $x=\{x_t\}_{t\in\mT}$ in $\ell^2(\mT)$. In this case we say that
{\it the matrix $\{a(s,t)\}_{(s,t)\in\mS\times\mT}$ induces bounded operator} $A:\ell^2(\mT)\to\ell^2(\mS)$. Put
$$
\|\{a(s,t)\}_{(s,t)\in\mS\times\mT}\|\df\|A\|
\quad\mbox{and}\quad \|\{a(s,t)\}_{(s,t)\in\mS\times\mT}\|_{\bS_1}\df\|A\|_{\bS_1}.
$$
If $A\not\in\bS_1(\ell^2(\mT),\ell^2(\mS))$, we assume that the last norm equals $\be$. 
If the matrix \lb$\{a(s,t)\}_{(s,t)\in\mS\times\mT}$ does not induce a bounded operator from
$\ell^2(\mT)$ to $\ell^2(\mS)$, we assume that its operator norm (as well as its trace norm) equals $\be$.
We denote by $\mB(\mS\times\mT)$ {\it the set of matrices $\{a(s,t)\}_{(s,t)\in\mS\times\mT}$ that induce bounded operators from $\ell^2(\mT)$ to $\ell^2(\mS)$}. Sometimes we will write $\|\{a(s,t)\}_{(s,t)\in\mS\times\mT}\|_{\mB(\mS\times\mT)}$ in place of 
$\|\{a(s,t)\}_{(s,t)\in\mS\times\mT}\|$
and $\|\{a(s,t)\}_{(s,t)\in\mS\times\mT}\|_{\bS_1(\mS\times\mT)}$ in place of 
$\|\{a(s,t)\}_{(s,t)\in\mS\times\mT}\|_{\bS_1}$.

A matrix ${\Phi}=\{\Phi(s,t)\}_{(s,t)\in\mS\times\mT}$ is called a {\it Schur multiplier} 
of the space $\mB(\mS\times\mT)$
if for every matrix
$A=\{a(s,t)\}_{(s,t)\in\mS\times\mT}$ in $\mB(\mS\times\mT)$, the matrix 
${\Phi}\star A\df\{\Phi(s,t)a(s,t)\}_{(s,t)\in\mS\times\mT}$ also belongs to 
$\mB(\mS\times\mT)$.

{\it We denote by $\fM(\mS\times\mT)$ the set of Schur multipliers of 
$\mB(\mS\times\mT)$}. It is easy to deduce from the closed graph theorem that  the Schur multipliers  induce bounded operators on $\mB(\mS\times\mT)$. Put
$$
\|\Phi\|_{\fM(\mS\times\mT)}\df
\sup\{\|\Phi\star A\|:~A\in\mB(\mS\times\mT),~  \|A\|_\mB\le1\}.
$$
Hence, by duality
\bay
\label{S1}
\|\Phi\|_{\fM(\mS\times\mT)}=\sup\{\|\Phi\star A\|_{\bS_1}:~A\in\mB(\mS\times\mT),~
\|A\|_{\bS_1}\le1\}.
\ey

It is easy to see that
$$
\|A\|_{\mB(\mS\times\mT)}=\sup\|A\|_{\mB(\mS_0\times\mT_0)}, \quad
\|A\|_{\bS_1(S\times\mT)}=\sup\|A\|_{\bS_1(\mS_0\times\mT_0)}
$$
and
$$
\|\Phi\|_{\fM(\mS\times\mT)}=\sup\|{\bs\f}\|_{\fM(\mS_0\times\mT_0)},
$$
where the suprema are taken over all finite subsets $\mS_0$ and $\mT_0$
of the sets $\mS$ and $\mT$.

Note also that 
$\|\Phi\|_{\ell^\be(\mS\times\mT)}\le\|\Phi\|_{\fM(\mS\times\mT)}$.
It is easy to see that the inequality turns into equality for every matrix $\{\Phi(s,t)\}_{(s,t)\in\mS\times\mT}$ of rank $1$.
There are other classes of matrices, for which this inequality turns into equality.
For example, if each row (or each column) of $\Phi$ has at most one nonzero entry, then  
$\|\Phi\|_{\ell^\be(\mS\times\mT)}=\|\Phi\|_{\fM(\mS\times\mT)}$.

We need one more characteristic of the matrix
$\Phi$. Put
$$
\|\Phi\|_{\fM_0(\mS\times\mT)}\df
\sup\{\|\Phi\star A\|:~A\in\mB(\mS\times\mT),~
\|A\|\le1,~a(t,t)=0~\;\mbox{for}\;~t\in\mS\cap\mT\}.
$$
We denote by $\fM_0(\mS\times\mT)$  the set of matrices
$\Phi=\{\Phi(s,t)\}_{(s,t)\in\mS\times\mT}$ such that $\|\Phi\|_{\fM_0(\mS\times\mT)}<\be$.
Obviously, $\|\Phi\|_{\fM_0(\mS\times\mT)}\le\|\Phi\|_{\fM(\mS\times\mT)}$.
It is easy to see that
$$
\|\Phi\|_{\fM_0(\mS\times\mT)}=\sup\|\Phi\|_{\fM_0(\mS_0\times\mT_0)},
$$
where the supremum is taken over all finite subsets $\mS_0$ and $\mT_0$
of $\mS$ and $\mT$.

Let us also observe that if the matrices $\Phi=\{\Phi(s,t)\}_{(s,t)\in\mS\times\mT}$
and $\Psi=\{\Psi(s,t)\}_{(s,t)\in\mS\times\mT}$ coincide off the ``diagonal''
$\{(t,t):t\in\mS\cap\mT\}$,  then
$$
\|\Phi-\Psi\|_{\fM_0(\mS\times\mT)}=0\quad
\text{and}\quad\|\Phi\|_{\fM_0(\mS\times\mT)}=\|\Psi\|_{\fM_0(\mS\times\mT)}.
$$
Note that
$\|\Phi\|_{\fM_0(\mS\times\mT)}=\|\Phi\|_{\fM(\mS\times\mT)}$ if
$\mS\cap\mT=\varnothing$.

\begin{lem} 
Let $\Phi\in\ell^\be(\mS\times\mT)$, where $\mS$ and $\mT$ are arbitrary sets such that
 $\mS\cap\mT\ne\varnothing$. Then
$$
\max\left\{\|\Phi\|_{\fM_0(\mS\times\mT)},\|\Phi(t,t)\|_{\ell^\be(\mS\cap\mT)}\right\}
\le\|\Phi\|_{\fM(\mS\times\mT)}
\le2\|\Phi\|_{\fM_0(\mS\times\mT)}+\|\Phi(t,t)\|_{\ell^\be(\mS\cap\mT)}.
$$
\end{lem}

\Pf  The first inequality is obvious. Let us prove the second one.
We denote by $\chi$ the characteristic function of the set
$\{(s,t)\in\mS\times\mT:s=t\}$. It is easy to see that $\|\chi\|_{\fM(\mS\times\mT)}=1$,
 whence $\|1-\chi\|_{\fM(\mS\times\mT)}\le\|1\|_{\fM(\mS\times\mT)}+\|\chi\|_{\fM(\mS\times\mT)}=2$.
 Let $A\in\mB(S\times T)$ and $\|A\|\le1$. Then
 $
 \Phi\star A=\Phi\star(1-\chi)\star A+\Phi\star\chi\star A.
 $
 It remains to observe that
 $$
 \|\Phi\star(1-\chi)\star A\|\le\|\Phi\|_{\fM_0(\mS\times\mT)}
 \|(1-\chi)\star A\|\le2\|\Phi\|_{\fM_0(\mS\times\mT)}
 $$
 and
 $$
 \|\Phi\star\chi\star A\|\le\|\Phi\star\chi\|_{\fM(\mS\times\mT)}
 =\|\Phi(t,t)\|_{\ell^\be(\mS\cap\mT)}.\quad \bl
 $$

 \begin{cor}
 \label{32}
 If $\Phi(t,t)=0$ for every $t\in\mS\cap\mT$, then
 $$
\|\Phi\|_{\fM_0(\mS\times\mT)}\le\|\Phi\|_{\fM(\mS\times\mT)}
\le2\|\Phi\|_{\fM_0(\mS\times\mT)}.
$$
 \end{cor}

 \begin{lem}
 \label{33}
 Let $\mS$ be a Hausdorff topological space. Suppose that
 the set $\mS\cap\mT$ has no isolated points in $\mS$.
 Then  $\|\Phi\|_{\fM_0(\mS\times\mT)}=\|\Phi\|_{\fM(\mS\times\mT)}$ for every
 function $\Phi\in\ell^\be(\mS\times\mT)$ continuous in the variable $s\in\mS$.
 \end{lem}

 \Pf It suffices to prove that
 $\|\Phi\|_{\fM(\mS\times\mT)}\le\|\Phi\|_{\fM_0(\mS\times\mT)}$
 or, which is the same, $\|\Phi\|_{\fM(\mS_0\times\mT_0)}
 \le\|\Phi\|_{\fM_0(\mS\times\mT)}$
 for all finite subsets $\mS_0$ and $\mT_0$ of the sets $\mS$ and $\mT$.
 Let us fix finite subsets $\mS_0$ and $\mT_0$ of $\mS$ and $\mT$.
 Clearly, for every positive number $\e$, there exists a perturbation $\widetilde\mS_0$
 of the set $\mS_0$ such that $\widetilde\mS_0\cap\mT_0=\varnothing$ and
 $\|\Phi\|_{\fM(\mS_0\times\mT_0)}<\e+\|\Phi\|_{\fM(\widetilde\mS_0\times\mT_0)}$.
 Consequently,
 $$
 \|\Phi\|_{\fM(\mS_0\times\mT_0)}<\e+\|\Phi\|_{\fM(\widetilde\mS_0\times\mT_0)}=
 \e+\|\Phi\|_{\fM_0(\widetilde\mS_0\times\mT_0)}\le\e+\|\Phi\|_{\fM_0(\mS\times\mT)}
 $$
 for every $\e>0$. $\bl$
 
 \medskip

We are going to consider an analog of the space $\fM_{0}(\mS\times\mT)$
that is defined in terms of the $\bS_1$ norm rather than the operator norm.
We set
$$
\|\Phi\|_{\fM_{0,\bS_1}(\mS\times\mT)}\df
\sup\{\|\Phi\star A\|_{\bS_1}:~A\in\mB(\mS\times\mT),~
\|A\|_{\bS_1}\le1,~a(t,t)=0~\mbox{for}~t\in\mS\cap\mT\}
$$
and $\fM_{0,\bS_1}(\mS\times\mT)\df\{\Phi: \|\Phi\|_{\fM_{0,\bS_1}(\mS\times\mT)}<\be\}$.
We would like to observe that there is no need to define a similar analog of 
$\fM(\mS\times\mT)$ for it coincides with the space
$\fM(\mS\times\mT)$ itself by \rf{S1}.

In the same way as for the operator norm one can prove the following facts.

\begin{lem}
\label{32s1l}
Let $\Phi\in\ell^\be(\mS\times\mT)$, where $\mS$ and $\mT$ are arbitrary
sets such that $\mS\cap\mT\ne\varnothing$. Then
\begin{align*}
\max\left\{\|\Phi\|_{\fM_{0,\bS_1}(\mS\times\mT)},\|\Phi(t,t)\|_{\ell^\be(\mS\cap\mT)}\right\}
&\le\|\Phi\|_{\fM(\mS\times\mT)}\\[.2cm]
&\le2\|\Phi\|_{\fM_{0,\bS_1}(\mS\times\mT)}+\|\Phi(t,t)\|_{\ell^\be(\mS\cap\mT)}.
\end{align*}
If $\Phi(t,t)=0$ for every $t\in\mS\cap\mT$, then
 $$
\|\Phi\|_{\fM_{0,\bS_1}(\mS\times\mT)}\le\|\Phi\|_{\fM(\mS\times\mT)}
\le2\|\Phi\|_{\fM_{0,\bS_1}(\mS\times\mT)}.
$$
\end{lem}

\begin{cor}
\label{32s1}
If $\Phi(t,t)=0$ for every $t\in\mS\cap\mT$, then
$$
\|\Phi\|_{\fM_{0,\bS_1}(\mS\times\mT)}\le\|\Phi\|_{\fM(\mS\times\mT)}
\le2\|\Phi\|_{\fM_{0,\bS_1}(\mS\times\mT)}.
$$
\end{cor}

\begin{lem}
\label{33s1}
Let $\mS$ be a Hausdorff topological space. Suppose that the set 
$\mS\cap\mT$ has no isolated points in $\mS$.
Then  $\|\Phi\|_{\fM_{0,\bS_1}(\mS\times\mT)}=\|\Phi\|_{\fM(\mS\times\mT)}$ for every
function $\Phi$ in $\ell^\be(\mS\times\mT)$ that is continuous in the variable $s\in\mS$.
\end{lem}

\medskip

\section{\bf A description of the discrete Schur multipliers}
\setcounter{equation}{0}
\label{dmsh+}

\medskip

\begin{thm}
\label{41}
Let $\{u_s\}_{s\in\mS}$ and $\{v_t\}_{t\in\mT}$ be families of vectors in a
(not necessarily separable) Hilbert space $\h$
such that $\|u_s\|\cdot\| v_t\|\le 1$ for all $s$ in $\mS$ and
 $t$ in $\mT$. Put
$\Phi(s,t)\df(u_s, v_t)$, $s\in\mS$, $t\in\mT$.  Then $\Phi \in\fM(\mS\times\mT)$
and $\|\Phi\|_{\fM(\mS\times\mT)}\le 1$.
\end{thm}

\Pf By \rf{S1}, it suffices to prove that
$$
\|\{a(s,t)(u_s,v_t)\}\|_{\bS_1} \le\|\{a(s,t)\}\|_{\bS_1}
$$
for every matrix $\{a(s,t)\}$ that induces a trace class operator. Clearly, it suffices to consider the case when $\rank\{a(s,t)\}=1$. Moreover, one can assume that
$\|a(s,t)\|_{\bS_1}=1$. Then $a(s,t)=\a_s\b_t$ for some $\a\in\ell^2(\mS)$
and $\b\in\ell^2(\mT)$ such that $\|\a\|_{\ell^2(\mS)}=\|\b\|_{\ell^2(\mT)}=1$.
Let $\{e_j\}_{j\in J}$ be an orthonormal basis in
$\h$. Put $\hat x(j)\df(x,e_j)$, $j\in J$. Then
\begin{align*}
\|\{\a_s\b_t(u_s,v_t)\}\|_{\bS_1}&\le\sum_{j\in J}\|\{\a_s\b_t\hat u_s(j)\ov{\hat v_t(j)})\}\|_{\bS_1}\\[.2cm]
&=\sum_{j\in J}\|\{\a_s\hat u_s(j)\}\|_{\ell^2(\mS)}\|\{\b_t\ov{\hat v_t(j)}\}\|_{\ell^2(\mT)}\\[.2cm]
&\le\Big(\sum_{j\in J}\|\{\a_s\hat u_s(j)\}\|_{\ell^2(\mS)}^2\Big)^{\frac12}
\Big(\sum_{j\in J}\|\{\b_t\ov{\hat v_t(j)}\}\|_{\ell^2(\mT)}^2\Big)^{\frac12}.
\end{align*}
Clearly,
\begin{align*}
\sum_{j\in J}\|\{\a_s\hat u_s(j)\}\|_{\ell^2(\mS)}^2
&=\sum_{j\in J}\sum_{s\in\mS}|\a_s|^2|\hat u_s(j)|^2\\[.2cm]
&=\sum_{s\in\mS}|\a_s|^2\sum_{j\in J}|\hat u_s(j)|^2
=\sum_{s\in\mS}|\a_s|^2\|u_s\|^2
\le\sup_{s\in\mS}\|u_s\|^2.
\end{align*}
Similarly,
$$
\sum_{j\in J}\|\{\b_t\ov{\hat v_t(j)}\}\|_{\ell^2(\mT)}^2\le\sup_{t\in\mT}\|v_t\|^2.
$$
Hence,
$$
\|\{\a_s\b_t(u_s,v_t)\}\|_{\bS_1}
\le\sup_{s\in\mS}\|u_s\|\sup_{t\in\mT}\|u_t\|\le 1. \quad  \bl
$$

It is very nontrivial that the converse also holds, see 
Theorem 5.1 of the monograph \cite{Pi} and also \cite{Pi2}.
We state this result without a proof.

\begin{thm}
\label{Pi}
Let $\Phi\df\{\Phi(s,t)\}$ be a Schur multiplier of $\mB(\mS\times\mT)$ and let
$\|\Phi\|_\fM\le1$. Then there are families of vectors
$\{u_s\}_{s\in\mS}$ and $\{v_t\}_{t\in\mT}$  in a
(not necessarily separable) Hilbert space $\h$
such that $\|u_s\|\le1$ for every $s\in\mS$, $\| v_t\|\le 1$  for every
 $t\in\mT$ and
$$
\Phi(s,t)=(u_s, v_t),\quad s\in\mS,\quad t\in\mT.
$$
\end{thm}

{\bf A remark to Theorem \ref{Pi}.}  {\it In this theorem we can additionally require that the linear spans of the family $\{u_s\}_{s\in\mS}$ and of the 
family $\{ v_t\}_{t\in\mT}$ are dense in $\h$}. Indeed,
let $\h_1$ be the closed linear span of the family $\{ v_t\}_{t\in\mT}$ and let $P_1$ be the orthogonal projection onto $\h_1$. Then
$\{P_1u_s\}_{s\in\mS}$ and $\{ v_t\}_{t\in\mT}$ are families in the Hilbert space $\h_1$ such that $\Phi(s,t)=(P_1u_s,v_t)$ for all $(s,t)\in\mS\times\mT$.
Let now $\h_2$ be the closed linear span of the family $\{P_1u_s\}_{s\in\mS}$ and let $P_2$ be the orthogonal projection onto $\h_2$.
Then $\{P_1u_s\}_{s\in\mS}$ and $\{P_2 v_t\}_{t\in\mT}$ are families in $\h_2$ such that
 $\Phi(s,t)=(P_1u_s,P_2 v_t)$ for $(s,t)\in\mS\times\mT$. It is clear that the linear spans of the families $\{P_1u_s\}_{s\in\mS}$ and $\{P_2 v_t\}_{t\in\mT}$ are dense in 
 $\h_2$. 

\medskip

The following theorem is contained in the results of \cite{KS3} and \cite{A4}.

\begin{thm}
\label{topolx}
Let $\Phi\in\fM(\mS\times\mT)$, where $\mS$ and $\mT$ are topological spaces.
Suppose that $\Phi$ is continuous in each variable. Then there are families
$\{u_s\}_{s\in\mS}$ and $\{ v_t\}_{t\in\mT}$ in a (not necessarily separable) Hilbert space $\h$ such that

{\rm a)} the linear span of $\{u_s\}_{s\in\mS}$ is dense in $\h$;

{\rm b)} the linear span of $\{ v_t\}_{t\in\mT}$ is dense in $\h$;

{\rm c)} the map $s\mapsto u_s$ is weakly continuous;

{\rm d)} the map $t\mapsto v_t$ is weakly continuous;

{\rm e)} $\|u_s\|^2\le\|\Phi\|_{\fM(\mS\times\mT)}$ for every $s\in\mS$;

{\rm f)} $\| v_t\|^2\le\|\Phi\|_{\fM(\mS\times\mT)}$ for every $t\in\mT$;

{\rm g)} $\Phi (s,t)=(u_s, v_t)$ for every $(s,t)\in\mS\times\mT$.
\end{thm}

\Pf By Theorem \ref{Pi} and by the remark to it, there are
families $\{u_s\}_{s\in\mS}$ and $\{ v_t\}_{t\in\mT}$
satisfying conditions a), b), e), f) and g).
Clearly, the function $s\mapsto(u_s,h)$ is continuous for $h= v_t$, where
$t\in\mT$. Thus, the function $s\mapsto(u_s,h)$ is continuous for every $h\in\h$ by b).
Therefore, the map $s\mapsto u_s$ is weakly continuous. Similarly, one can deduce from a) the weak continuity of the map $t\mapsto v_t$. $\bl$

\medskip

{\bf Remark.}  If at least one of the spaces $\mS$ and $\mT$ is separable, then the space $\h$ is also separable.
Indeed, it suffices to observe that if, for example, $\mS$ is separable, the the closed linear span of the family $\{u_s\}_{s\in S}$ is separable.

\medskip

This remark allows us to obtain the following assertion:

\begin{thm}
\label{npkp}
Let $\Phi\in\fM(\mS\times \mT)$, where $ \mathcal S$ and $ \mathcal T$ are
topological spaces such that at least one of them is separable.
Suppose that $\Phi$ is continuous in each variable. 
Then there exist a sequence $\{\f_n\}_{n\ge1}$ of continuous functions on $ \mathcal S$
and a sequence $\{\psi_n\}_{n\ge1}$ of continuous functions on $ \mathcal T$ such that 
$$
\sum_{n=1}^\be|\f_n(s)|^2\le\|\Phi\|_{\fM(\mathcal \mathcal S\times \mathcal T)},\quad
\sum_{n=1}^\be|\psi_n(t)|^2\le\|\Phi\|_{\fM(\mathcal \mathcal S\times \mathcal T)}
$$
and
$$
\sum_{n=1}^\be\f_n(s)\psi_n(t)=\Phi(s,t)
\quad\mbox{for all}~s\in\mathcal S~\mbox{~and~}~t\in\mathcal T.
$$
\end{thm}

\Pf  Let $\{u_s\}_{s\in \mathcal S}$ and $\{ v_t\}_{t\in \mathcal T}$ be families in a Hilbert space $\h$ whose existence is guaranteed by Theorem \ref{topolx}. It follows from the remark to that theorem that the space $\h$
is separable. Let $\{e_n\}_{n=1}^N$ be an orthonormal basis in $\h$, where 
$0\le N\le\be$. It remains to define $\f_n(s)\df(u_s,e_n)$ and
$\psi_n(t)\df(e_n,v_t)$; if $N<\be$, then $\f_n(s)\df\psi_n(t)\df0$
for $n>N$. $\bl$

\medskip

{\bf Definition.}
A map $g$ from a topological space $\mT$ to a Hilbert space $\h$ is called {\it weakly Borel measurable} if
the function $t\mapsto(g(t),u)$ is Borel measurable on $T$ for every $u$ in $\h$.

\medskip

It is easy to see that it suffices to verify that the measurability of the function 
$t\mapsto(g(t),u)$ only for vectors $u$ in a subset of $\h$ whose linear span is dense in 
$\h$.

\begin{thm}
\label{topolbor}
Let $\mS$ and $\mT$ be topological spaces and let $\Phi\in\fM(\mS\times\mT)$.
Suppose that $\Phi$ is Borel measurable in each variable. Then there exist families $\{u_s\}_{s\in\mS}$ and $\{ v_t\}_{t\in\mT}$ in a (not necessarily separable) Hilbert space 
$\h$ such that

{\rm a)} the linear span of $\{u_s\}_{s\in\mS}$ is dense in $\h$;

{\rm b)} the linear span of $\{ v_t\}_{t\in\mT}$ is dense in $\h$;

{\rm c)} the map $s\mapsto u_s$ is weakly Borel measurable;

{\rm d)} the map $t\mapsto v_t$ is weakly Borel measurable;

{\rm e)} $\|u_s\|^2\le\|\Phi\|_{\fM(\mS\times\mT)}$ for every $s\in\mS$;

{\rm f)} $\| v_t\|^2\le\|\Phi\|_{\fM(\mS\times\mT)}$ for every $t\in\mT$;

{\rm g)} $\Phi(s,t)=(u_s, v_t)$ for all $(s,t)\in\mS\times\mT$.
\end{thm}

The proof of Theorem \ref{topolbor} practically repeats word for word the proof of Theorem \ref{topolx}.

\begin{thm}
\label{topolborizm}
Let $\mS$ and $\mT$ be topological spaces and let 
$\Phi$ be a Borel function in $\fM(\mS\times\mT)$. 
Suppose that $\mu$ and $\nu$ are Borel
$\s$-finite measures on $\mS$ and $\mT$. Then there exist sequences
$\{\f_k\}_{k\ge1}$ and $\{\psi_k\}_{k\ge1}$ such that

{\rm a)} $\f_k\in L^\be(\mu)$ and $\psi_k\in L^\be(\nu)$ for every $k\ge1$;

{\rm b)} $\sum_{k=1}^\be|\f_k(s)|^2\le\|\Phi\|_{\fM(S\times T)}$ for $\mu$-almost all $s\in\mS$;

{\rm c)} $\sum_{k=1}^\be|\psi_k(t)|^2\le\|\Phi\|_{\fM(S\times T)}$ for $\nu$-almost all $t\in\mT$;

{\rm d)} $\Phi(s,t)=\sum_{k=1}^\be\f_k(s)\psi_k(t)$ for $\mu\otimes\nu$-almost all $(s,t)\in\mS\times\mT$.
\end{thm}

\Pf Clearly, we can assume that $\|\Phi\|_{\fM(S\times T)}=1$. Let 
$\{u_s\}_{s\in\mS}$ and $\{ v_t\}_{t\in\mT}$ be families in a Hilbert space $\h$
whose existence is guaranteed by Theorem \ref{topolbor}. Let $\{e_j\}_{j\in J}$ be an orthonormal basis in $\h$. Put
$\f_j(s)\df(u_s,e_j)$ and $\psi_j(t)\df(e_j,v_t)$. Then the $\f_j$ and $\psi_j$ are Borel measurable,
$$
\sum\limits_{j\in J}|\f_j(s)|^2\le1\quad\mbox{for every}~s\in\mS,
\quad\sum\limits_{j\in J}|\psi_j(t)|^2\le1\quad\mbox{for every}~t\in\mT
$$
and
$$
\Phi(s,t)=\sum_{j\in J}\f_j(s)\psi_j(t)\quad\mbox{for every}~(s,t)\in\mS\times\mT.
$$

This immediately completes the proof of the theorem in the case when $J$ is at most countable. To consider the case of an arbitrary set $J$, we set
$\Psi(s,t)\df\sum\limits_{j\in J}|\f_j(s)|\cdot|\psi_j(t)|$. By the Cauchy--Bunyakovsky inequality,
$$
\Psi(s,t)\le\Big(\sum\limits_{j\in J}|\f_j(s)|^2\Big)^{1/2}\Big(\sum\limits_{j\in J}|\psi_j(t)|^2\Big)^{1/2}
\le1.
$$
We may assume that $\mu$ and $\nu$ are probability measures. 
Put $J_s\df\{j\in J: \f_j(s)\ne0\}$, where $s\in\mS$. Note that
$J_s$ is at most countable for every $s\in\mS$ because $\sum_{j\in J}|\f_j(s)|^2\le1$.
It is easy to see that that for every $s$ in $\mS$, 
\begin{align*}
\sum_{j\in J}|\f_j(s)|\int_\mT|\psi_j(t)|\,d\nu(t)&=\sum_{j\in J_s}|\f_j(s)|\int_\mT|\psi_j(t)|\,d\nu(t)\\[.2cm]
&=\int_\mT\Big(\sum_{j\in J_s}|\f_j(s)|\cdot|\psi_j(t)|\Big)\,d\nu(t)=\int_\mT\Psi(s,t)\,d\nu(t).
\end{align*}
To integrate in $s$, we consider the at most countable set
$J_\flat\df\Big\{j\in J:\lb\int_T|\psi_j(t)|\,d\nu(t)\ne0\Big\}$. Then
\begin{align*}
\sum_{j\in J}\int_\mS|\f_j(s)|\,d\mu(s)\int_\mT|\psi_j(t)|\,d\nu(t)&=
\sum_{j\in J_\flat}\int_\mS|\f_j(s)|\,d\mu(s)\int_\mT|\psi_j(t)|\,d\nu(t)\\[.2cm]
&=
\int_\mS\left(\int_\mT\Psi(s,t)\,d\nu(t)\right)\,d\mu(s).
\end{align*}
Now it is clear that
$$
\int_\mS\Big(\int_\mT\sum_{j\in J\setminus J_\flat}|\f_j(s)|\cdot|\psi_j(t)|\,d\nu(t)\Big)\,d\mu(s)=0.
$$
This together with the inequality
$$
\Big|\Phi(s,t)-\sum_{j\in J_\flat}\f_j(s)\psi_j(t)\Big|\le\sum_{j\in J\setminus J_\flat}|\f_j(s)|\cdot|\psi_j(t)|
$$
implies that
$$
\sum_{j\in J_\flat}\f_j(s)\psi_j(t)=\Phi(s,t)
$$
for $\mu\otimes\nu$-almost all $(s,t)\in\mS\times\mT$. $\bl$

\medskip

Consider now some examples of Schur multipliers.
Let $\mM(\T^2)$ be the space if complex Borel measure on the two-dimensional 
torus $\T^2$ with norm $\|\mu\|_{\mM(\T^2)}\df|\mu|(\T)$.

\medskip

{\bf Example.}  Let $\mu\in\mM(\T^2)$.
Then $\{\widehat\mu(m,n)\}\in\fM(\Z^2)$ and $\|\widehat\mu\|_{\fM(\Z\times\Z)}\le\|\mu\|$.

\medskip

This fact is an obvious consequence of Theorem \ref{41}. It is clear that not every 
Schur multiplier ${\bs a}\in\fM(\Z\times\Z)$ can be represented as
${\bs a}=\widehat\mu$, where $\mu\in\mM(\T^2)$.
Consider, for example, the case when the matrix
${\bs a}=\{a_{mn}\}_{m,n\in\Z}$ consists of the same columns (or rows). To be definite, suppose that 
$a_{mn}=t_n$ for all $m,n\in\Z$. Then
${\bs a}\in\fM(\Z\times\Z)$ if and only if
${\bs a}\in\ell^\be(\Z\times\Z)$ and $\|{\bs a}\|_{\fM(\Z\times\Z)}=\|{\bs a}\|_{\ell^\be(\Z\times\Z)}=\|\{t_n\}\|_{\ell^\be}$.
Obviously, not all such matrices ${\bs a}$ with bounded entries
can be represented as ${\bs a}=\widehat\mu$, where $\mu\in\mM(\T^2)$.

On the other hand, if we assume that a matrix ${\bs a}=\{a_{mn}\}_{m,n\in\Z}$
is a {\it Laurent matrix}, i.e., $a_{mn}=t_{m-n}$, the situation considerably changes.

\begin{thm}
\label{Boch}
Let $A=\{a_{mn}\}_{m,n\in\Z}$ a Laurent matrix.
Then $A\in\fM(\Z^2)$ if and only if  $a_{mn}=\widehat\mu(m-n)$
for a measure $\mu$ in $\mM(\T)$
and $\|{\bs a}\|_{\fM(\Z^2)}=\|\mu\|_{\mM(\T)}$.
\end{thm}

All these results can be generalized to locally compact abelian groups. Herewith in the case of a nondiscrete abelian group $G$ in the statement of the analog of
Theorem \ref{Boch} we have to require that the functions in question are continuous.

\begin{thm}
\label{BochR}
Let $h$ be a continuous function on $\R$.
Then the matrix $A=\{h(s-t)\}_{s,t\in\R}$ belongs
$\fM(\R\times\R)$ if and only if there exists a complex Borel measure $\mu$ such that 
$h=\F\mu$. Moreover,
$\|A\|_{\fM(\R\times\R)}=\|\mu\|_{\mM(\R)}$.
\end{thm}

\medskip

\section{\bf Double operator integrals}
\setcounter{equation}{0}
\label{mSid}

\medskip

Double operator integrals are expressions of the form
\bay
\label{doi}
\int_\mS\int_\mT\Phi(s,t)\,dE_1(s)T\,dE_2(t),
\ey
where $E_1$ and $E_2$ are spectral measures on a separable Hilbert space $\h$, 
$\Phi$ is a bounded measurable function and $T$ is a bounded linear operator on $\h$.

Double operator integrals appeared in the paper \cite{DK}. In the papers
\cite{BS1}, \cite{BS2} and \cite{BS3} Birman and Solomyak created a harmonious theory of double operator integrals. The idea of Birman and Solomyak was to define first double operator integrals of the form \rf{doi} for arbitrary bounded measurable functions $\Phi$ and for operators $T$ of Hilbert--Schmidt class $\bS_2$. For this purpose they introduced a spectral measure $\E$ that takes values in the set of orthogonal projections on the Hilbert space $\bS_2$ and is defined by
$$
\E(\L\times\D)T=E_1(\L)TE_2(\D),\quad T\in\bS_2,
$$
where $\L$ and $\D$ are measurable subsets of $\mS$ and $\mT$. 
It is clear that left multiplication by $E_1(\L)$ commutes with right multiplication by
$E_2(\D)$. In \cite{BS4}  it was shown that $\E$ extends to a spectral measure on
 $\mS\times\mT$. In this situation the double operator integral \rf{doi} is defined by
$$
\int_\mS\int_\mT\Phi(s,t)\,d E_1(s)T\,dE_2(t)\df
\left(\,\,\int_{\mS\times\mT}\Phi\,d\E\right)T.
$$
It follows immediately from this definition that
$$
\left\|\int_\mS\int_\mT\Phi(s,t)\,dE_1(s)T\,dE_2(t)\right\|_{\bS_2}
\le\|\Phi\|_{L^\be}\|T\|_{\bS_2}.
$$

If a function $\Phi$ possesses the property
$$
T\in\bS_1\quad\Longrightarrow\quad\int_\mS\int_\mT\Phi(s,t)\,d E_1(s)T\,dE_2(t)\in\bS_1,
$$
$\Phi$ is said to be a {\it Schur multiplier of $\bS_1$ with respect to the spectral measures $E_1$ and $E_2$}.

To define double operator integrals \rf{doi} for bounded operators $T$,
we consider the transformer
$$
Q\mapsto\int_{\mT}\int_{\mS}\Phi(t,s)\,d E_2(t)\,Q\,dE_1(s),\quad Q\in\bS_1,
$$
and assume that the function $(y,x)\mapsto\Phi(y,x)$ is a Schur multiplier of 
$\bS_1$ with respect to $E_2$ and $E_1$.
In this case the transformer
\bay
\label{tra}
T\mapsto\int_\mS\int_\mT\Phi(s,t)\,d E_1(s)T\,dE_2(t),\quad T\in \bS_2,
\ey
extends by duality to a bounded linear transformer on the space of bounded linear operators on $\h$. In this case $\Phi$ is said to be a {\it Schur multiplier (with respect to $E_1$ and $E_2$) of the space of bounded linear operators}.
We denote the space of such Schur multipliers by $\fM(E_1,E_2)$.
The norm of $\Phi$ in $\fM(E_1,E_2)$ is defined as the norm of the transformer
\rf{tra} on the space of bounded linear operators.

It is easy to see that if a function $\Phi$ on $\mS\times\mT$ belongs to the  {\it projective tensor product} $L^\be(E_1)\hat\otimes L^\be(E_2)$ of the spaces $L^\be(E_1)$ 
and $L^\be(E_2)$ (i.e., $\Phi$ admits a representation
$$
\Phi(s,t)=\sum_{n\ge0}\f_n(s)\psi_n(t),\quad\mbox{where}\quad
\sum_{n\ge0}\|\f_n\|_{L^\be(E_1)}\|\psi_n\|_{L^\be(E_2)}<\be),
$$
then $\Phi\in\fM(E_1,E_2)$.
For such functions $\Phi$ we have
$$
\int_\mS\int_\mT\Phi(s,t)\,dE_1(s)T\,dE_2(t)=
\sum_{n\ge0}\left(\,\int_\mS\f_n\,dE_1\right)T\left(\,\int_\mT\psi_n\,dE_2\right).
$$

More generally, $\Phi\in\fM(E_1,E_2)$ if $\Phi$
belongs to the {\it integral projective tensor product} 
$L^\be(E_1)\hat\otimes_{\rm i}L^\be(E_2)$ of the spaces $L^\be(E_1)$ and $L^\be(E_2)$, i.e.,  $\Phi$ admits a representation
\bay
\label{ipt}
\Phi(s,t)=\int_\O \f(s,w)\psi(t,w)\,d\s(w),
\ey
where $(\O,\s)$ is a space with a $\s$-finite measure, $\f$ is a measurable function on $\mS\times \O$,\lb
$\psi$ is a measurable function on $\mT\times \O$ and
$$
\int_\O\|\f(\cdot,w)\|_{L^\be(E_1)}\|\psi(\cdot,w)\|_{L^\be(E_2)}\,d\s(w)<\be.
$$

It turns out that all Schur multipliers can be obtained in this way (see Theorem \ref{tomSc}).

Another sufficient condition for a function to be a Schur multiplier can be stated in terms of the Haagerup tensor product  $L^\be(E_1)\!\otimes_{\rm h}\!L^\be(E_2)$, which is defined as the space of functions $\Phi$ of the form
\bay
\label{FiH}
\Phi(s,t)=\sum_{n\ge0}\f_n(s)\psi_n(t),
\ey
where $\{\f_n\}_{n\ge0}\in L_{E_1}^\be(\ell^2)$ and $\{\psi_n\}_{n\ge0}\in L_{E_2}^\be(\ell^2)$. Put
$$
\|\Phi\|_{L^\be(E_1)\otimes_{\rm h}\!L^\be(E_2)}\df\inf
\Big\|\sum_{n\ge0}|\f_n|^2\Big\|_{L^\be(E_1)}^{1/2}
\Big\|\sum_{n\ge0}|\psi_n|^2\Big\|_{L^\be(E_2)}^{1/2},
$$
where the infimum is taken over all representations of $\Phi$ in the form \rf{FiH}. 
It is easy to verify that if $\Phi\in L^\be(E_1)\!\otimes_{\rm h}\!L^\be(E_2)$, then $\Phi\in\fM(E_1,E_2)$ and
\bay
\label{Haagrazl}
\iint\Phi(s,t)\,dE_1(s)T\,dE_2(t)=
\sum_{n\ge0}\Big(\int\f_n\,dE_1\Big)T\Big(\int\psi_n\,dE_2\Big),
\ey
the series on the right being convergent in the weak operator topology and
$$
\|\Phi\|_{\fM(E_1,E_2)}\le\|\Phi\|_{L^\be(E_1)\otimes_{\rm h}L^\be(E_2)}.
$$

As we can see from the following theorem, the condition
$\Phi\in L^\be(E_1)\!\otimes_{\rm h}\!L^\be(E_2)$ is not only sufficient, but also necessary. 

\begin{thm}
\label{tomSc} 
Let $\Phi$ be a measurable function on
$\mS\times\mT$ and  $\mu$ and $\nu$ be positive $\s$-finite measures on $\mS$ and $\mT$ they are mutually absolutely continuous with respect to $E_1$ and $E_2$. The following statements are equivalent:

{\rm (a)} $\Phi\in\fM(E_1,E_2)$;

{\rm (b)} $\Phi\in L^\be(E_1)\hat\otimes_{\rm i}L^\be(E_2)$;

{\rm (c)} $\Phi\in L^\be(E_1)\!\otimes_{\rm h}\!L^\be(E_2)$;

{\rm (d)} there exist a $\s$-finite measure $\s$ on a set $\O$, measurable functions $\f$ on $\mS\times\O$ and $\psi$ on $\mT\times\O$ such that
{\em\rf{ipt}} holds and
\bay
\label{bs}
\left\|\left(\int_\O|\f(\cdot,w)|^2\,d\s(w)\right)^{1/2}\right\|_{L^\be(E_1)}
\left\|\left(\int_\O|\psi(\cdot,w)|^2\,d\s(w)\right)^{1/2}\right\|_{L^\be(E_2)}<\be;
\ey

{\rm (e)} if an integral operator
$f\mapsto\int k(x,y)f(y)\,d\nu(y)$ from $L^2(\nu)$ to $L^2(\mu)$ belongs to $\bS_1$, then the integral operator $f\mapsto\int\Phi(x,y)k(x,y)f(y)\,d\nu(y)$ belongs to the same class.
\end{thm}

The implications (d)$\imp$(a)$\eq$(e) were established in \cite{BS3}.
In the case of matrix Schur multipliers the implication (a)$\imp$(b) was proved in  \cite{Ben}. We refer the reader to \cite{Pe1} for the proof of the equivalence of (a), (b) and (d). The equivalence of (c) and (d) can be proved elementarily.

It is easy to see that conditions (a) -- (e) are also equivalent to the fact that $\Phi$ is a Schur multiplier of $\bS_1$. 

Note that one can also define double operator integrals of the form \rf{doi}  in the case when $E_1$ and $E_2$ are spectral measures on different Hilbert spaces and $T$ is an operator from one Hilbert space to another one.

\medskip

{\bf Remark.} It follows easily from Theorems \ref{topolborizm} and \ref{tomSc} that if $\mS$ and $\mT$ are topological spaces and $\Phi$ is a Borel function on $\mS\times\mT$ of class $\fM(\mS\times\mT)$ (i.e., $\Phi$ is a discrete Schur multiplier), then $\Phi\in\fM(E_1,E_2)$ for all Borel spectral measures $E_1$ and $E_2$ on $\mS$ and $\mT$.

\medskip

Double operator integrals can also be defined with respect to semispectral measures. Recall that a {\it semispectral measure} $\mE$ 
on a measurable space $(\X,{\frak B})$ is a map on the $\s$-algebra ${\frak B}$ with values in the set of bounded linear operators on a Hilbert space $\h$ that is countably additive in the strong operator topology and such that
$$
\mE(\D)\ge\0,\quad\D\in{\frak B},\quad\mE(\varnothing)=\0
\quad\mbox{and}\quad\mE(\X)=I.
$$
By Naimark's theorem \cite{Nai}, each semispectral measure $\mE$ has a
{\it spectral dilation}, i.e., a spectral measure $E$ on the same measurable space
$(\X,{\frak B})$ that takes values in the set of orthogonal projections on a Hilbert space $\K$ containing $\h$, and such that$$
\mE(\D)=P_\h E(\D)\big|\h,\quad\D\in{\frak B},
$$
where $P_\h$ is the orthogonal projection on $\K$ onto $\h$. Such a spectral dilation can be chosen {\it minimal} in the sense that
$$
\K=\clos\spn\{E(\D)\h:~\D\in{\frak B}\}.
$$ 
Note that it was shown in \cite{MM} that if $E$ is a minimal spectral dilation of a semispectral measure $\E$, then $E$ and $\E$ are mutually absolutely continuous.

Integrals with respect to semispectral  measures are defined in the following way:
$$
\int_\X \f(x)\,d\mE(x)=P_\h\left(\int_\X \f(x)\,d E(x)\right)\Big|\h,\quad
\f\in L^\be(\mE)\df L^\be(E).
$$

If $\E_1$ and $\E_2$ are semispectral measures on
$(\X_1,{\frak B}_1)$ and $(\X_2,{\frak B}_2)$, $E_1$
and $E_2$ are their spectral dilations on Hilbert spaces $\K_1$ and $\K_2$, and a function 
$\Phi$ on $\X_1\times\X_2$ satisfies the equivalent statements of Theorem \ref{tomSc}, then the double operator integral with respect to  $\mE_1$ and $\mE_2$ is defined by 
$$
\int_{\!\X_1}\!\!\int_{\!\X_2}\Phi(x_1,x_2)\,d\E_1(x_1)Q\,d\E_2(x_2)=
P^{[1]}_\h\int_{\!\X_1}\!\!\int_{\!\X_2}\Phi(x_1,x_2)\,dE_1(x_1)\big(QP^{[2]}_{\h}\big)\,dE_2(x_2)\Big|\h
$$
for an arbitrary bounded linear operator $Q$ on $\h$. 
Here $P^{[1]}_\h$ and $P^{[2]}_\h$ are the orthogonal projections from $\K_1$ and $\K_2$
onto $\h$.
If $\Phi\in L^\be(E_1)\!\otimes_{\rm h}\!L^\be(E_2)$, then
the following equality holds:
\bay
\label{doipolcpm}
\iint\Phi(x_1,x_2)\,d\mE_1(x_1)T\,d\mE_2(x_2)=
\sum_{n\ge0}\Big(\int\f_n\,d\mE_1\Big)T\Big(\int\psi_n\,d\mE_2\Big),
\ey
where $T\in\mB(\h)$, and $\f_n$ and $\psi_n$ are functions in the representation \rf{FiH}.

Double operator integrals with respect to semispectral measures were introduced in  \cite{Pe*}, see also \cite{Pe9}.

\medskip

\begin{center}
\bf\Large Chapter III
\end{center}

\begin{center}
\bf\Large 
Operator Lipschitz function on subsets of the plane
\end{center}
\label{OL2per}

\addtocontents{toc}{\vspace*{.3cm}\hspace*{-.35cm}\textbf{{\bf Chapter III}. Operator Lipschitz function on subsets of the plane}\hfill\pageref{OL2per}}

\renewcommand{\thesection}{3.\arabic{section}}
\setcounter{section}{0}

\medskip

In this chapter we study operator Lipschitz and commutator Lipschitz functions on closed subsets of the complex plane. A significant role will be played by Schur multipliers. 
We offer two methods of obtaining difference and commutator estimates.
The first method uses discrete Schur multipliers and approximation by operators with finite spectrum. The second method is based on double operator integrals.

\

\section{\bf Operator Lipschitz and commutator Lipschitz functions on closed subsets of the plane}
\setcounter{equation}{0}
\label{oplip}

\

We define here the classes of operator Lipschitz functions and of commutator Lipschitz functions on closed subsets of the plane. We will see that unlike in the case of functions on the line, these two classes by no means have to coincide.
When defining these classes, we consider only bounded operators. In the next section 
we will see that if we admit not necessarily bounded operators, we obtain the same classes of functions.

Let $\fF$ be a nonempty subset of the complex plane $\C$. We denote by $\Li(\fF)$  the space of functions $f:\fF\to\C$ satisfying the {\it Lipschitz condition}:
\bay
\label{lipsc}
|f(z)-f(w)|\le C|z-w|,\quad z,w\in\C.
\ey
The smallest constant $C\ge0$ satisfying \rf{lipsc} is denoted by 
$\|f\|_{\Li(\fF)}=\|f\|_{\Li}$. Put $\|f\|_{\Li}\df\be$ if $f\not\in\Li(\fF)$.

Usually we require the closedness of $\fF$.

It follows easily from the spectral theory for pairs  {\it commuting} normal operators that the inequality
\bay
\label{1comlip}
\|f(N_1)-f(N_2)\|\le \|f\|_{\Li(\fF)}\|N_1-N_2\|
\ey
holds for arbitrary commuting normal operators $N_1$ and $N_2$ whose spectra are contained in $\fF$.

A complex continuous function $f$ on a nonempty closed set $\fF$,
$\fF\subset\C$, will be called {\it operator Lipschitz} if there exists a positive number
$C$ such that
\bay
\label{lipid}
\|f(N_1)-f(N_2)\|\le C\|N_1-N_2\|
\ey
for arbitrary normal operators $N_1$ and $N_2$
with spectra in $\fF$. We denote the space of operator Lipschitz functions on $\fF$ by 
$\OL(\fF)$. The smallest constant $C$ satisfying \rf{lipid} is denoted by 
$\|f\|_{\OL(\fF)}=\|f\|_{\OL}$. Put $\|f\|_{\OL}=\be$ if $f\not\in\OL(\fF)$.

If a function $f$ is defined on a bigger set ${\frak G}\supset\fF$, we will usually write for brevity $f\in\OL(\fF)$ and $\|f\|_{\OL(\fF)}$ rather than $f|\fF\in\OL(\fF)$ and $\|f|\fF\|_{\OL(\fF)}$. We will also use the same convention for other function spaces.

It is easy to see that $\OL(\fF)\subset\Li(\fF)$ and 
$\|f\|_{\Li(\fF)}\le\|f\|_{\OL(\fF)}$
for every $f\in\OL(\fF)$. We will see in \S\:\ref{orav} that the equality $\OL(\fF)=\Li(\fF)$ holds only for finite sets $\fF$.

Note that if $f\in\OL(\fF)$ and $\|f\|_{\OL}\le1$,  then
\bay
\label{unit}
\|f(N_1)U-Uf(N_2)\|\le\|N_1U-UN_2\|
\ey
for all unitary operators $U$ and all
normal operators $N_1$ and $N_2$ such that $\s(N_1),\s(N_2)\subset\fF$.
To see this, it suffices to apply inequality \rf{lipid}
with $C=1$ to normal operators $U^*N_1U$ and $N_2$. Conversely, if \rf{unit} holds
for all unitary operators $U$ and all normal operators
 $N_1$  and $N_2$ such that $N_1=N_2$ and $\s(N_1)=\s(N_2)\subset\fF$, then
$f\in\OL(\fF)$ and $\|f\|_{\OL(\fF)}\le1$.
Indeed, applying inequality \rf{unit} to the operators
$$
\mathcal N_1=\mathcal N_2= \left(\begin{array}{cc}
    N_1 & \0 \\
    \0 & N_2 \\
  \end{array}\right)
\quad\text{and}\quad
\mathcal U= \left(\begin{array}{cc}
    \0 & I \\
    I & \0 \\
  \end{array}\right),
$$
we obtain
$$
\|f(N_1)-f(N_2)\|\le\|N_1-N_2\|.
$$
Note that in this reasoning we dealt only with self-adjoint unitary operators
$\mathcal U$, i.e., normal operators $\mathcal U$ such that $\mathcal U^2=I$
or, which is the same, unitary operators with spectra in $\{-1,1\}$.

\begin{thm}
\label{Nlip0}
Let $f$ be a continuous function on a subset $\fF$ of $\C$.
The following statements are equivalent:

{\em(a)} $\|f(N_1)-f(N_2)\|\le\|N_1-N_2\|$ for all normal operators $N_1$ and $N_2$ such that
$\s(N_1),\s(N_2)\subset\fF$;

{\em(b)} $\|f(N_1)U-Uf(N_2)\|\le\|N_1U-UN_2\|$ for all unitary operators $U$
and all normal operators $N_1$ and $N_2$ such that $\s(N_1),\s(N_2)\subset\fF$;

{\em(c)} $\|f(N)U-Uf(N)\|\le\|NU-UN\|$ for all self-adjoint unitary operators $U$ and all
normal operators $N$ such that $\s(N)\subset\fF$;

{\em(d)} $\|f(N)A-Af(N)\|\le\|NA-AN\|$ for all self-adjoint operators $A$ and all normal operators $N$ such that $\s(N)\subset\fF$.
\end{thm}

\Pf The equivalence of (a), (b) and (c) was proved as a matter of fact
before the statement of the theorem. The implication (d)$\Longrightarrow$(c) is obvious.
It remains to prove that (c) implies (d). Denote by $\frak X$ the set of operators $R$
such that $\|f(N)R-Rf(N)\|\le\|NR-RN\|$ for all normal operators 
$N$ with spectrum in $\fF$. It is clear that the set $\frak X$ is closed in the norm and
$\a U+\b I\in\frak X$ for every unitary operator $U$ and for all $\a,\b\in\C$. to prove that an arbitrary self-adjoint operator $A$ belongs to $\frak X$,
it suffices to observe that the operator $(I-\e{\rm i}A)(I+\e{\rm i}A)^{-1}$ is
unitary for all $\e$ in $(-\|A\|^{-1},\|A\|^{-1})$ and
$$
A=\lim_{\e\to0}\frac1{2\e{\rm i}}\left(I-(I-\e{\rm i}A)(I+\e{\rm i}A)^{-1}\right).\quad\bl
$$

\medskip

{\bf Remark.}  A unitary operator $U$ is self-adjoint if and only if it can be represented in the form $U=2P-I$, where $P$ is an orthogonal projection.
Statement (c) of Theorem \ref{Nlip0} can be rewritten in the following:

{\it$\|f(N)P-Pf(N)\|\le\|NP-PN\|$ for all orthogonal projections $P$ and all 
normal operators $N$ such that $\s(N)\subset\fF$}.

\begin{thm}
\label{Nlip}
Let $f$ be a continuous function on a closed subset $\fF$ of $\C$.
The following statements are equivalent:

{\em(a)} $\|f(N_1)-f(N_2)\|\le\|N_1-N_2\|$ for all normal operators $N_1$ and $N_2$ such that
$\s(N_1),\s(N_2)\subset\fF$;

{\em(b)} $\|f(N)R-Rf(N)\|\le\max\big\{\|NR-RN\|,\|N^*R-RN^*\|\big\}$ for all operators
$R\in\mB(\h)$ and all normal operators $N$ such that $\s(N)\subset\fF$;

{\em(c)} $\|f(N_1)R-Rf(N_2)\|\le\max\big\{\|N_1R-RN_2\|,\|N_1^*R-RN_2^*\|\big\}$ for all 
$R\in\mB(\h)$ and all normal operators $N_1$ and $N_2$ such that 
$\s(N_1),\s(N_2)\subset\fF$.
\end{thm}

\Pf Let us first prove the implication (a)$\Longrightarrow$(b). Suppose that
 (a) holds. Then it follows from Theorem \ref{Nlip0} that $\|f(N)A-Af(N)\|\le\|NA-AN\|$ for all self-adjoint operators $A$ and all normal operators
$N$ such that $\s(N)\subset\fF$. Applying this assertion to the normal operator
$\left(
  \begin{array}{cc}
     N & \0 \\
     \0 & N \\
  \end{array}
\right)$
and the self-adjoint operator
$\left(
  \begin{array}{cc}
     \0 & R \\
      R^* & \0 \\
  \end{array}
\right)$, we obtain
$$
\max\big\{\|f(N)R-Rf(N)\|,\|f(N)R^*-R^*f(N)\|\big\}
\le\max\big\{\|NR-RN)\|,\|NR^*-R^*N\|\big\}
$$
which implies (b).

Applying (b) to the normal operator
$\left(
  \begin{array}{cc}
     N_1 & \0 \\
     \0 & N_2 \\
  \end{array}
\right)$
and the operator
$\left(
  \begin{array}{cc}
     \0 & R \\
     \0 & \0 \\
  \end{array}
\right)$, we obtain (c). 
The implication (c)$\Longrightarrow$(a) is obvious. $\bl$

\medskip

In the proofs of Theorems \ref{Nlip0} and \ref{Nlip} we have used the standard technique of passing to block matrix operators; it allows one in certain cases to proceed from one operator to a pair of operators. This technique will be useful in what follows.
Kittaneh \cite{Ki} calls it the Berberian trick keeping apparently in mind the 
paper \cite{Ber} by Berberian.

Theorems \ref{Nlip0} and \ref{Nlip} are contained in Theorem 3.1 of
 \cite{AP6}, but to a certain extent in a certain form they can be  
extracted from the paper \cite{KS2}, in which together with the operator norm
arbitrary symmetric norms are also considered. 

Note that the equality
$\|N_1^*R-RN_2^*\|=\|N_1R-RN_2\|,$
and so the equality
$$
\max\big(\|N_1R-RN_2\|,\|N_1^*R-RN_2^*\|\big)=\|N_1R-RN_2\|,
$$
holds in each of the following special cases:

1) the operators $N_1$ and $N_2$ are self-adjoint (this is the case if $\fF\subset\R$);

2) the operators $N_1$ and $N_2$ are unitary (this is the case if $\fF\subset\T$);

3) $R$ is self-adjoint and $N_1=N_2$;

4) $R$ is a unitary operator.

A complex function $f$
continuous on a closed set $\fF$, $\fF\subset\C$, is called
{ \it commutator Lipschitz} if there is a number $C\ge0$
such that
\bay
\label{cld}
\|f(N)R-Rf(N)\|\le C\|NR-RN\|
\ey
for every $R\in\mB(\h)$ and every normal operator $N$ with spectrum in $\fF$.
We denote the set of commutator Lipschitz functions on
$\fF$ by $\CL(\fF)$. The smallest constant $C$, satisfying \rf{cld} is denoted by $\|f\|_{\CL(\fF)}=\|f\|_{\CL}$. Put $\|f\|_{\CL(\fF)}=\be$ if $f\not\in\CL(\fF)$.

\begin{thm}
\label{Clip}
Let $f$ be a continuous function on a closed subset $\fF$ of  $\C$.
The following statements are equivalent:

{\em(a)} $\|f(N)R-Rf(N)\|\le\|NR-RN\|$ for every $R\in\mB(\h)$ and for every normal operator
$N$ such that $\s(N)\subset\fF$;

{\em(b)} $\|f(N_1)R-Rf(N_2)\|\le\|N_1R-RN_2\|$ for every $R\in\mB(\h)$ and for all normal operators
$N_1$ and $N_2$ such that
$\s(N_1),\:\s(N_2)\subset\fF$;

{\em(c)} $\|f(N_1)A-Af(N_2)\|\le\|N_1A-AN_2\|$ for every self-adjoint operator $A$
and for all normal operators $N_1$ and $N_2$ such that $\s(N_1),\:\s(N_2)\subset\fF$.
\end{thm}

\Pf To prove the implication (a)$\Longrightarrow$(b), it suffices to apply (a)
to the normal operator
$\left(
  \begin{array}{cc}
     N_1 & \0 \\
     \0 & N_2 \\
  \end{array}
\right)$
and the operator
$\left(
  \begin{array}{cc}
     \0 & R \\
     \0 & \0 \\
  \end{array}
\right)$.
The implication (b)$\Longrightarrow$(c) is evident.
It remains to prove that (c) implies (a).
 Applying (c) to $N_1=U^{*}NU$ and $N_2=N$, where $U$ is a unitary operator, we obtain
 $$
 \|f(N)UA-UAf(N)\|=\|f(U^{*}NU)A-Af(N)\|\le\|NUA-UAN\|
 $$
 for arbitrary self-adjoint operator $A$, unitary operator $U$
 and normal operator $N$ such that $\s(N)\subset\fF$.
 Note that if (a) is satisfied for an operator $R\in\mB(\h)$,
 then it is also satisfied for the operator $R+\l I$, where $\l\in\C$. Thus, we may assume that $R$ is invertible. Then applying the polar decomposition to the invertible operator $R$, we obtain
 $R=UA$, where $U$ is a unitary operator and $A$ is a (positive) self-adjoint operator. $\bl$

It follows immediately from Theorems \ref{Nlip0} and \ref{Clip} that
$\CL(\fF)\subset\OL(\fF)$ and $\|f\|_{\OL(\fF)}\le\|f\|_{\CL(\fF)}$
for every $f$ in $\CL(\fF)$.

\medskip

{\bf Remark.} 
In statements  (b) of Theorem \ref{Nlip0},  (c) of Theorem \ref{Nlip} and 
(b) of Theorem \ref{Clip}) we may assume that the normal operators $N_1$ and $N_2$ act on different Hilbert spaces (herewith the unitary operator 
$U$  can act from one Hilbert space to another one).   
This can be seen from the proofs. We give here, as an illustration, a relevant
reformulation of statement (b) of Theorem  \ref{Clip}:

\medskip

{\em $\|f(N_1)R-Rf(N_2)\|\le\|N_1R-RN_2\|$ for all operators  $R\in\mB(\mathcal \h_2,\h_1)$  and all normal operators $N_1$ and $N_2$ acting on Hilbert spaces $\h_1$
and $\h_2$ such that $\s(N_1),\:\s(N_2)\subset\fF$.}

\medskip

Analogs of Theorems \ref{Nlip0}, \ref{Nlip}  and \ref{Clip} with appropriate modifications hold for symmetrically normed ideals with practically the same proofs. 
We consider here only trace class $\bS_1$.
With each closed set $\fF$, $\fF\subset\C$, we associate the space of
{\it trace class Lipschitz} (or {\it$\bS_1$-Lipschitz}) functions $\OL_{\bS_1}(\fF)$ and 
the space of {\it trace class commutator Lipschitz}
(or {\it$\bS_1$-commutator Lipschitz}) functions 
$\CL_{\bS_1}(\fF)$. To define the spaces $\OL_{\bS_1}(\fF)$
and $\CL_{\bS_1}(\fF)$, 
we have to replace in \rf{lipid} and \rf{cld} the operator norm with the trace norm.

The corresponding ''united'' analog of Theorems \ref{Nlip0} and \ref{Nlip} for the trace norm can be stated in the following way. 

\begin{thm}
\label{yaNlip01}
Let $f$ be a continuous function on a closed subset $\fF$ of the complex plane $\C$.
The following statements are equivalent:

{\em(a)} $\|f(N_1)-f(N_2)\|_{\bS_1}\le\|N_1-N_2\|_{\bS_1}$ for all normal operators $N_1$ and $N_2$ such that
$\s(N_1),\:\s(N_2)\subset\fF$;

{\em(b)} $\|f(N_1)U-Uf(N_2)\|_{\bS_1}\le\|N_1U-UN_2\|_{\bS_1}$ for all unitary operators $U$ and all normal operators $N_1$ and $N_2$ such that $\s(N_1),\:\s(N_2)\subset\fF$;

{\em(c)} $\|f(N)U-Uf(N)\|_{\bS_1}\le\|NU-UN\|_{\bS_1}$ for all self-adjoint unitary operators $U$ and all normal operators $N$ such that $\s(N)\subset\fF$;

{\em(d)} $\|f(N)A-Af(N)\|_{\bS_1}\le\|NA-AN\|_{\bS_1}$ for all self-adjoint operators $A$ and all normal operators $N$ such that $\s(N)\subset\fF$;

{\em(e)} $\|f(N)R-Rf(N)\|_{\bS_1}+\|\ov f(N)R-R\ov f(N)\|_{\bS_1}\le\|NR-RN\|_{\bS_1}+\|N^*R-RN^*\|_{\bS_1}$ for all
$R\in\mB(\h)$ and all normal operators $N$ such that $\s(N)\subset\fF$;

{\em(f)} $\|f(N_1)R-Rf(N_2)\|_{\bS_1}+\|\ov f(N_1)R-R\ov f(N_2)\|_{\bS_1}
\le\|N_1R-RN_2\|_{\bS_1}+\lb\|N_1^*R-RN_2^*\|_{\bS_1}$ for all $R\in\mB(\h)$
and all normal operators $N_1$ and $N_2$ such that $\s(N_1),\:\s(N_2)\subset\fF$.
\end{thm}

Let us state now the analog of Theorem \ref{Clip}.

\begin{thm}
\label{ClipS1}
Let $f$ be a continuous function on a closed subset $\fF$ of $\C$.
The following statements are equivalent:

{\em(a)} $\|f(N)R-Rf(N)\|_{\bS_1}\le\|NR-RN\|_{\bS_1}$ for all 
$R\in\mB(\h)$ and all normal operators $N$ such that $\s(N)\subset\fF$;

{\em(b)} $\|f(N_1)R-Rf(N_2)\|_{\bS_1}\le\|N_1R-RN_2\|_{\bS_1}$ for all $R\in\mB(\h)$ and all normal operators
$N_1$ and $N_2$ such that
$\s(N_1),\:\s(N_2)\subset\fF$;

{\em(c)} $\|f(N_1)A-Af(N_2)\|_{\bS_1}\le\|N_1A-AN_2\|_{\bS_1}$ for all self-adjoint operators $A$
and all normal operators $N_1$ and $N_2$ such that $\s(N_1),\:\s(N_2)\subset\fF$.
\end{thm}

Note that one can reformulate Theorem \ref{yaNlip01} for self-adjoint operators
in the following way:

\begin{thm}
\label{yaNlip01A}
Let $f$ be a real continuous function on a closed subset $\fF$ of $\R$.
The following statements are equivalent:

{\em(a)} $\|f(A)-f(B)\|_{\bS_1}\le\|A-B\|_{\bS_1}$ for all self-adjoint operators $A$ and $B$ such that
$\s(A),\s(B)\subset\fF$;

{\em(b)} $\|f(A)U-Uf(B)\|_{\bS_1}\le\|AU-UB\|_{\bS_1}$ for all unitary operators $U$
and all self-adjoint operators $A$ and $B$ such that $\s(A),\:\s(B)\subset\fF$;

{\em(c)} $\|f(A)U-Uf(A)\|_{\bS_1}\le\|AU-UA\|_{\bS_1}$ for all self-adjoint unitary operators $U$ and all self-adjoint operators $A$ such that $\s(A)\subset\fF$;

{\em(d)} $\|f(A)R-Rf(A)\|_{\bS_1}\le\|AR-RA\|_{\bS_1}$ for all self-adjoint operators  $A$ and $R$ such that $\s(A)\subset\fF$;

{\em(e)} $\|f(A)R-Rf(A)\|_{\bS_1}\le\|AR-RA\|_{\bS_1}$ for all
$R\in\mB(\h)$ and all self-adjoint operators $A$ such that $\s(A)\subset\fF$;

{\em(f)} $\|f(A)R-Rf(B)\|_{\bS_1}
\le\|AR-RB\|_{\bS_1}$ for all $R\in\mB(\h)$
and all self-adjoint operators $A$ and $B$ such that $\s(A),\:\s(B)\subset\fF$.
\end{thm}

\begin{cor}
\label{sledoveshch}
If $f$ is a continuous real function on a closed subset
$\fF$ of the real line $\R$, then
$\|f\|_{\OL_{\bS_1}(\fF)}=\|f\|_{\CL_{\bS_1}(\fF)}$.
\end{cor}

It follows that for a complex continuous function $f$ the following inequalities hold:
$$
\|f\|_{\OL_{\bS_1}(\fF)}\le\|f\|_{\CL_{\bS_1}(\fF)}\le2\|f\|_{\OL_{\bS_1}(\fF)}.
$$

The same also can be said in the case of unitary operators $N_1$ and $N_2$,
i.e., in the case when $\fF$ is contained in $\T$.

Clearly, $\ov z\in\OL(\fF)$ for every closed set $\fF$ in $\C$ and
$\|\ov z\|_{\OL(\fF)}=1$ if $\fF$ has at least two points.

\medskip

{\bf Definition.}
A closed subset $\fF$ of $\C$ is called a
{\it Fuglede set} if $\CL(\fF)=\OL(\fF)$. 

\medskip

This notion was introduced by Kissin and Shulman in \cite{KS2}. 

Johnson and Williams \cite{JW} proved that each function $f\in\CL(\fF)$
is differentiable in the complex sense at each nonisolated point of 
$\fF$, see Theorems \ref{olm} and \ref{pro} below. Note that $\ov z\in\OL(\fF)$.
Therefore, a Fuglede set cannot have interior points and even cannot contain
two intersecting intervals not contained in the same straight line. 
Kissin and Shulman proved in \cite{KS2} that each compact curve of class 
$C^2$ is a Fuglede set.

The following theorem is substantially contained in Proposition 4.5 of \cite{KS2}.

\begin{thm} 
\label{Fug}
A closed subset $\fF$ of $\C$ is a Fuglede set if and only if
$\ov z\in\CL(\fF)$. If $\ov z\in\CL(\fF)$, then
$\|f\|_{\CL(\fF)}\le\|\ov z\|_{\CL(\fF)}\|f\|_{\OL(\fF)}$ for every $f\in\OL(\fF)$.
\end{thm}

The theorem is a straightforward consequence of Theorems \ref{Nlip} and \ref{Clip}.

\begin{cor} 
Let $\fF$ be a closed subset of $\C$. Then the equality $\CL(\fF)=\OL(\fF)$ holds together with the equality of the seminorms $\|\cdot\|_{\CL(\fF)}=\|\cdot\|_{\OL(\fF)}$ if and only if $\|\ov z\|_{\CL(\fF)}\le1$.
\end{cor}

Note that $\|\ov z\|_{\CL(\fF)}\ge\|\ov z\|_{\OL(\fF)}=1$ if
$\fF$ has at least two points.
This, the condition $\|\ov z\|_{\CL(\fF)}\le1$ can be replaced with the condition
$\|\ov z\|_{\CL(\fF)}=1$ whenever $\fF$ has at least two points.

\begin{thm}
 \label{su}
If a closed subset $\fF$ of $\C$ is contained in a straight line or in a circle, then
$\fF$ is a Fuglede set and $\|\cdot\|_{\CL(\fF)}=\|\cdot\|_{\OL(\fF)}$.
\end{thm}

\Pf It is easy to see that the (semi)norms in $\CL(\fF)$ and $\OL(\fF)$ coincide if and only if
\bay
\label{Km}
\|N^*R-RN^*\|=\|NR-RN\|
\ey
for every normal operator $N$ with spectrum in $\fF$ and for every bounded operator $R$.
As we have observed above, see special cases 1) and 2) in front of Theorem \ref{Clip}, this equality obviously holds for both self-adjoint and unitary operators  $N$.
This proves the theorem in the two following special cases: $\fF\subset\R$ and $\fF\subset\T$.
The general case reduces to these special cases with the help of affine transforms of
the complex plane. $\bl$

\medskip

Kamowitz \cite{Ka} proved that for a given operator $N$, \rf{Km} holds for all bounded operators  $R$ if and only if $N$ is a normal operator whose spectrum is contained in a circle or in a line. It follows from this result of Kamowitz that
Theorem \ref{su} has a converse. In other words, the equality 
$\|\cdot\|_{\CL(\fF)}=\|\cdot\|_{\OL(\fF)}$
holds if and only if $\fF$ is contained in a circle or in a line.

Let $\fF_1$ and $\fF_2$ be nonempty closed subsets of $\C$.  Denote by $\CL(\fF_1,\fF_2)$
the space of continuous functions $f$ on $\fF=\fF_1\cup\fF_2$,
for which there exists $C\ge0$ such that
\bay
\label{cld2}
\|f(N_1)R-Rf(N_2)\|\le C\|N_1R-RN_2\|
\ey
for all $R\in\mB(\h)$ and all normal operators $N_1$ and $N_2$ with spectra
in $\fF_1$ and $\fF_2$. We denote by $\|f\|_{\CL(\fF_1,\fF_2)}$ the smallest constant
$C$ satisfying \rf{cld2}
Put $\|f\|_{\CL(\fF_1,\fF_2)}=\be$ if $f\not\in\CL(\fF_1,\fF_2)$.

Passing to the adjoint operators, we see that \rf{cld2}
is equivalent to the condition
$
\|R^*\ov f(N_1)-\ov f(N_2)R^*\|\le C\|R^*N_1^*-N_2^*R^*\|.
$
It follows that $f\in\CL(\fF_1,\fF_2)$ if and only if 
$\ov f(\ov z)\in\CL(\ov\fF_2,\ov\fF_1)$ and
$\|\ov f(\ov z)\|_{\CL(\ov\fF_2,\ov\fF_1)}=\|f\|_{\CL(\fF_1,\fF_2)}$,
where for a subset $\fF$ of $\C$, we denote by $\ov\fF$ the set $\{\ov\z:~\z\in\fF\}$.

If we rewrite \rf{cld2} in terms of matrices and consider the transposed matrices,
we find that $\CL(\fF_1,\fF_2)=\CL(\fF_2,\fF_1)$ and
$\|\cdot\|_{\CL(\fF_1,\fF_2)}=\|\cdot\|_{\CL(\fF_2,\fF_1)}$.

\begin{thm}
\label{28}
Let $f$ be a continuous function on the union
$\fF_1\cup\fF_2$ of closed subsets $\fF_1$ and $\fF_2$
of $\C$. The following statements are equivalent:

{\em(a)} $\|f(N_1)R-Rf(N_2)\|\le\|N_1R-RN_2\|$ for all $R\in\mB(\h)$ and all
normal operators $N_1$ and $N_2$ such that $\s(N_1)\subset\fF_1$ and
$\s(N_2)\subset\fF_2$;

{\em(b)} $\|f(N_1)R-Rf(N_2)\|\le\|N_1R-RN_2\|$ for an arbitrary operator
$R$ from a Hilbert space $\h_2$
to a Hilbert space $\h_1$ and arbitrary normal operators
$N_1$ and $N_2$ on $\h_1$
and $\h_2$ such that $\s(N_1)\subset\fF_1$ and $\s(N_2)\subset\fF_2$;

{\em(c)} statement {\em(b)} holds under the additional assumption that the normal operators $N_1$ and $N_2$ have simple spectra;

{\em(d)} $\|f(N_1)A-Af(N_2)\|\le\|N_1A-AN_2\|$ for all self-adjoint operators $A$
and all normal operators $N_1$ and $N_2$ such that $\s(N_1)\subset\fF_1$ and
$\s(N_2)\subset\fF_2$.
\end{thm}

\Pf 
The implications (b)$\Longrightarrow$(a) and  (b)$\Longrightarrow$(c) are trivial.
Let us prove that (a)$\Longrightarrow$(b).
If the spaces $\h_1$ and $\h_2$ are isomorphic, then there exists a unitary operator
 $U:\h_1\to\h_2$.
The operator $RU$ and the normal operators $N_1$ and $U^*N_2U$ are 
operators on the same Hilbert space $\h_1$, and so
$$
\|f(N_1)(RU)-(RU)f(U^*N_2U)\|\le\|N_1(RU)-(RU)U^*N_2U\|
$$
by (a) which immediately implies the desired estimate.
To reduce the general case to the special case considered above,
we introduce the operators $\mathcal R\df\bigoplus_{j\ge1}R$, 
$\mathcal N_1\df\bigoplus_{j\ge1}N_1$
and $\mathcal N_2\df\bigoplus_{j\ge1}N_2$. It is easy to see that the inequality
$\|f(N_1)R-Rf(N_2)\|\le\lb\|N_1R-RN_2\|$
is equivalent to the inequality
$
\|f(\mathcal N_1)\mathcal R-\mathcal Rf(\mathcal N_2)\|\le
\|\mathcal N_1\mathcal R-\mathcal R\mathcal N_2\|.
$

Let us prove now that (c) implies (b).
Assume the contrary. Then there exist $R\in\mB(\h_2,\h_1)$ and normal operators $N_1$ and $N_2$ on $\h_1$ and $\h_2$ such that
$\s(N_1)\subset\fF_1$, $\s(N_2)\subset\fF_2$, $\|N_1R-RN_2\|=1$ and $\|f(N_1)R-Rf(N_2)\|>1$.
Thus, there are vectors $u_0\in\h_2$ and $v_0\in\h_1$ such that $\|u_0\|=1$,
$\|v_0\|=1$ and
$\left|\left((f(N_1)R-Rf(N_2))u_0,v_0\right)\right|>1$. Let $\h_1^0$ and
$\h_2^0$ be the smaller reducing subspaces of 
$N_1$ and $N_2$ that contain $v_0$ and $u_0$.
Let $P$ be the orthogonal projection onto $\h_1^0$
and let $Q$ be the orthogonal projection onto $\h_2^0$.
Note that $\|f(N_1)PRQ-PRQf(N_2)\|>1$ for $((f(N_1)PRQ-PRQf(N_2))u_0,v_0)=
((f(N_1)R-Rf(N_2))u_0,v_0)$. Moreover,
$\|N_1PRQ-PRQN_2\|=\|P(N_1R-RN_2)Q\|\le1$.
Put $N_1^0\df N|{\h_1^0}$ and $N_2^0\df N|{\h_2^0}$.
Then the operators $N_1^0$ and $N_2^0$ can be considered as normal operators on
$\h_1^0$ and $\h_2^0$.
Clearly, $N_1^0$ and $N_2^0$ are normal operators with simple spectra.
To get a contradiction, it suffices to observe that
$\|f(N_1^0)PRQ-PRQf(N_2^0)\|>1$ and $\|N_1^0PRQ-PRQN_2^0\|\le1$.

The implication (a)$\Longrightarrow$(d) is trivial.
It remains to prove that (d) implies (a).
Applying (d) to the normal operators $U^{*}N_1U$ and $N_2$, where $U$ is a unitary operator, we obtain
 $$
 \|f(N_1)UA-UAf(N_2)\|\le\|N_1UA-UAN_2\|
 $$
 for every self-adjoint operator $A$, every unitary operator $U$
 and for arbitrary normal operators $N_1$ and $N_2$ such that $\s(N_1)\subset\fF_1$ and
 $\s(N_2)\subset\fF_2$. With the help of polar decomposition this implies
 (a) for invertible operators $R$. 
 Therefore, (a) holds for 
 the operators $R$ that belong to the closure of the set of invertible operators in the operator norm. It remains to observe that in the general case the block operator 
 $\mathcal R=\left(\begin{array}{cc}
    R & \0 \\
    \0 & \0 \\
  \end{array}\right)$
  on $\h\oplus\h$ can obviously be approximated with an arbitrary accuracy in the operator norm by invertible operators in $\mB(\h\oplus\h)$.
  One can proceed from the operator $\mathcal R$ to the operator $R$ by using the Berberian trick, which was discussed after the proof of Theorem \ref{Nlip}.
 $\bl$
 
We need the following well-known elementary result.

\begin{lem}
\label{apr}
Let $N$ be a bounded normal operator. Suppose that the subset  $\L$ of $\C$
is an $\e$-net of the spectrum $\s(N)$ of $N$, i.e., for each $\z\in\s(N)$,
there is $\l\in\L$ such that $|\l-\z|<\e$. Then there exists a normal operator
$N_0$ such that $NN_0=N_0N$, $\|N-N_0\|<\e$ and $\s(N_0)$ is a finite subset of
$\L$.
\end{lem}

\Pf Since the spectrum of $N$ is compact, there exists a finite $\e$-net $\L_0$
of $\s(N)$ such that $\L_0\subset \L$. Then we can find a Borel function 
$\eta:\s(N)\to\L_0$ such that $\sup\{|z-\eta(z)|:z\in\s(N)\}<\e$.
It remains to put $N_0\df\eta(N)$. $\bl$

\medskip

It follows easily from this lemma and from inequality \rf{1comlip} that if inequality
\rf{lipid} holds for all normal operators $N_1$ and  $N_2$ with finite spectra
in  $\fF$, then $f\in\OL(\fF)$ and $\|f\|_{\OL(\fF)}\le C$.

In other words, for every continuous function $f$ on a closed subset $\fF$
of $\fF\subset\C$, the following equality holds:
\bay
\label{finite}
\|f\|_{\OL(\frak F)}=\sup\big\{\|f\|_{\OL(\L)}:\L\subset\frak F,\,\,\,
\L\,\,\text{is finite}\big\}.
\ey
Moreover, 
\bay
\label{finite1}
\|f\|_{\OL(\frak F)}=\sup\big\{\|f\|_{\OL(\L)}:\L\subset\fF_0,\,\,\,
\L\,\,\text{is finite}\big\},
\ey
where $\fF_0$ is a dense subset of $\fF$.

Similar equalities also hold for the commutator Lipschitz seminorms.

Hence, we would obtain nothing new if we 
tried to define the spaces $\OL(\fF)$ and $\CL(\fF)$ for an arbitrary subset
$\fF$ of $\C$.

To be definite, we dwell on the space $\OL(\fF)$ (the space $\CL(\fF)$  can be considered in a similar way).
We say that an arbitrary function $f:\fF\to\C$ belongs to
$\OL(\fF)$ if there is $C\ge0$ such that 
inequality \rf{lipid} holds for all normal operators 
$N_1$ and $N_2$ with finite spectra in $\fF$.
Note that since the spectra are finite, we can define
$f(N_1)$ and $f(N_2)$ for an arbitrary function $f$.
Clearly, $\OL(\fF)\subset\Li(\fF)$. Thus, each such function
$f$ admits a Lipschitz extension to the closure $\clos\fF$ of $\fF$.
It is easy to see from \rf{finite1} that this extension belongs to $\OL(\clos\fF)$ and its $\OL$-seminorm does not change.
Therefore, the space $\OL(\fF)$ can be identified in a natural way with the space
 $\OL(\clos\fF)$. 

Taking this remark into account, we will write for brevity
$\OL(\dd)$, $\OL(\C_+)$, $\CL(\dd)$ and $\CL(\C_+)$ rather than
$\OL(\clos\dd)$, $\OL(\clos\C_+)$, $\CL(\clos\dd)$ and $\CL(\clos\C_+)$.

\medskip

\section{\bf Bounded and unbounded normal operators}
\setcounter{equation}{0}
\label{ogrneogr}

\medskip

We prove in this section certain auxiliary results that imply that in the definitions of operator and commutator Lipschitz (as well as operator H\"older) functions 
one can consider only bounded normal operators or admit unbounded operators. In both cases we get the same classes with the same norms.

Let $N_1$ and $N_2$ be not necessarily bounded normal operators on Hilbert spaces 
$\h_1$ and $\h_2$ with domains $\mD_{N_1}$ and $\mD_{N_2}$.
Let $R$ be a bounded operator from $\h_2$ to $\h_1$.
We say that that $N_1R-RN_2$ is a {\it bounded operator} if
$R(\mD_{N_2})\subset \mD_{N_1}$ and $\|N_1Ru-RN_2u\|\le\const\|u\|$ for every $u\in \mD_{N_2}$.
Then there exists a unique bounded operator $K$ such that 
$Ku=N_1Ru-RN_2u$ for every $u\in \cd_{N_2}$. In this case we write $K=N_1R-RN_2$.
Thus, $N_1R-RN_2$ is a bounded operator if and only if
\bay
\label{MN}
\big|(Ru,{N_1}^*v)-(N_2u,R^*v)\big|\le\const\|u\|\cdot\|v\|
\ey
for all $u\in \mD_{N_2}$ and $v\in \mD_{N_1^*}=\mD_{N_1}$. It is easy to see that
$N_1R-RN_2$ is a bounded operator if and only if $N_2^*R^*-R^*N_1^*$ is a bounded operator. Herewith
$(N_1R-RN_2)^*=-(N_2^*R^*-R^*N_1^*)$.
In particular, we write $N_1R=RN_2$ if $R(\mD_{N_2})\subset\mD_{N_1}$ and $N_1Ru=RN_2u$ for every $u\in\mD_{N_2}$.
We write $\|N_1R-RN_2\|=\be$ if $N_1R-RN_2$ is not a bounded operator.

\medskip

{\bf Remark.} Let $N_1$ and $N_2$ be normal operators. Suppose that $N_1^*$ is the closure of ${N_1}_\flat$ and $N_2$ is the closure of ${N_2}_\sharp$.
Then if \rf{MN} holds for all $u\in \mD_{{N_2}_\sharp}$ and $v\in\mD_{{N_1}_\flat}$,
then it holds for all $u\in \mD_{N_2}$ and $v\in \mD_{N_1}$.

\begin{thm}
\label{anbnrn}
Let $N_1$ and $N_2$ be normal operators on Hilbert spaces $\mathscr H_1$
and $\mathscr H_2$ and let $R$ be a bounded operator from $\h_2$ to $\h_1$. Then there exist sequences  
$\{N_{1,n}\}_{n\ge1}$ and $\{N_{2,n}\}_{n\ge1}$ of bounded normal operators on Hilbert spaces $\h_{1,n}$
and $\h_{2,n}$ and a sequence of bounded operators $\{R_n\}_{n\ge1}$
from $\h_{2,n}$ to $\h_{1,n}$ such that 

{\em(a)} the sequence $\{\|R_n\|\}_{n\ge1}$ is nondecreasing and
$
\lim_{n\to\be}\|R_n\|=\|R\|;
$

{\em(b)} $\s(N_{1,n})\subset\s(N_1)$ and $\s(N_{2,n})\subset\s(N_2)$ for every $n\ge1$;

{\em(c)} for every continuous function $f$ on $\s(N_1)\cup\s(N_2)$,
the sequence 
$\big\{\big\|f(N_{1,n})R_n-R_nf(N_{2,n})\big\|\big\}_{n\ge1}$ is nondecreasing
and
$$
\lim_{n\to\be}\big\|f(N_{1,n})R_n-R_nf(N_{2,n})\big\|=\|f(N_1)R-Rf(N_2)\|;
$$

{\em(d)} for every continuous function $f$ on $\s(N_1)\cup\s(N_2)$ such that
$\|f(N_1)R-Rf(N_2)\|<\be$ and for every natural number $j$, the sequence 
$\big\{s_j\big(f(N_{1,n})R_n-R_nf(N_{2,n})\big)\big\}_{n\ge0}$ of singular values  is nondecreasing and
$$
\lim_{n\to\be}s_j\big(f(N_{1,n})R_n-R_nf(N_{2,n})\big)=s_j\big(f(N_1)R-Rf(N_2)\big).
$$
\end{thm}

\Pf Without loss of generality we may assume that $0\in\s(N_1)\cup\s(N_2)$.
Put $P_{1,n}\df E_{N_1}\big(\{|\l|\le n\})$ and $P_{2,n}\df E_{N_2}\big(\{|\l|\le n\})$,
where $E_{N_1}$ and $E_{N_2}$ are the spectral measures of the normal operators $N_1$ and $N_2$. 
Put 
$$
\widetilde N_{1,n}\df P_{1,n}N_1=N_1P_{1,n}=P_{1,n}N_1P_{1,n},\quad
\widetilde N_{2,n}\df P_{2,n}N_2=N_2P_{2,n}=P_{2,n}N_2P_{2,n},
$$
$$
\mathscr H_{1,n}\df P_{1,n}\mathscr H_1\quad\mbox{and}\quad\mathscr H_{2,n}\df P_{2,n}\mathscr H_2.
$$
Clearly, $\widetilde N_{1,n}$ and $\widetilde N_{2,n}$ are bounded normal operators on $\mathscr H_1$ and $\mathscr H_2$, and $\mathscr H_{1,n}$
and $\mathscr H_{2,n}$ are reducing subspaces of $\widetilde N_{1,n}$ and $\widetilde N_{2,n}$. 

Put $N_{1,n}\df\widetilde N_{1,n}\big|\mathscr H_{1,n}$ and $N_{2,n}\df\widetilde N_{2,n}\big|\mathscr H_{2,n}$.
Then $N_{1,n}$ and $N_{2,n}$ are normal operators on $\mathscr H_{1,n}$ and $\mathscr H_{2,n}$. The operator $R_n$ in $\mB(\h_{2,n},\h_{1,n})$ is defined by $R_nu\df P_{1,n}Ru=P_{1,n}R P_{2,n}u$ for $u\in\mathscr H_{2,n}$.
Then (a) and (b) are obvious. To prove the remaining assertions,
it suffices to observe that
$$
P_{1,n}\big(f(N_1)R-Rf(N_2)\big)P_{2,n}u=\Big(f(N_{1,n})R_n-R_nf(N_{2,n})\Big)u
$$
for every $u$ in $\mathscr H_{2,n}$. $\bl$

\medskip

\section{\bf Divided difference and commutator Lipschitzness}
\setcounter{equation}{0}
\label{razdrazn}

\medskip

With each function $f$ on a closed set $\fF$, $\fF\subset\C$,
we associate the function
\lb$\dg_0 f:\fF\times\fF\to\C$,
\bay
\label{razdra0}
(\dg_0 f)(z,w)\df\left\{\begin{array}{ll}\dfrac{f(z)-f(w)}{z-w},&\text {if}\,\,\,\,z\ne w,\\[.2cm]
0,&\text{if}\,\,\,\,z=w.
\end{array}\right.
\ey
If $\fF$ has no isolated points and for each point
$z$ in $\fF$, there exists a finite derivative $f'(z)$
in the complex sense, then we can define the {\it divided difference}
$\dg f:\fF\times\fF\to\C$ by
$$
(\dg f)(z,w)\df\left\{\begin{array}{ll}\dfrac{f(z)-f(w)}{z-w},&\text {if}\,\,\,\,z\ne w,\\[.2cm]
f'(z),&\text{if}\,\,\,\,z=w.
\end{array}\right.
$$

\begin{thm}
 \label{f1f2}
Let $f$ be a continuous function on the union
$\fF_1\cup\fF_2$ of closed subsets $\fF_1$ and $\fF_2$
of $\C$. Then $f\in\CL(\fF_1,\fF_2)$ if and only if 
$\dg_0 f\in\fM(\fF_1\times\fF_2)$. Moreover,
$$
\|f\|_{\CL(\fF_1,\fF_2)}=\|\dg_0 f\|_{\fM_0(\fF_1\times\fF_2)}\le\|\dg_0 f\|_{\fM(\fF_1\times\fF_2)}\le2\|f\|_{\CL(\fF_1,\fF_2)}.
$$
\end{thm}

\Pf We prove only the equality $\|f\|_{\CL(\fF_1,\fF_2)}=\|\dg_0 f\|_{\fM_0(\fF_1\times\fF_2)}$ because everything else follows from Corollary \ref{32}. Consider first the case of finite sets $\fF_1$ and $\fF_2$.
Let $N_1$ and $N_2$ be normal operators such that $\s(N_1)\subset\fF_1$
and $\s(N_2)\subset\fF_2$. By Theorem \ref{28}, we can assume that $N_1$ and $N_2$
have simple spectra. Then there exist orthonormal bases
$\{u_\l\}_{\l\in\s(N_1)}$ and $\{v_\mu\}_{\mu\in\s(N_2)}$ in $\h_1$ and $\h_2$
such that $N_1u_\l=\l u_\l$ for every 
$\l\in\s(N_1)$ and $N_2v_\mu=\mu v_\mu$ for every $\mu\in\s(N_2)$.  With each operator $X:\h_2\to\h_1$ we associate the matrix
$\{(Xv_\mu,u_\l)\}_{(\l,\mu)\in\s(N_1)\times\s(N_2)}$. We have
$$
\left((N_1R-RN_2)v_\mu,u_\l\right)=
\big(Rv_\mu,N_1^*u_\l)-(RN_2v_\mu,u_\l)=
(\l-\mu)(Rv_\mu,u_\l).
$$
Similarly,
$$
\left((f(N_1)R-Rf(N_2))v_\mu,u_\l\right)=(f(\l)-f(\mu))(Rv_\mu,u_\l).
$$
Clearly,
$$
\{(f(\l)-f(\mu))(Rv_\mu,u_\l)\}
=\{(\dg_0 f)(\l,\mu)\}\star
\{(\l-\mu)(Rv_\mu,u_\l)\}.
$$
Note that the matrix $\{a_{\l\mu}\}_{(\l,\mu)\in\s(N_1)\times\s(N_2)}$ can be represented in the form
$$
\{a_{\l\mu}\}_{(\l,\mu)\in\s(N_1)\times\s(N_2)}=
\{(\l-\mu)(Rv_\mu,u_\l)\}_{(\l,\mu)\in\s(N_1)\times\s(N_2)},
$$
where $R$ is an operator from $\h_2$ to $\h_1$, if and only if $a_{\l\mu}=0$ for $\l=\mu$.
The equality $\|f\|_{\CL(\fF_1,\fF_2)}=\|\dg_0 f\|_{\fM_0(\fF_1\times\fF_2)}$
for finite sets $\fF_1$ and $\fF_2$ is now obvious.

The inequality $\|f\|_{\CL(\fF_1,\fF_2)}\ge\|\dg_0 f\|_{\fM_0(\fF_1\times\fF_2)}$ easily reduces to the case of finite sets $\fF_1$ and $\fF_2$.

Let us proceed to the inequality $\|f\|_{\CL(\fF_1,\fF_2)}\le\|\dg_0 f\|_{\fM_0(\fF_1\times\fF_2)}$, which means that $\|f(N_1)R-Rf(N_2)\|\le\|\dg_0 f\|_{\fM_0(\fF_1\times\fF_2)}\|N_1R-RN_2\|$
for every bounded operator $R$  and for arbitrary bounded normal operators
$N_1$ and $N_2$  such that $\s(N_1)\subset\fF_1$ and $\s(N_2)\subset\fF_1$. 
It follows from the special case treated above that this inequality certainly holds in the case when the normal operators
$N_1$ and $N_2$ have finite spectra. The case of arbitrary normal operators $N_1$ and $N_2$ with spectra in $\fF_1$ and $\fF_2$ reduces to this special case
with the help of Lemma \ref{apr}. $\bl$

\medskip

{\bf Remark.} The inequality $\|f\|_{\CL(\fF_1,\fF_2)}\le\|\dg_0 f\|_{\fM(\fF_1\times\fF_2)}$ can also be proved with the help of double operator integrals, see the remark to Theorem \ref{doikomlip}.

\medskip

In the case when $\fF_1=\fF_2$, Theorem
\ref{f1f2} reduces to the following result:

\begin{thm}
\label{olm}
Let $f$ be a function on a nonempty closed subset
$\fF$ of $\C$. Then $f\in\CL(\fF)$
if and only if $\dg_0 f\in\fM(\fF\times\fF)$. Herewith
$$
\|f\|_{\CL(\fF)}=\|\dg_0 f\|_{\fM_0(\fF\times\fF)}\le\|\dg_0 f\|_{\fM(\fF\times\fF)}\le2\|f\|_{\CL(\fF)}.
$$
\end{thm}

Note that if $\dg_0 f\in\fM(\fF\times\fF)$ for a function $f$ on
$\F$, then it is continuous and even satisfies the Lipschitz condition.
Indeed, if $\z,\t\in\fF$, then
$|(\dg_0 f)(\t,\z)|\le\|\dg_0 f\|_{\fM_0(\fF\times\fF)}$,
whence $|f(\z)-f(\t)|\le\|\dg_0 f\|_{\fM_0(\fF\times\fF)}|\z-\t|$.

The following assertion was obtained in \cite{JW}. 

\begin{thm}
\label{pro}
Let $f$ be a function on a closed subset
 $\fF$ of $\C$ such that $\dg_0 f\in\fM(\fF\times\fF)$.
Then $f$ is differentiable in the complex sense at each nonisolated point of 
$\fF$. Moreover, if $\F$ is unbounded, then there exists a final limit
$\lim\limits_{|z|\to\be}z^{-1}f(z)$.
\end{thm}

We will need an elementary lemma, which is given here without a proof.

\begin{lem}
\label{37}
Let $S$ and $T$ be arbitrary sets.
Suppose that a sequence  
$\{{\bs\f}_n\}$ of functions on $S\times T$
converges pointwise to a function $\bs\f$. Then
$\|\bs\f\|_{\fM(S\times T)}\le\varliminf\limits_{n\to\be}\|{\bs\f}_n\|_{\fM(S\times T)}$.
\end{lem}

{\bf Proof of Theorem \ref{pro}.} Let us first prove differentiability at each nonisolated point $a$ of $\fF$.
Without loss of generality we may assume that $a=0$ and $f(0)=0$. We have to show that the function $z^{-1}f(z)$ has finite limit as $z\to0$. 
Suppose that this function has at least two finite (because $f$ is Lipschitz) limit values as $z\to0$. Clearly, we may assume that these limit values are $1$ and $-1$.
Thus, there exist sequences $\{\l_n\}_{n\ge1}$ and $\{\mu_n\}_{n\ge1}$  of points of $\fF\setminus\{0\}$ that tend to zero and such that $\lim\limits_{n\to\be}\l_n^{-1}f(\l_n)=1$  and
$\lim\limits_{n\to\be}\mu_n^{-1}f(\mu_n)=-1$. Passing, if necessary, to subsequences, we can achieve the following conditions:

a) $|\l_n|>|\mu_n|>|\l_{n+1}|$ for every $n\ge1$;

b) $\lim\limits_{n\to\be}\mu_n^{-1}\l_n=0$ and $\lim\limits_{n\to\be}\l_{n+1}^{-1}\mu_n=0$.

Clearly,
$\|\{(\dg_0 f)(\l_m,\mu_n)\}\|_{\fM(\nn\times\nn)}\le\|\dg_0 f\|_{\fM(\fF\times\fF)}$.
Note that the sequence \lb$\left\{\|\{(\dg_0 f)(\l_{m+k},\mu_{n+k})\}\|
_{\fM(\nn\times\nn)}\right\}_{k\ge1}$  is nonincreasing and
$$
\lim_{k\to\be}(\dg_0 f)(\l_{m+k},\mu_{n+k})=\sgn(m-n+1/2).
$$
It follows now from Lemma \ref{37} that  
$\|\{\sgn(m-n+1/2)\}\|_{\fM(\nn\times\nn)}<+\be$
which contradicts Theorem \ref{Boch}.

The existence of a finite limit $\lim\limits_{|z|\to\be}z^{-1}f(z)$
in the case of an unbounded set $\fF$ can be proved in a similar way
with the only difference that now we should select sequences  
$\{\l_n\}_{n\ge1}$ and $\{\mu_n\}_{n\ge1}$ that tend to infinity. $\bl$

\begin{cor} 
\label{tollin}
The space $\CL(\C)$ coincides with the set of linear functions 
$az+b$, where $a,b\in\C$.
\end{cor}

\Pf Clearly, every function of the form $az+b$ with $a,b\in\C$ belongs to
$\CL(\C)$. Conversely, it follows from Theorem \ref{pro} that
$f$ is an entire function. Clearly, $f'$ is bounded
because $\CL(\C)\subset\OL(\C)\subset\Li(\C)$. It remains to use Liouville's theorem. $\bl$

\begin{thm}
\label{muld}
Let $f$ be a function on a perfect set $\fF$ in $\C$. Then $f\in\CL(\fF)$ if and only if $f$ is differentiable in the complex sense at each point of the set
$\fF$ and $\dg f\in\fM(\fF\times\fF)$. Moreover, $\|f\|_{\CL(\fF)}=\|\dg f\|_{\fM(\fF\times\fF)}$.
\end{thm}

\Pf If $f\in\CL(\fF)$,  then $\dg_0 f\in\fM(\fF\times\fF)$ by Theorem \ref{olm}
and Corollary \ref{32}. The differentiability of $f$ follows from Theorem
\ref{pro}. Conversely, if $\dg f\in\fM(\fF\times\fF)$, then
$\dg_0 f\in\fM_0(\fF\times\fF)$ and we can apply Theorem \ref{olm}.
The equality  $\|f\|_{\CL(\fF)}=\|\dg f\|_{\fM(\fF\times\fF)}$ follows from Theorem \ref{olm},
Lemma \ref{33} and the obvious equality
$\|\dg f\|_{\fM_0(\fF\times\fF)}=\|\dg_0 f\|_{\fM_0(\fF\times\fF)}$. $\bl$

\medskip

The following theorem shows that to estimate quasicommutator norms, there is no need to consider all normal operators 
$N_1$ and $N_2$, but it suffices to consider only one pair of normal operators $N_1$ and $N_2$ such that $\s(N_1)=\fF_1$ and $\s(N_2)=\fF_2$.
In particular, when keeping $\CL(\fF)$ in mind,  we may assume that $N_1=N_2$,
i.e., consider only one normal operator  $N=N_1=N_2$ such that $\s(N)=\fF$.

\begin{thm} 
\label{n1n2}
Let $N_1$ and $N_2$ be normal operators
on Hilbert spaces $\h_1$
and $\h_2$. Suppose that a continuous function $f$ 
on $\s(N_1)\cup\s(N_2)$ has the property:
\bay
\label{nerkomlip1}
\|f(N_1)R-Rf(N_2)\|\le\|N_1R-RN_2\|
\ey
for every $R\in\mB(\h_2,\h_1)$. Then
$f\in\CL(\s(N_1),\s(N_2))$  and $\|f\|_{\CL(\s(N_1),\s(N_2))}\le1$. 
\end{thm}

Let $f$ be a continuous function on a subset of the complex plane. Suppose that 
$N_1$ and $N_2$ are normal operators on Hilbert spaces $\h_1$ and $\h_2$ and the union of their spectra is contained in the domain of $f$. We say that the pair $(N_1,N_2)$ is {\it$f$-regular} if inequality \rf{nerkomlip1} holds for all
$R\in\mB(\h_2,\h_1)$.

Theorem \ref{n1n2} can be reformulated in the following way:

\medskip

{\em If a pair $(N_1,N_2)$ of normal operators is 
$f$-regular, then an arbitrary pair $(M_1,M_2)$ of normal operators 
with $\s(M_1)\subset\s(N_1)$ and $\s(M_2)\subset\s(N_2)$} is also $f$-regular.

\medskip

Let us first prove a lemma.

\begin{lem} 
\label{n1n2lem}
Let $(N_1,N_2)$ be an $f$-regular pair of bounded normal operators on
$\h_1$ and $\h_2$, and let $\K_1$ and $\K_2$ be reducing
subspaces of these operators. 
If $M_1$ is unitarily equivalent to $N_1\big|\K_1$ and $M_2$ is unitarily equivalent to $N_2\big|\K_2$, then $(M_1,M_2)$ is $f$-regular.
\end{lem}

\Pf Let $M_1\in\mB(\widetilde{\h_1})$ and 
$M_2\in\mB(\widetilde{\h_2})$.
It suffices to consider the following two partial cases:

1. $\K_1=\h_1$ and $\K_2=\h_2$.
Then $M_1=U_1^*N_1U_1$ and $M_2=U_2^*N_2U_2$
for some unitary operators $U_1$ and $U_2$.
We have
\begin{align*}
\|f(M_1)R-Rf(M_2)\|&=\|U_1^*f(N_1)U_1R-RU_2^*f(N_2)U_2\|\\[.2cm]
&=\|f(N_1)U_1RU_2^*-U_1RU_2^*f(N_2)\|\\[.2cm]
&\le\|N_1U_1RU_2^*-U_1RU_2^*N_2\|
=\|M_1R-RM_2\|
\end{align*}
for every $R$ in $\mB(\widetilde{\h_2},\widetilde{\h_1})$.

2. $M_1=N_1\big|\K_1$ and $M_2=N_2\big|\K_2$. Let
$P_1$ be the orthogonal projection from $\h_1$ onto $\K_1$
and let  $P_2$ be the orthogonal projection from $\h_2$ onto $\K_2$.
If $R\in\mB(\K_2,\K_1)$, then
\begin{align*}
\|f(M_1)R-Rf(M_2)\|&=\|P_1(f(M_1)R-Rf(M_2))P_2\|\\[.2cm]
&=\|P_1(f(N_1)R-Rf(N_2))P_2\|
=\|f(N_1)P_1RP_2-P_1RP_2f(N_2)\|\\[.2cm]
&\le\|N_1P_1RP_2-P_1RP_2N_2\|
=\|M_1R-RM_2\|.\quad\bl
\end{align*}

\medskip

{\bf Proof of Theorem \ref{n1n2}.} By Theorem  
\ref{anbnrn}, it suffices to consider the case of bounded operators $N_1$ and $N_2$. In the case when the spectra of $N_1$ and $N_2$ are finite,
Theorem \ref{n1n2} follows immediately from Lemma \ref{n1n2lem} and Theorem \ref{28}. It follows from Lemma \ref{apr} that it remains to prove that
for arbitrary finite subsets $\D_1$ and $\D_2$ of $\s(N_1)$ and $\s(N_2)$,
there are normal operators $M_1\in\mB(\K_1)$ and $M_2\in\mB(\K_2)$
such that $\s(M_1)=\D_1$, $\s(M_2)=\D_2$ and
$
\|f(M_1)R-Rf(M_2)\|\le\|M_1R-RM_2\|
$
for every $R$ in $\mB(\K_2,\K_1)$.
With each normal operator $N$ we associate the function $\a_N$ such that
$\a_N(\z)$ is the spectral multiplicity of $N$ at an isolated point $\z$ of its spectrum $\s(N)$ and
$\a_N(\z)=\be$ at each nonisolated point $\z$ of its spectrum.

We can take for operators $M_1$ and $M_2$ normal operators
on Hilbert spaces $\K_1$ and
$\K_2$ that satisfy the properties:

1)  $\s(M_1)=\D_1$ and $\s(M_2)=\D_2$;

2)  the functions $\a_{M_1}$ and $\a_{M_2}$ are the restrictions of $\a_{N_1}$ and $\a_{N_2}$.

Let $\D_1^{(\e)}$ and $\D_2^{(\e)}$ be the closed $\e$-neighborhoods of
$\D_1$ and $\D_2$. Let $N_1^{(\e)}$ be the restriction of $N_1$  to the subspace
$E_{N_1}\big(\D_1^{(\e)}\cap\s(N_1)\big)$ and let $N_2^{(\e)}$ be the restriction of $N_2$  to the subspace
$E_{N_2}\big(\D_2^{(\e)}\cap\s(N_2)\big)$, where $E_{N_1}$ and $E_{N_2}$ are  the spectral measures of $N_1$ and $N_2$. It is easy to see that there exist operators $M_1^{(\e)}$ in $\mB(\K_1)$ and $M_2^{(\e)}$ in $\mB(\K_2)$ such that $M_1^{(\e)}$ is unitarily equivalent to
$N_1^{(\e)}$,  $M_2^{(\e)}$ is unitarily equivalent to $N_2^{(\e)}$, 
$\|M_1-M_1^{(\e)}\|\le\e$ and $\|M_2-M_2^{(\e)}\|\le\e$.
Then for every $R\in\mB(\K_2,\K_1)$, we have
\begin{align*}
\|f(&M_1)R-Rf(M_2)\|
\le\|R\|\cdot\big\|f(M_1)-f\big(M_1^{(\e)}\big)\big\|+
\|R\|\cdot\big\|f(M_2)-f\big(M_2^{(\e)}\big)\big\|\\[.2cm]
+&\big\|f\big(M_1^{(\e)}\big)R-Rf\big(M_2^{(\e)}\big)\big\|
\le\|R\|\cdot\big\|f(M_1)-f\big(M_1^{(\e)}\big)\big\|
+\|R\|\cdot\big\|f(M_2)-f\big(M_2^{(\e)}\big)\big\|\\[.2cm]
+&\big\|M_1^{(\e)}R-RM_2^{(\e)}\big\|
\le\|R\|\cdot\big\|f(M_1)-f\big(M_1^{(\e)}\big)\big\|
+\|R\|\cdot\big\|f(M_2)-f\big(M_2^{(\e)}\big)\big\|\\[.2cm]
+&2\e\|R\|
+\|M_1R-RM_2\|.
\end{align*}
It remains to pass to the limit as $\e\to0$. $\bl$

\medskip

The following theorem is contained in \cite{JW}.

\begin{thm} 
Let $M$ and $N$ be operators on a Hilbert space $\h$
such that $N$ is normal. The following statements are equivalent:

{\em(a)} $M=f(N)$ for some $f$ in $\CL(\s(N))$;

{\em(b)} there exists a constant $c$ such that $\|MR-RM\|\le c\|NR-RN\|$ for every bounded operator $R$;

{\em(c)} there exists a constant $c$ such that $\|MR-RM\|_{\bS_1}\le c\|NR-RN\|_{\bS_1}$ for every bounded operator $R$;

{\em(d)} for each bounded operator $T$ there exists a bounded operator $S$
such that $SN-NS=TM-MT$;

{\em(e)} for each compact operator $T$ there exists a bounded operator $S$
such that $SN-NS=TM-MT$;

{\em(f)} for each $T$ in $\bS_1(\h)$ there exists an operator $S$
in $\bS_1(\h)$ such that $SN-NS=TM-MT$.
\end{thm}

\medskip

\section{\bf Schur multipliers and operator Lipschitzness}
\setcounter{equation}{0}
\label{dmshol}

\medskip

Note that if a closed set $\fF$ is a Fuglede set,
then $\OL(\fF)=\CL(\fF)$ by Theorem \ref{Fug}. Therefore, in this case Theorem \ref{olm}
gives a complete description of $\OL(\fF)$ in terms of Schur multipliers.

In particular, for subsets $\fF$ of a line or a circle we have a complete
description of $\OL(\fF)$ in terms of Schur multipliers.
Moreover, in the last case the seminorm of an operator Lipschitz function can be expressed in terms of the norm of the corresponding Schur multiplier. 

In the case when $\fF$ is not a Fuglede set we do not know a complete description of the operator Lipschitz functions on $\fF$ in terms of Schur multipliers.

In this case we can offer the following sufficient condition for operator Lipschitzness.

\begin{thm}
\label{teoroboperlipots} 
Let $f$ be a continuous function on a closed subset $\fF$
of $\C$. Suppose that there are Schur multipliers 
$\Phi_1,\Phi_2\in\fM(\fF\times\fF)$ such that 
$$
f(z)-f(w)=(z-w)\Phi_1(z,w)+(\ov z-\ov w)\Phi_2(z,w).
$$
Then $f\in\OL(\fF)$ and 
$
\|f\|_{\OL(\fF)}\le\|\Phi_1\|_{\fM(\fF\times\fF)}+\|\Phi_2\|_{\fM(\fF\times\fF)}.
$
\end{thm}

This theorem can be proved with the help of approximation by operators with finite 
spectra as it has been done in the proof of Theorem \ref{f1f2}. We omit this proof and instead give a proof based on double operator integrals, see Theorem \ref{olivsdoi} and the remark to it.

\medskip

{\bf Remark.} Sometimes it is more convenient to use Theorem \ref{teoroboperlipots} in terms of the real variables:
$z=x_1+{\rm i}y_1$ and $w=x_2+{\rm i}y_2$.  Suppose that there are
Schur multipliers $F_1,F_2\in\fM(\fF\times\fF)$ such that 
$$
f(z)-f(w)=(x_1-x_2)F_1(z,w)+(y_1-y_2)F_2(z,w).
$$
Then $f\in\OL(\fF)$ and 
$$
\|f\|_{\OL(\fF)}\le\frac12\|F_1+{\rm i}F_2\|_{\fM(\fF\times\fF)}+\frac12\|F_1-{\rm i}F_2\|_{\fM(\fF\times\fF)}
\le\|F_1\|_{\fM(\fF\times\fF)}+\|F_2\|_{\fM(\fF\times\fF)}.
$$

\begin{thm} 
Let $f\in\OL(\fF)$, where $\fF$ is a closed set in
$\C$. Then for every line $l$, the restriction $f\big|l\cap\fF$
is differentiable at each nonisolated point of $l\cap\fF$
and at $\be$ if the set $l\cap\fF$ is unbounded.
\end{thm}

\Pf Clearly, $f\big|l\cap\fF\in\OL(l\cap\fF)$. It remains to observe that
$\CL(l\cap\fF)=\OL(l\cap\fF)$ by Theorem \ref{su} and apply Theorem
\ref{muld} to the function $f\big|l\cap\fF$. $\bl$

\begin{cor}  Let $f\in\OL(\fF)$, where $\fF$ is a closed subset of the complex plane. Then $f$ is differentiable in an arbitrary direction at each interior point of $\fF$.
\end{cor}

{\bf Remark.}  A function $f$ in $\OL(\fF)$ does not have to be differentiable as a function of two real variables. For example, it is easy to verify that the function $f$ defined in the polar coordinates by $f(r,\theta)=re^{3{\rm i}\theta}$
belongs to $\OL(\C)$, but it is not differentiable at the origin as a function of two real variables in the cartesian coordinates.
This was observed in \cite{AP6},  see also \cite{A1}.

\medskip

\section{\bf The role of double operator integrals}
\setcounter{equation}{0}
\label{doinastene}

\medskip

In this section we demonstrate the role of double operator integrals in estimates of operator differences and (quasi)commutators. We start with estimates of operator differences under a perturbation of a self-adjoint operator by a Hilbert--Schmidt operator and discuss the Birman--Solomyak formula.

Next, we return to the results of two previous sections where we have obtained conditions for commutator Lipschitzness and operator Lipschitzness in terms of
the membership of certain functions in the space of discrete Schur multipliers.
In this section we give another proof of the sufficiency of those conditions with the help of double operator integrals. We obtain useful formulae that express operator differences and commutators in terms of double operator integrals. 

Finally, we obtain formulae for operator derivatives in terms of double operator integrals.

The following theorem was obtained by M.S. Birman and M.Z. Solomyak in \cite{BS3}.

\begin{thm}
\label{BiSoS2}
Let $f$ be a Lipschitz function on $\R$, and let $A$ and $B$ be self-adjoint operators on Hilbert space whose difference $A-B$ belongs to the Hilbert--Schmidt class $\bS_2$. Then the following formula holds:
\bay
\label{infoABS2}
f(A)-f(B)=\int_\R\!\int_\R\big(\dg_0 f\big)(x,y)\,dE_A(x)(A-B)\,dE_B(y).
\ey
\end{thm}

Note that formula \rf{infoABS2} implies straightforwardly the inequality
$$
\|f(A)-f(B)\|_{\bS_2}\le\|f\|_{\Li}\|A-B\|_{\bS_2}.
$$
In other words, Lipschitz functions are $\bS_2$-Lipschitz. It turns out that 
Lipschitz functions are also $\bS_p$-Lipschitz for $p\in(1,\be)$. This was proved recently in \cite{PS}. Recall that for $p=1$, the corresponding statement is false. This was first proved in \cite{F2}. Moreover, the class of $\bS_1$-Lipschitz functions coincides with the class of operator Lipschitz functions, see Theorem \ref{tsentrez}.

We proceed now to commutator Lipschitzness. 

\begin{thm}
\label{doikomlip}
Let $\fF_1$ and $\fF_2$ be closed subsets of $\C$. Suppose that $f$ is a continuous function on $\fF_1\cup\fF_2$ such that 
the function $\dg_0 f$ defined by {\em\rf{razdra0}} belongs to the class
$\fM(\fF_1\times\fF_2)$ of Schur multipliers. If $N_1$ and $N_2$ are normal operators such that $\s(N_j)\subset\fF_j$, $j=1,\,2$, and $R$ is a bounded linear operator, then the following formula holds:
\bay
\label{fkommudoi}
f(N_1)R-Rf(N_2)=
\int_{\!\fF_1}\!\!\int_{\!\fF_2}
\big(\dg_0 f\big)(\z_1,\z_2)\,dE_1(\z_1)(N_1R-RN_2)\,dE_2(\z_2),
\ey
where $E_j$ is the spectral measure of $N_j$.
\end{thm}

{\bf Remark.}
It follows immediately from \rf{fkommudoi} that
$$
\|f(N_1)R-Rf(N_2)\|\le\|\dg_0 f\|_{\fM(E_1,E_2)}\|N_1R-RN_2\|
\le\|\dg_0 f\|_{\fM(\fF_1\times\fF_2)}\|N_1R-RN_2\|,
$$
and, in particular, $f$ is a commutator Lipschitz function.

\medskip 

In the special case when $R$ is the identity operator, we obtain the following result:

\begin{thm}
\label{doioplip}
Let $\fF$ be a closed subset of $\C$ and let $f$ be a continuous function on $\fF$ such that $\dg_0 f\in\fM(\fF\times\fF)$. If $N_1$ and $N_2$ are normal operators with spectra in $\fF$, then the following formula holds:
\bay
\label{frazndoi}
f(N_1)-f(N_2)=
\int_\fF\int_\fF
\big(\dg_0 f\big)(\z_1,\z_2)\,dE_1(\z_1)(N_1-N_2)\,dE_2(\z_2).
\ey
\end{thm}

\medskip

{\bf Proof of Theorem \ref{doikomlip}.} Suppose first that $N_1$ and $N_2$ are bounded. We have
\begin{align*}
\int_{\fF_1}\!\!\int_{\fF_2}\big(\dg_0 f\big)(\z_1,\z_2)\,dE_1(\z_1)&(N_1R-RN_2)\,dE_2(\z_2)\\[.2cm]
&=\int_{\fF_1}\!\!\int_{\fF_2}
\big(\dg_0 f\big)(\z_1,\z_2)\,dE_1(\z_1)N_1R\,dE_2(\z_2)\\[.2cm]
&-\int_{\fF_1}\!\!\int_{\fF_2}
\big(\dg_0 f\big)(\z_1,\z_2)\,dE_1(\z_1)RN_2\,dE_2(\z_2).
\end{align*}

It follows from the definition of double operator integrals that
$$
\int_{\fF_1}\!\!\int_{\fF_2}
\big(\dg_0 f\big)(\z_1,\z_2)\,dE_1(\z_1)N_1R\,dE_2(\z_2)
=
\int_{\fF_1}\!\!\int_{\fF_2}
\z_1\big(\dg_0 f\big)(\z_1,\z_2)\,dE_1(\z_1)R\,dE_2(\z_2)
$$
and
$$
\int_{\fF_1}\!\!\int_{\fF_2}
\big(\dg_0 f\big)(\z_1,\z_2)\,dE_1(\z_1)RN_2\,dE_2(\z_2)
=
\int_{\fF_1}\!\!\int_{\fF_2}
\z_2\big(\dg_0 f\big)(\z_1,\z_2)\,dE_1(\z_1)R\,dE_2(\z_2).
$$
Since $(\z_1-\z_2)(\dg_0 f\big)(\z_1,\z_2)=f(\z_1)-f(\z_2)$, $\z_1\in\fF_1$, $\z_2\in\fF_2$, we obtain
\begin{align*}
\int_{\fF_1}\!\!\int_{\fF_2}
\big(\dg_0 f\big)(\z_1,\z_2)&\,dE_1(\z_1)(N_1R-RN_2)\,dE_2(\z_2)\\
=&\int_{\fF_1}\!\!\int_{\fF_2}
f(\z_1)\,dE_1(\z_1)R\,dE_2(\z_2)-
\int_{\fF_1}\!\!\int_{\fF_2}
f(\z_2)\,dE_1(\z_1)R\,dE_2(\z_2).
\end{align*}
Again, it is easy to see from the definition of double operator integrals that
$$
\int_{\fF_1}\!\!\int_{\fF_2}f(\z_1)\,dE_1(\z_1)R\,dE_2(\z_2)=
\left(\,\,\int_{\fF_1}f(\z_1)\,dE_1(\z_1)\right)R=f(N_1)R
$$
and
$$
\int_{\fF_1}\!\!\int_{\fF_2}f(\z_2)\,dE_1(\z_1)R\,dE_2(\z_2)=
R\int\limits_{\fF_1}f(\z_1)\,dE_1(\z_1)=Rf(N_2),
$$
which implies \rf{fkommudoi}.

Suppose now that $N_1$ and $N_2$ are not necessarily bounded normal operators.
The special case of Theorem \ref{doikomlip} that has been proved above and  Theorem \ref{anbnrn} imply the commutator estimate, and so the operator
$f(N_1)R-Rf(N_2)$ is bounded.

Put
$$
P_k\df E_1\big(\{\z\in\C:|\z|\le k\}\big)\quad\mbox{and}\quad
Q_k\df E_2\big(\{\z\in\C:|\z|\le k\}\big),\quad k>0.
$$
Then
$N_{1,k}\df P_kN_1$ and $N_{2,k}\df Q_kN_2$ are bounded normal operators. Let $E_{j,k}$ be the spectral measure of
$N_{j,k}$, $j=1,\,2$. We have
\begin{align*}
P_k&\left(\int_{\fF_1}\!\!\int_{\fF_2}\big(\dg_0 f\big)(\z_1,\z_2)\,
dE_1(\z_1)(N_1R-RN_2)\,dE_2(\z_2)\right)Q_k\\[.2cm]
&=P_k\left(\int_{\fF_1}\!\!\int_{\fF_2}\big(\dg_0 f\big)(\z_1,\z_2)\,
dE_{1,k}(\z_1)(P_kf(N_1)R-Rf(N_2)Q_k)\,dE_{2,k}(\z_2)\right)Q_k.
\end{align*}

Applying \rf{fkommudoi} to the bounded normal operators $N_{1,k}$ and $N_{2,k}$, we obtain
\begin{align*}
P_k\big(f(N_{1,k})R&-Rf(N_{2,k})\big)Q_k=\\[.2cm]
=&P_k\left(\int_{\fF_1}\!\!\int_{\fF_2}\big(\dg_0 f\big)(\z_1,\z_2)\,
dE_{1,k}(\z_1)(P_kN_1R-RN_2Q_k)\,dE_{2,k}(\z_2)\right)Q_k.
\end{align*}
Since
$$
P_k\big(f(N_{1,k})R-Rf(N_{2,k})\big)Q_k=P_k\big(f(N_{1})R-Rf(N_{2})\big)Q_k,
$$
we have
\begin{align*}
P_k\big(f(N_{1})R&-Rf(N_{2})\big)Q_k=\\[.2cm]
=&P_k\left(\int_{\fF_1}\!\!\int_{\fF_2}\big(\dg_0 f\big)(\z_1,\z_2)\,
dE_1(\z_1)(N_1R-RN_2)\,dE_2(\z_2)\right)Q_k.
\end{align*}
It remains to pass to the limit in the strong operator topology. $\bl$

\medskip

It is easy to verify that in all formulae of this section one can replace
the function $\dg_0f(\z_1,\z_2)$ in double operator integrals
with an arbitrary bounded measurable function $F(\z_1,\z_2)$ that
coincides with $\dg_0f(\z_1,\z_2)$ for all $\z_1$ and $\z_2$ such that $\z_1\ne\z_2$.

In particular, in the case when $\fF_1=\fF_2$ and the set $\fF_1$ is perfect,
by Theorem \ref{muld} one can replace in \rf{fkommudoi} replace $\dg_0f$ with the divided difference $\dg f$.

\begin{thm}
\label{infsovmn}
Let $\fF$  be a closed subset of $\C$ and 
$f\in\CL(\fF)$. If $N_1$ and $N_2$ are normal operators such that $\s(N_j)\subset\fF$, $j=1,\,2$, and $R$ is a bounded linear operator, then the following formula holds:
\bay
\label{fkommudoisov}
f(N_1)R-Rf(N_2)=
\int_{\fF}\!\int_{\fF}
\big(\dg f\big)(\z_1,\z_2)\,dE_1(\z_1)(N_1R-RN_2)\,dE_2(\z_2),
\ey
where $E_j$ is the spectral measure of $N_j$.
\end{thm}

Let us proceed now to the inerpretation of the results of \S\:\ref{dmshol} in terms of double operator integrals. The following assertion holds:

\begin{thm}
\label{olivsdoi}
Let $f$ be a continuous function on a closed subset $\fF$
of $\C$ and suppose that there exist 
Schur multipliers $\Phi_1,\Phi_2\in\fM(\fF\times\fF)$ such that 
$$
f(\z_1)-f(\z_2)=(\z_1-\z_2)\Phi_1(\z_1,\z_2)+(\ov\z_1-\ov\z_2)\Phi_2(\z_1,\z_2),
\quad\z_1,~\z_2\in\fF.
$$
Let $N_1$ and $N_2$ be normal operators whose spectra are contained in $\fF$.
Then 
\begin{align}
\label{intforadoi}
f(N_1)-f(N_2)&=\int_\fF\int_\fF
\Phi_1(\z_1,\z_2)\,
dE_1(\z_1)(N_1-N_2)\,dE_2(\z_2)\nonumber\\[.2cm]
&+\int_\fF\int_\fF
\Phi_2(\z_1,\z_2)(\z_1,\z_2)\,
dE_1(\z_1)(N^*_1-N^*_2)\,dE_2(\z_2).
\end{align}
\end{thm}

{\bf Remark.}
It follows easily from \rf{intforadoi} that 
$$
\|f(N_1)-f(N_2)\|\le\big(\|\Phi_1\|_{\fM(\fF\times\fF)}
+\|\Phi_2\|_{\fM(\fF\times\fF)}\big)\|N_1-N_2\|
$$
and, in particular, $f$ is an operator Lipschitz function.

\medskip

\Pf As in the proof of Theorem \ref{doikomlip}, we assume first that the operators $N_1$ and $N_2$ are bounded. Then
\begin{align*}
\int_\fF\int_\fF
\Phi_1(\z_1,\z_2)\,
&dE_1(\z_1)(N_1-N_2)\,dE_2(\z_2)\\[.2cm]
&=\int_\fF\int_\fF
\Phi_1(\z_1,\z_2)
dE_1(\z_1)N_1dE_2(\z_2)-
\int_\fF\int_\fF
\Phi_1(\z_1,\z_2)
dE_1(\z_1)N_2dE_2(\z_2)\\[.2cm]
&=\int_\fF\int_\fF
\z_1\Phi_1(\z_1,\z_2)
dE_1(\z_1)dE_2(\z_2)-
\int_\fF\int_\fF
\z_2\Phi_1(\z_1,\z_2)
dE_1(\z_1)dE_2(\z_2)\\[.2cm]
&=\int_\fF\int_\fF
(\z_1-\z_2)\Phi_1(\z_1,\z_2)\,dE_1(\z_1)\,dE_2(\z_2).
\end{align*}

Similarly,
$$
\int_\fF\int_\fF
\Phi_2(\z_1,\z_2)(\z_1,\z_2)
dE_1(\z_1)(N^*_1-N^*_2)dE_2(\z_2)=
\int_\fF\int_\fF(\ov\z_1-\ov\z_2)\Phi_2(\z_1,\z_2)
\,dE_1(\z_1)\,dE_2(\z_2).
$$
Therefore, the right-hand side of \rf{intforadoi} is equal to 
\begin{align*}
\int_\fF\int_\fF(f(\z_1)&-f(\z_2))\,dE_1(\z_1)\,dE_2(\z_2)\\[.2cm]
&=\int_\fF f(\z_1)\,dE_1(\z_1)
-\int_\fF f(\z_2)\,dE_2(\z_2)=f(N_1)-f(N_2).
\end{align*}
The passage from bounded operators to unbounded ones can be done in the same way as in the proof of Theorem \ref{doikomlip}. $\bl$

\medskip

Let us proceed now to applications of double operator integrals in problems of operator differentiability.

\begin{thm} 
\label{sildif}
Let $f$ be an operator Lipschitz function on $\R$, and let
$A$ and $K$ be self-adjoint operators and suppose that $K$ is bounded.
Then
$$
\lim_{t\to0}\frac1t(f(A+tK)-f(A))=\int_\R\int_\R(\dg f)(x,y)\,dE_A(x)K\,dE_A(y),
$$
where the limit is taken in the strong operator topology.
\end{thm}

We need several auxiliary results.
Let $\widehat\R\df\R\cup\{\be\}$ denote the {\it one point compactification
of the real line} $\R$. Recall that every function $f\in\OL(\R)$ is everywhere differentiable on $\widehat\R$, see Theorem \ref{pro}.

\begin{lem}
\label{abchag} 
Let $f\in\OL(\R)$. Then there exist sequences $\{\f_n\}_{n\ge0}$ and 
$\{\psi_n\}_{n\ge0}$ of continuous functions
on $\widehat\R$ such that 

{\rm a)} $\sum\limits_{n\ge0}|\f_n|^2\le\|f\|_{\OL(\R)}$ everywhere on $\widehat\R$,

{\rm b)} $\sum\limits_{n\ge0}|\psi_n|^2\le\|f\|_{\OL(\R)}$ everywhere on $\widehat\R$,

{\rm c)} $(\dg f)(x,y)=\sum\limits_{n\ge0}\f_n(x)\psi_n(y)$ for all $x,\,y\in\R$.
\end{lem}

\Pf By Theorem \ref{muld}, we have $\dg f\in\fM(\R\times\R)$ and $\|\dg f\|_{\fM(\R\times\R)}=\|f\|_{\OL(\R)}$.
We extend the function $\dg f$ to the set $\widehat\R\times\widehat\R$ by putting
$(\dg f)(x,y)=f'(\be)=\lim_{t\to\be}t^{-1}f(t)$ in the case when $|x|+|y|=\be$. Clearly, this extended function $\dg f$ on $\widehat\R\times\widehat\R$
is continuous in each variable. Hence,
\begin{align*}
\|\dg f\|_{\fM(\widehat\R\times\widehat\R)}&=
\sup\{\|\dg f\|_{\fM(\L_1\times\L_2)}:~\L_1,\,\L_2\subset\widehat\R, \,\,
\L_1\,\,\text {and} \,\,\L_2\,\,\,\text{are finite}\}\\[.2cm]
&=\sup\{\|\dg f\|_{\fM(\L_1\times\L_2)}:~\L_1,\,\L_2\subset\R, \,\,
\L_1\,\,\text {and} \,\,\L_2\,\,\,\text{are finite}\}=\|\dg f\|_{\fM(\R\times\R)}.
\end{align*}
It remains to apply Theorem \ref{npkp} to the function $\dg f:\widehat\R\times\widehat\R\to\C$. $\bl$

\begin{lem} 
\label{nepratk-}
Let $A$ and $K$ be self-adjoint operators such that $K$ is bounded.
Then for every function $f$ in $C(\widehat\R)$ the function $H$, $H(t)\df f(A+tK)$, acts continuously from $\R$ to the space $\mB(\h)$ with the normed topology.
\end{lem}

Note that in \cite{AP2} a considerably stronger result was obtained.

\medskip

\Pf We may assume that $f(\be)=0$. Then we can construct a sequence 
$\{f_n\}_{n\ge0}$
of functions of class $C^\be$ with compact support such that $f_n\to f$ uniformly.
Each function $H_n$, $H_n(t)\df f_n(A+tK)$, is continuous because $f_n\in\OL(\R)$
for $n\ge0$. It remains to observe that $H_n\to H$ uniformly. $\bl$

\begin{lem} 
\label{Xnun-}
Let $\{X_n\}_{n\ge0}$ be a sequence in $\mB(\h)$ and let
$\{u_n\}_{n\ge0}$ be a sequence in $\h$. Suppose that 
$\sum_{n\ge0}X_nX_n^*\le a^2I$ and $\sum_{n\ge0}\|u_n\|^2\le b^2$
for some nonnegative numbers $a$ and $b$. Then the series $\sum_{n\ge0}X_nu_n$
converges weakly and 
$$
\Big\|\sum_{n\ge0}X_nu_n\Big\|\le ab.
$$
\end{lem}

\Pf Let $v\in\h$ and $\|v\|=1$. Then
$$
\sum_{n\ge0}|(X_nu_n,v)|=\sum_{n\ge0}|(u_n,X_n^*v)|\le
\Big(\sum_{n\ge0}\|u_n\|^2\Big)^{1/2}\Big(\sum_{n\ge0}\|X_n^*v\|^2\Big)^{1/2}\le ab
$$
which implies the result. $\bl$

\medskip

{\bf Proof of Theorem \ref{sildif}.} 
By formulae \rf{fkommudoisov} and \rf{Haagrazl}, Theorem \ref{sildif} can be reformulated in the following way:
$$
\lim_{t\to0}\sum_{n\ge0}\f_n(A+tK)K\psi_n(A)=\sum_{n\ge0}\f_n(A)K\psi_n(A)
$$
in the strong operator topology, where $\f_n$ and $\psi_n$ are functions from the conclusion of Lemma \ref{abchag}.
In other words, we have to prove that for every $u\in\h$, we have
$$
\lim_{t\to0}\sum_{n\ge0}(\f_n(A+tK)-\f_n(A))K\psi_n(A)u=\0,
$$
where the series is understood in the weak topology of $\h$ while the limit is taken in the norm of $\h$.
Assume that $\|u\|=1$ and $\|f\|_{\OL(\R)}=1$.  Then $\sum_{n\ge0}|\f_n|^2\le1$ and $\sum_{n\ge0}|\psi_n|^2\le1$ everywhere on $\R$.

Put $u_n\df K\psi_n(A)u$. We have
$$
\sum_{n\ge0}\|u_n\|^2\le\|K\|^2\sum_{n\ge0}\|\psi_n(A)u\|^2=
\|K\|^2\sum_{n\ge0}(|\psi_n|^2(A)u,u)\le\|K\|^2<+\be.
$$
Let $\e>0$. Let us choose a natural number $N$ such that 
$\sum_{n>N}\|u_n\|^2<\e^2$. Then it follows from Lemma \ref{Xnun-}
that
$$
\Big\|\sum_{n>N}(\f_n(A+tK)-\f_n(A))u_n\Big\|\le2\e
$$
for all $t\in\R$. Note that by Lemma \ref{nepratk-},
$$
\Big\|\sum_{n=0}^N(\f_n(A+tK)-\f_n(A))u_n\Big\|\le\|K\|\sum_{n=0}^N\|\f_n(A+tK)-\f_n(A)\|<\e
$$
for all $t$ sufficiently close to zero.
Thus, 
$$
\Big\|\sum_{n\ge0}(\f_n(A+tK)-\f_n(A))u_n\Big\|<3\e
$$
for all $t$ sufficiently close to zero. $\bl$

\medskip

By analogy with Theorem \ref{sildif}, we can prove the following fact:

\begin{thm} 
\label{sildifloc-}
Let $A$ and $K$ be bounded self-adjoint operators.
Then
$$
\lim_{t\to0}\frac1t(f(A+tK)-f(A))=\int_\R\int_\R(\dg f)(x,y)\,dE_A(x)K\,dE_A(y)
$$
for every $f$ in $\OL_{\rm loc}(\R)$, where the limit is taken in the strong operator topology.
\end{thm}

Theorem \ref{sildif} implies the following result:

\begin{thm}
Let $f$ be an operator differentiable function on $\R$, and let $A$ and $K$ be self-adjoint operators such $K$ is bounded. Then for the derivative of the function $t\mapsto f(A+tK)-f(A)$ in the operator norm, the following formula holds:
\bay
\label{propoopno}
\frac{d}{dt}\big(f(A+tK)-f(A)\big)\Big|_{t=0}=
\int_\R\int_\R\frac{f(x)-f(y)}{x-y}\,dE_A(x)K\,dE_A(y).
\ey
\end{thm} 

In particular, formula \rf{propoopno} holds for an arbitrary function $f$ of Besov class $B_{\be,1}^1(\R)$, see Theorem \ref{Besdiffer}.

Similar results hold for functions on the unit circle.

\medskip

\section{\bf Trace class Lipschitzness and trace class commutator Lipschitzness}
\setcounter{equation}{0}
\label{yadiyadcom}

\medskip

The purpose of this section is to prove that for an arbitrary closed set $\fF$ in the plane, the classes $\CL(\fF)$ and $\CL_{\bS_1}(\fF)$ coincide. In particular, if $\fF\subset\R$, then  $\OL(\fF)=\OL_{\bS_1}(\fF)$ (see \S\:\ref{oplip} where the classes $\CL_{\bS_1}(\fF)$ and $\OL_{\bS_1}(\fF)$ are defined). 

Note that the definition of $\CL_{\bS_1}(\fF)$ can be extended naturally to the definition of the class $\CL_{\bS_1}(\fF_1,\fF_2)$,
where $\fF_1$ and $\fF_2$ are nonempty closed subsets of $\C$.

\begin{lem}
\label{f1f2s1}
Let $f$ be a continuous function on the union
$\fF_1\cup\fF_2$ of closed subsets $\fF_1$ and $\fF_2$
of $\C$. Then
$$
\|f\|_{\CL_{\bS_1}(\fF_1,\fF_2)}\ge\|\dg_0 f\|_{\fM_{0,\bS_1}(\fF_1\times\fF_2)}
\ge\frac12\|\dg_0 f\|_{\fM(\fF_1\times\fF_2)}.
$$
\end{lem}

\Pf The second inequality follows from Corollary \ref{32s1}. Let us prove the first one.
It suffices to consider the case of finite sets $\fF_1$ and $\fF_2$.
Let $N_1$ and $N_2$ be normal operators with simple spectra such that $\s(N_1)=\fF_1$ and $\s(N_2)=\fF_2$.
Then there exist orthonormal bases $\{u_\l\}_{\l\in\s(N_1)}$
and $\{v_\mu\}_{\mu\in\s(N_2)}$ in $\h_1$ and $\h_2$
such that $N_1u_\l=\l u_\l$ for $\l\in\s(N_1)$
and $N_2v_\mu=\mu v_\mu$ for $\mu$ in $\s(N_2)$.  With each operator
$X:\h_2\to\h_1$ we associate the matrix
$\{(Xv_\mu,u_\l)\}_{(\l,\mu)\in\s(N_1)\times\s(N_2)}$. We have
$$
\left((N_1R-RN_2)v_\mu,u_\l\right)=
\big(Rv_\mu,N_1^*u_\l)-(RN_2v_\mu,u_\l)=
(\l-\mu)(Rv_\mu,u_\l).
$$
Similarly,
$$
\left((f(N_1)R-Rf(N_2))v_\mu,u_\l\right)=(f(\l)-f(\mu))(Rv_\mu,u_\l).
$$
It is easy to see that
$$
\{(f(\l)-f(\mu))(Rv_\mu,u_\l)\}
=\{(\dg_0 f)(\l,\mu)\}\star
\{(\l-\mu)(Rv_\mu,u_\l)\}.
$$
Note that the matrix $\{a_{\l\mu}\}_{(\l,\mu)\in\s(N_1)\times\s(N_2)}$ can be represented in the form
$$
\{a_{\l\mu}\}_{(\l,\mu)\in\s(N_1)\times\s(N_2)}=
\{(\l-\mu)(Rv_\mu,u_\l)\}_{(\l,\mu)\in\s(N_1)\times\s(N_2)}
$$
for an operator $R$ from $\h_2$ to $\h_1$ if and only if
$a_{\l\mu}=0$ with $\l=\mu$.
Now the inequality $\|f\|_{\CL(\fF_1,\fF_2)}\ge\|\dg_0 f\|_{\fM_{0,\bS_1}(\fF_1\times\fF_2)}$
is obvious. $\bl$

\begin{cor}
\label{napryamoi}
Let $f$ be a real continuous function on a closed subset $\fF$  
of $\R$. Then
\bay
\label{tsepravnerav}
\|f\|_{\OL_{\bS_1}(\fF)}=\|f\|_{\CL_{\bS_1}(\fF)}\ge\|\dg_0 f\|_{\fM_{0,\bS_1}(\fF\times\fF)}
\ge\frac12\|\dg_0 f\|_{\fM(\fF\times\fF)}.
\ey
If $f$ is not necessarily real, then 
\bay
\label{vtortsepner}
\|f\|_{\OL_{\bS_1}(\fF)}\ge\frac12\|f\|_{\CL_{\bS_1}(\fF)}\ge\frac12\|\dg_0 f\|_{\fM_{0,\bS_1}(\fF\times\fF)}
\ge\frac14\|\dg_0 f\|_{\fM(\fF\times\fF)}.
\ey
\end{cor}

\Pf The equality in \rf{tsepravnerav} follows from Corollary \ref{sledoveshch}.
All the inequalities in \rf{tsepravnerav} have been already proved above. Obviously, \rf{vtortsepner} follows from \rf{tsepravnerav} $\bl$

\medskip

Let us proceed now to the main results of this section.

\begin{thm}
\label{sovkliyadkl}
Let $\fF$ be a closed set in $\C$. Then $\CL(\fF)\!=\!\CL_{\bS_1}(\fF)$ and
$$
\frac12\|f\|_{\CL_{\bS_1}(\fF)}\le\|f\|_{\CL(\fF)}\le2\|f\|_{\CL_{\bS_1}(\fF)},
\quad f\in\CL_(\fF).
$$
\end{thm}

\Pf Let $f\in\CL(\fF)$, and let $N_1$ and $N_2$ be normal operators with spectra in $\fF$ such that $N_1R-RN_2\in\bS_1$ and $R$ is a bounded operator. Then by the remark to Theorem \ref{doikomlip},
\begin{align*}
\|f(N_1)R-Rf(N_2)\|_{\bS_1}
&\le\big\|\dg_0f\|_{\fM(E_1,E_2)}\|N_1R-RN_2\|_{\bS_1}\\
&\le\big\|\dg_0f\|_{\fM(\fF\times\fF)}\|N_1R-RN_2\|_{\bS_1}
\le2\|f\|_{\CL(\fF)}\|N_1R-RN_2\|_{\bS_1}
\end{align*} by Theorem \ref{olm}.
This implies the inequality $\|f\|_{\CL_{\bS_1}(\fF)}\le2\|f\|_{\CL(\fF)}$.
On the other hand, we obtain from Lemma \ref{f1f2s1} and Theorem \ref{olm} 
$$
\|f\|_{\CL(\fF)}\le\big\|\dg_0f\big\|_{\fM(\fF\times\fF)}\le2\|f\|_{\CL_{\bS_1}(\fF)}.\quad\bl
$$

\medskip

If $\fF$ is a perfect set, we can improve the result.

\begin{thm}
\label{sovershenstvo}
Let $\fF$ be a perfect set in $\C$. Then 
$\|f\|_{\CL_{\bS_1}(\fF)}=\|f\|_{\CL(\fF)}$ for every 
$f$ in $\CL(\fF)=\CL_{\bS_1}(\fF)$.
\end{thm}

\Pf By \rf{fkommudoisov}, for $f\in\CL(\fF)$ we have
\begin{align*}
\|f(N_1)R-Rf(N_2)\|_{\bS_1}
&\le\big\|\dg f\|_{\fM(E_1,E_2)}\|N_1R-RN_2\|_{\bS_1}\\
&\le\big\|\dg f\|_{\fM(\fF\times\fF)}\|N_1R-RN_2\|_{\bS_1}
=\|f\|_{\CL(\fF)}\|N_1R-RN_2\|_{\bS_1}.
\end{align*}
The last equality is guaranteed by Theorem \ref{muld}. Thus, we have proved that
\lb$\|f\|_{\CL_{\bS_1}(\fF)}\le\|f\|_{\CL(\fF)}$. 
Applying now Lemmata \ref{f1f2s1} and \ref{33s1} as well as Theorem \ref{muld}, we obtain
$$
\|f\|_{\CL_{\bS_1}(\fF)}\ge\|\dg_0 f\|_{\fM_{0,\bS_1}(\fF\times\fF)}
=\|\dg f\|_{\fM(\fF\times\fF)}=\|f\|_{\CL(\fF)}.\quad\bl
$$

\medskip

It is time to proceed to the central result of this section.

\begin{thm}
\label{tsentrez}
Let $f$ be a continuous function on $\R$. The following statements are equivalent:

{\em(a)} $f$ is operator Lipschitz;

{\em(b)} $f$ is trace class Lipschitz;

{\em(c)} $f(A)-f(B)\in\bS_1$, whenever $A$ and $B$ are self-adjoint operators with $A-B$ in $\bS_1$.
\end{thm} 

Note that in (c) {\it it is necessary to consider not only bounded operators} $A$ and $B$.

\medskip

\Pf The equivalence of (а) and (b) is established in Corollary \ref{napryamoi}. The implication (b)$\Rightarrow$(c) is trivial. Let us show that (c)$\Rightarrow$(b). Suppose that $f\not\in\CL_{\bS_1}(\R)$. Then we can find sequences $A_n$ and $B_n$ of self-adjoint operators such that $A_n-B_n\in\bS_1$ and
$\|A_n-B_n\|_{\bS_1}^{-1}\|f(A_n)-f(B_n)\|_{\bS_1}\to\be$ as $n\to\be$.
Without loss of generality we may assume that $\lim_{n\to\be}\|A_n-B_n\|_{\bS_1}=0$.
Indeed, consider the increment $A_n\mapsto A_n+K_n$, where $K_n\df B_n-A_n$. 
Consider now the following increments: 
$A_n+(j/M_n)K_n\mapsto A_n+((j+1)/M_n)K_n$,
$0\le j\le M_n-1$, where $\{M_n\}$ is a sequence of natural numbers such that 
$\lim_{n\to\be}\|A_n-B_n\|_{\bS_1}/M_n=0$. We select now $j$ that maximizes the number
$$
\|f(A_n+((j+1)/M_n)K_n)-f(A_n+(j/M_n)K_n)\|_{\bS_1}
$$
and replace the pair $(A_n,B_n)$ with the pair 
$(A_n+(j/M_n)K_n, A_n+((j+1)/M_n)K_n)$. Then
$$
\lim_{n\to\be}\|A_n-B_n\|_{\bS_1}^{-1}\|f(A_n)-f(B_n)\|_{\bS_1}=\be\quad\mbox{and}
\quad\lim_{n\to\be}\|A_n-B_n\|_{\bS_1}=0.
$$
It suffices now, if necessary, to select a subsequence of the sequence
$(A_n,B_n)$ or repeat certain terms of the sequence to achieve the condition
$$
\sum_n\|B_n-A_n\|_{\bS_1}<\be\quad\mbox{but}\quad
\sum_n\|f(B_n)-f(A_n)\|_{\bS_1}=\be.
$$
Let $A$ be the orthogonal sum of the $A_n$ and let
$B$ be the orthogonal sum of the $B_n$. Then $B-A\in\bS_1$ but
$f(B)-f(A)\not\in\bS_1$. $\bl$

\medskip

{\bf Remark.} A similar result holds for functions on the unit circle and for unitary operators.

\medskip

\section{\bf Operator Lipschitz functions on the plane. Sufficient conditions}
\setcounter{equation}{0}
\label{OLnaplBes}

\medskip

In this section we obtain a sufficient condition for operator Lipschitzness in terms of the Besov class $B_{\be,1}^1(\R^2)$. It is similar to Theorem \ref{Besdost} for functions on the real line. The results of this section were obtained in \cite{APPS}.

Recall (see \rf{frazndoi}) that in the case of functions on the line, the operator Lipschitzness of $f$ can be obtained from the formula
$$
f(A)-f(B)=\int_\R\int_\R\frac{f(s)-f(t)}{s-t}\,dE_A(s)(A-B)\,dE_B(t).
$$
Here $A$ and $B$ are self-adjoint operators. This is the way the operator Lipschitzness of functions of class $B_{\be,1}^1(\R)$ was established in \cite{Pe1} and \cite{Pe3}.

It would be natural to try the same approach for functions on the plane. 
However, (see Corollary \ref{tollin}) if the divided difference is a Schur multiplier for arbitrary Borel spectral measures on $\C$, the function must be linear.

In \cite{APPS} another method was used: for normal operators $N_1$ and $N_2$, the difference $f(N_1)-f(N_2)$ is represented as the sum of double operator integrals, the integrands being the divided differences in each variable.

We introduce the following notation. Let $N_1$ and $N_2$ be normal operators on Hilbert space. Put
$A_j\df\re N_j$, $B_j\df\im N_j$, $j=1,\,2$, and let $E_j$ be the spectral measure of $N_j$.
In other words, $N_j=A_j+{\rm i}B_j$, $j=1,\,2$, where $A_j$ and $B_j$ are commuting self-adjoint operators. 

If $f$ is a function on $\R^2$ that has partial derivatives in each variable, we consider the divided differences in each variable
$$
\big(\dg_xf\big)(z_1,z_2)\df\frac{f(x_1,y_2)-f(x_2,y_2)}{x_1-x_2},
\quad z_1,\,z_2\in\C,
$$
and
$$
\big(\dg_yf\big)(z_1,z_2)\df\frac{f(x_1,y_1)-f(x_1,y_2)}{y_1-y_2},
\quad z_1,\,z_2\in\C,
$$
where
$x_j\df\re z_j$, $y_j\df\im z_j$, $j=1,\,2$.
On the sets $\{(z_1,z_2):~x_1=x_2\}$ and
$\{(z_1,z_2):~y_1=y_2\}$
the divided differences are understood as the corresponding partial derivatives of  $f$.

The following result gives us a key estimate.

\begin{thm}
\label{SMdd}
Let $f$ be a bounded continuous function on $\R^2$ whose Fourier transform 
$\F f$ has compact support. Then $\dg_xf$ and 
$\dg_yf$ are Schur multipliers of class $\fM(\C\times\C)$.
Moreover, if
$\supp\F f\subset\{\z\in\C:~|\z|\le\s\}$, $\s>0$,
then
\bay
\label{nerxy}
\|\dg_xf\|_{\fM(\C\times\C)}\le\const\s\|f\|_{L^\be}
\quad\mbox{and}\quad\|\dg_yf\|_{\fM(\C\times\C)}\le \const\s\|f\|_{L^\be}.
\ey
\end{thm}

It follows from the definition of the Besov class $B_{\be,1}^1(\R^2)$ and from Theorem \ref{SMdd} that for every $f\in B_{\be,1}^1(\R^2)$, the divided differences $\dg_xf$ and $\dg_yf$ are Schur multipliers and
\bay
\label{nerdlyaBes}
\|\dg_xf\|_{\fM(\C\times\C)}\le\const\|f\|_{B_{\be,1}^1}
\quad\mbox{and}\quad\|\dg_yf\|_{\fM(\C\times\C)}\le \const\|f\|_{B_{\be,1}^1}.
\ey

Inequalities \rf{nerdlyaBes} together with Theorem \ref{olivsdoi} imply the following result obtained in \cite{APPS} and being the central result of this section.

\begin{thm}
\label{doin}
Let $f$ be a function in $B_{\be,1}^1(\R^2)$. Suppose that $N_1$ and $N_2$ are normal operators such that $N_1-N_2$ is bounded. Then
\begin{align*}
f(N_1)-f(N_2)&=\iint\limits_{\C^2}\big(\dg_yf\big)(z_1,z_2)\,
dE_1(z_1)(B_1-B_2)\,dE_2(z_2)\nonumber\\[.2cm]
&+\iint\limits_{\C^2}\big(\dg_xf\big)(z_1,z_2)\,
dE_1(z_1)(A_1-A_2)\,dE_2(z_2)
\end{align*}
and the inequality
$
\|f(N_1)-f(N_2)\|\le\const\|f\|_{B_{\be,1}^1}\|N_1-N_2\|
$ holds, i.e., $f$ is an operator Lipschitz function on $\C$.
\end{thm}

To prove Theorem \ref{SMdd}, we use a formula for a representation of the divided difference as an element of the Haagerup tensor product. Recall that $\mathscr E_\s$ denotes the set of entire functions (of one complex variable)
of exponential type at most $\s$.

\begin{lem}
\label{predosnnaformKotSh}
Let $\f\in\mathscr E_\s\cap L^\be(\R)$. Then
\begin{align}
\label{haa1}
\frac{\f(x)-\f(y)}{x-y}&=\sum_{n\in\Z}\s\cdot\frac{\f(x)-\f\big(\pi n\s^{-1}\big)}{\s x-\pi n}
\cdot\frac{\sin(\s y-\pi n)}{\s y-\pi n}.
\end{align}
Moreover,
\bay
\label{2vy}
\sum_{n\in\Z}\frac{\big|\f(x)-\f\big(\pi n\s^{-1}\big)\big|^2}{(\s x-\pi n)^2}
\le 3\|\f\|_{L^\be(\R)}^2,
\quad x\in\R,
\ey
and
\bay
\label{hiz}
\sum_{n\in\Z}\frac{\sin^2(\s y-\pi n)}{(\s y-\pi n)^2}=1,\quad y\in\R.
\ey
\end{lem}

We refer the reader to \cite{APPS} where in Section 5 a proof of Lemma \ref{predosnnaformKotSh} is given that is based on the Kotel'nikov--Shannon formula, which, in turn, is based on the fact that the family of functions$\{(z-\pi n)^{-1}\sin(z-\pi n)\}_{n\in\Z}$
forms an orthonormal basis in $\mathscr E_1\cap L^2(\R)$,
see \cite{L}, Lecture 20, Section 3.

\medskip

{\bf Proof of Theorem \ref{SMdd}.} Clearly, $f$ is the restriction to $\R^2$ of an entire function of two complex variables.
Moreover, $f(\cdot,a),\,f(a,\cdot)\in \mathscr E_\s\cap L^\be(\R)$ for every $a\in\R$.
Without loss of generality we may assume that $\s=1$. By Lemma \ref{predosnnaformKotSh},
\bey
\big(\dg_xf\big)(z_1,z_2)=\frac{f(x_1,y_2)-f(x_2,y_2)}{x_1-x_2}=
\sum_{n\in\Z}(-1)^n\frac{f(\pi n,y_2)-f(x_2,y_2)}{\pi n-x_2}\cdot\frac{\sin(x_1-\pi n)}{x_1-\pi n}
\eey
and
\bey
\big(\dg_yf\big)(z_1,z_2)=\frac{f(x_1,y_1)-f(x_1,y_2)}{y_1-y_2}=\sum_{n\in\Z}(-1)^n\frac{f(x_1,y_1)-f(x_1,\pi n)}{y_1-\pi n}\cdot\frac{\sin(y_2-\pi n)}{y_2-\pi n}.
\eey

Note that the expressions $\dfrac{\sin(x_1-n\pi)}{x_1-\pi n}$ and $\dfrac{f(x_1,y_1)-f(x_1,\pi n)}{y_1-\pi n}$ 
depend on $z_1=(x_1,y_1)$ but do not depend on $z_2=(x_2,y_2)$, while
the expressions $\dfrac{f(\pi n,y_2)-f(x_2,y_2)}{\pi n-x_2}$ and $\dfrac{\sin(y_2-\pi n)}{y_2-\pi n}$
depend $z_2=(x_2,y_2)$ but do not depend on $z_1=(x_1,y_1)$. Moreover, by Lemma \ref{predosnnaformKotSh}
\bey
\sum_{n\in\Z}\frac{|f(x_1,y_1)-f(x_1,\pi n)|^2}{(y_1-\pi n)^2}\le3\|f(x_1,\cdot)\|_{L^\be(\R)}^2
\le3\|f\|_{L^\be(\C)}^2,\\
\sum_{n\in\Z}\frac{|f(\pi n,y_2)-f(x_2,y_2)|^2}{(\pi n-x_2)^2}\le3\|f(\cdot,y_2)\|_{L^\be(\R)}^2
\le3\|f\|_{L^\be(\C)}^2\\
\eey
and
$$
\sum_{n\in\Z}\frac{\sin^2(x_1-\pi n)}{(x_1-\pi n)^2}=
\sum_{n\in\Z}\frac{\sin^2(y_2-\pi n)}{(y_2-\pi n)^2}=1,
$$
which proves \rf{nerxy}. $\bl$

\medskip

Note that inequalities  \rf{nerxy} play the role of operator Bernstein's inequalities (see \S\:\ref{Bern}) as in the case of functions of self-adjoint operators, one can prove the following results:

\begin{thm}
\label{Hanormal}
Let $0<\a<1$ and let $f$ a function of the H\"older class $\L_\a(\R^2)$.
Then 
$$
\|f(N_1)-f(N_2)\|\le c(1-\a)^{-1}\|f\|_{\L_\a}\|N_1-N_2\|^\a
$$
for some constant $c>0$ and
for arbitrary normal operators $N_1$ and $N_2$ with bounded difference $N_1-N_2$.
\end{thm}

One can generalize Theorem \ref{Hanormal}  to the case of arbitrary moduli of continuity and obtain an analog of Theorem \ref{prmone}.

\begin{thm}
\label{Spnormal}
Let $0<\a<1$, $p>1$ and let $f$ be a function of the H\"older class 
$\L_\a(\R^2)$. Then there exists a positive number $c$ such that
$$
\|f(N_1)-f(N_2)\|_{S_{p/\a}}\le c\,\|f\|_{\L_\a}\|N_1-N_2\|_{\bS_p}^\a
$$
for arbitrary normal operators $N_1$ and $N_2$ whose difference belongs to the Schatten--von Neumann class $\bS_p$.
\end{thm}

We refer the reader to \cite{APPS} where he can find proofs of these results as well as other related results.

\medskip

\section{\bf A sufficient condition for commutator Lipschitzness in terms of Cauchy integrals}
\setcounter{equation}{0}
\label{dostol}

\medskip

In this section we give a sufficient condition for commutator Lipschitzness that was obtained in \cite{A2}.

Let $\fF$ be a nonempty closed subset of $\C$ such that 
$\fF\ne\C$. We denote by $\M(\C\setminus \fF)$
the space of complex Radon measures $\mu$ on $\C\setminus\fF$ such that
\bay
\label{uslname}
\|\mu\|_{\M(\C\setminus \fF)}\df
\sup_{z\in\fF}\int_{\C\setminus\fF}\frac{d|\mu|(\z)}{|\z-z|^2}<+\be.
\ey
For $\mu\in\M(\C\setminus \fF)$, the Cauchy integral
$$
\widehat\mu(z)=\int_{\C\setminus\fF}\frac{d\mu(\z)}{\z-z}\,\,,
$$
is not defined in general even for $z\in\fF$ because the function 
$\z\mapsto(\z-z)^{-1}$ does not have to be integrable with respect to the measure $|\mu|$.
With each fixed point $z_0\in\fF$ we associate the modified Cauchy integral
$$
\widehat\mu_{z_0}(z)\df\int_{\C\setminus\fF}\Big(\frac1{\z-z}-\frac1{\z-z_0}\Big)\,d\mu(\z).
$$
It follows from the Cauchy--Bunyakovsky inequality that
$\widehat\mu_{z_0}(z)$ is well defined for $z\in\fF$ and
\mbox{$|\widehat\mu_{z_0}(z)|\le\|\mu\|_{\M(\C\setminus \fF)}|z-z_0|$}.
Moreover, $\widehat\mu_{z_0}(z_1)=-\widehat\mu_{z_1}(z_0)$ and
$$
|\widehat\mu_{z_0}(z_1)-\widehat\mu_{z_0}(z_2)|=|\widehat\mu_{z_1}(z_1)-\widehat\mu_{z_1}(z_2)|
=|\widehat\mu_{z_1}(z_2)|
\le\|\mu\|_{\M(\C\setminus \fF))}|z_1-z_2|
$$
for all $z_1,z_2\in\fF$.

Note that $z\mapsto(\z-z)^{-1}$ is a continuous map from $\fF$ to the Hilbert space $L^2(|\mu|)$ endowed with the weak topology.
This allows one to verify easily that the function $\widehat\mu_{z_0}(z)$
is differentiable as a function of complex variable at every nonisolated point of  $\fF$. In particular, $\widehat\mu_{z_0}(z)$ is analytic in the interior of $\fF$.

We denote by $\widehat\M(\fF)$
the set of functions
$f$ on $\fF$ that can be represented in the form $f=c+\widehat\mu_{z_0}$, where 
$c$ is a constant. Put
$$
\|f\|_{\widehat\M(\fF)}\df
\inf\{\|\mu\|_{\M(\C\setminus \fF)}:~\mu\in\M(\C\setminus \fF),\,\,\,
f-\widehat\mu_{z_0}=\const\,\,\, \text{on}\,\,\,\fF\}.
$$
It is easy to see that the definition of $\widehat\M(\fF)$ and the seminorm 
$\|f\|_{\widehat\M(\fF)}$
do not depend on the choice of $z_0\in\fF$.

\begin{thm}
\label{coint7}
Let $\fF$ be a proper closed subset of $\C$.
Then $\widehat\M(\fF)\subset\CL(\fF)$ and $\|f\|_{\CL(\fF)}\le\|f\|_{\widehat\M(\fF)}$ for every $f$ in $\widehat\M(\fF)$.
\end{thm}

\Pf Let $\mu\in\M(\C\setminus \fF)$ and $f=\widehat\mu_{z_0}$. 
Consider the divided difference
\begin{align*}
\frac{f(z)-f(w)}{z-w}&=
\frac1{z-w}\int_{\C\setminus\fF}\left(\frac1{\z-z}-\frac1{\z-w}\right)\,d\mu(\z)
\\[.2cm]
&=\int_{\C\setminus\fF}\frac{d\mu(\z)}{(\z-z)(\z-w)}.
\end{align*}
Inequality \rf{uslname} means that this divided difference satisfies condition (d) of Theorem \ref{tomSc}. Thus, it is a Schur multiplier for arbitrary Borel spectral measures on $\fF$ and its multiplier norm is at most $\|\mu\|_{{\mathscr M}(\C\setminus\fF)}$. It remains to refer to Theorem \ref{olm}. $\bl$

\medskip

Let $\widehat\M_\be(\fF)$ denote the space of functions of the form
$f+az$, where $f\in\widehat\M(\fF)$, $a\in\C$. It is easy to see that the linear function 
$az$ belongs to $\widehat\M(\fF)$ if $\fF$ is compact. Therefore,
$\widehat\M_\be(\fF)=\widehat\M(\fF)$ for compact $\fF$.
In the case of an unbounded set $\fF$, it is easy to verify that $f'(\be)=0$ for every $f\in\widehat\M(\fF)$. Thus, $\widehat\M_\be(\fF)\ne\widehat\M(\fF)$
for noncompact sets $\fF$. It follows from Theorem \ref{coint7} that
$\widehat\M_\be(\fF)\subset\CL(\fF)$.

The authors do not know whether the equality $\widehat\M_\be(\fF)=\CL(\fF)$
holds even for such simple sets $\fF$ as the circle or the line.

\

\section{\bf Commutator Lipschitz functions on the disk and on the half-plane}
\setcounter{equation}{0}
\label{full}

\

We consider here the spaces of commutator Lipschitz functions in the unit disk  
$\dd$ and in the upper halfplane  $\C_+\df\{\z\in\C:~\re\z>0\}$. In particular, we present the results by E. Kissin and V.S. Shulman \cite{KS3} and their analogs for the upper half-plane.

Let $\CA$ denote the disk-algebra, i.e., the space of functions $f$ analytic in the open disk $\dd$ and continuous in its closure. It was proved in \cite{KS3} that $\CL(\dd)=\{f\in \CA:~f\in\OL(\T)\}$. The following result shows that this equality is isometric.

\begin{thm}
\label{KST}
Let $f\in\CL(\dd)$. Then $f\in \CA$ and 
 $\|f\|_{\CL(\dd)}=\|f\|_{\CL(\T)}=\|f\|_{\OL(\T)}$. If $f\in \CA$, then $f\in\CL(\dd)$ if and only if $f\in\OL(\T)$.
\end{thm}

\Pf The equality $\|f\|_{\CL(\T)}=\|f\|_{\OL(\T)}$ follows from Theorem \ref{su}.
The inequality $\|f\|_{\CL(\T)}\le\|f\|_{\CL(\dd)}$ is obvious. It remains to prove that
$\|f\|_{\CL(\dd)}\le\|f\|_{\CL(\T)}$. We may assume that $\|f\|_{\CL(\T)}=1$.
Then $\|\dg f\|_{\fM(\T\times\T)}=1$ by Theorem \ref{muld}. Let us apply Theorem \ref{topolx}.
We obtain families $\{u_\z\}_{\z\in\T}$ and $\{v_\t\}_{\t\in\T}$
in a Hilbert space $\h$ that depend on parameters continuously in the weak topology and such that $\|u_\z\|\le1$, 
$\|v_\t\|\le1$  and
$(\dg f)(\z,\t)=(u_\z,v_\t)$ for all $\z,\,\t\in\T$. Consider the harmonic extensions of the functions
$\z\mapsto u_\z$ and $\t\mapsto v_\t$ to the unit disk by putting
$$
u_z\df\int_{\T}\frac{1-|z|^2}{|z-\z|^2}u_\z\,d\m(\z)\quad\text{and}\quad
v_w\df\int_{\T}\frac{1-|w|^2}{|w-\t|^2}v_\t\,d\m(\t)
$$
for $z,w\in\dd$. The integrals are understood as integrals of $\h$-valued functions continuous in the weak topology.
Applying the Poisson integral in the variable $\z$ to both parts of the equality
$(\dg f)(\z,\t)=(u_\z,v_\t)$, we obtain $(\dg f)(z,\t)=(u_z,v_\t)$ for 
$z\in\clos\dd$ and $\t\in\T$. Applying now the Poisson integral to the last equality, we obtain $(\dg f)(z,w)=(u_z,v_w)$ for all $z\in\clos\dd$
and $w\in\clos\dd$. Now it is clear that
$$
\|f\|_{\CL(\dd)}=\|\dg f\|_{\fM(\clos\dd\times\clos\dd)}\le\sup_{z\in\clos\dd}\|u_z\|\sup_{w\in\clos\dd}\|v_w\|
=\sup_{\z\in\T}\|u_\z\|\sup_{\t\in\T}\|v_\t\|=1.\quad \bl
$$

\medskip

We give now an analog of Theorem \ref{npkp} for functions in the unit disk.

\begin{thm} 
\label{KSTc}
Let $f\in\CL(\dd)$. Then there are sequences $\{\f_n\}_{n\ge1}$ and 
$\{\psi_n\}_{n\ge1}$ in the disk-algebra $\CA$ such that
$$
\left(\sup_{z\in\dd}\sum_{n=1}^\be|\f_n(z)|^2\right)\left(\sup_{w\in\dd}\sum_{n=1}^\be|\psi_n(w)|^2\right)
=\|f\|_{\CL(\dd)}^2\!\!\quad\text{and}\!\!\quad
(\dg f)(z,w)=\sum_{n=1}^\be\f_n(z)\psi_n(w),
$$
wherein the series converge uniformly while $z$ and $w$ range over compact subsets of the unit disk.
\end{thm}

\Pf We may assume that $\|f\|_{\CL(\dd)}=1$. To prove the first equality, it suffices to prove the inequality $\le$, for the inequality $\ge$ follows from Theorems \ref{muld} and 
\ref{41}. Let $\h$, $u_z$ and $v_w$ denote the same as in the proof of Theorem
\ref{KST}. Consider an orthonormal basis 
$\{e_n\}_{n=1}^\be$ in $\h$.
Put $\f_n(z)\df(u_z,e_n)$ and $\psi_n(w)\df(e_n,v_w)$ for $n\ge1$. Let us prove that $\f_n\in \CA$ and $\psi_n\in \CA$. Denote by $X$ the set of vectors 
$e\in\h$ such that $(u_z,e)\in \CA$. Clearly, $X$ is a closed subspace
of $\h$.  Note that $v_\t\in X$ for every $\t\in\T$ because
$(\dg f)(\cdot,\t)\in \CA$ for every $\t\in\T$. Thus, $X=\h$, for
the linear span of $\{v_\t\}_{\t\in\T}$ is dense in $\h$.
Therefore, $(u_z,e)\in \CA$ for every $e\in\h$. Similarly, one can prove that
$(e,v_w)\in \CA$ for every $e\in\h$. It remains to prove uniform convergence on compacta. Note that
$$
\left|\sum_{n=N}^\be\f_n(z)\psi_n(w)\right|\le\left(\sum_{n=N}^\be|\f_n(z)|^2\right)^{\frac12}
\left(\sum_{n=N}^\be|\psi_n(w)|^2\right)^{\frac12}.
$$
Thus, it suffices to establish uniform convergence on the compact subsets of $\dd$  for the series $\sum_{n=1}^\be|\f_n(z)|^2$ and $\sum_{n=1}^\be|\psi_n(z)|^2$.
This is a consequence of the following elementary lemma.

\begin{lem} 
\label{kner}
Let $\{h_k\}_{k=1}^\be$ be a sequence of analytic functions in $\dd$.
Suppose that the function $\sum_{k=1}^\be|h_k(z)|$ is bounded in $\dd$. Then the series
$\sum_{k=1}^\be|h_k(z)|$ converges uniformly on compact subsets of the open unit disk.
\end{lem}

We denote by $(\OL)_+(\T)$ the space of functions $f$ in $\OL(\T)$ that admit
an analytic extension to the unit disk $\dd$ and are continuous in its closure.
It follows from Theorem \ref{pro} that every function $f\in\CL(\dd)$ is analytic in $\dd$. Thus, Theorem \ref{KST} implies the following result of \cite{KS3}.

\begin{thm} 
\label{olclt}
The operator of restriction $f\mapsto f\big|\T$ is a linear isometry of
$\CL(\dd)$ onto $(\OL)_+(\T)$.
\end{thm}

Similar results also hold for the space $\CL(\C_+)$. 

\begin{thm} 
\label{KSR}
Let $f$ be a continuous function in the closed upper half-plane $\clos\C_+$.
Suppose that $f$ is analytic in the open half-plane $\C_+$. 
Then $\|f\|_{\CL(\C_+)}=\|f\|_{\CL(\R)}=\|f\|_{\OL(\R)}$. In particular, $f\in\CL(\C_+)$ if and only if $f\in\OL(\R)$.
\end{thm}

We denote by $\CA(\C_+)$ the set of functions analytic in $\C_+$ continuous in $\clos\C_+$ and having finite limit at infinity.

\begin{thm} 
\label{KSRc}
Let $f\in\CL(\C_+)$. Then there are sequences $\{\f_n\}_{n=1}^\be$ and 
$\{\psi_n\}_{n=1}^\be$ in $\CA(\C_+)$ such that$$
\left(\sup_{z\in\C_+}\sum_{n=1}^\be|\f_n(z)|^2\right)\left(\sup_{w\in\C_+}\sum_{n=1}^\be|\psi_n(w)|^2\right)
=\|f\|_{\CL(\C_+)}^2
$$
and
$$
(\dg f)(z,w)=\sum_{n=1}^\be\f_n(z)\psi_n(w).
$$
Herewith the series converge uniformly while $z$ and $w$ range over compact subsets of the open upper half-plane.
\end{thm}

We skip here the proofs of Theorems \ref{KSR} and \ref{KSRc}. They are similar to the proofs of the corresponding results for functions in the unit disk.

We also state the following analog for the line of Theorem \ref{olclt}.

\begin{thm} 
\label{olclr}
The operator of restriction $f\mapsto f\big|\R$ is a linear isometry of
$\CL(\C_+)$ onto $(\OL)_+(\R)$.
\end{thm}

Note that in \cite{A4} the following result was obtained that contains as a matter of fact both Theorem \ref{KST} and Theorem \ref{KSR}.

\begin{thm}
Let $\fF_0$ and $\fF$ be nonempty perfect subsets of $\C$
such that $\fF_0\subset\fF$ and $\Omega\df\fF\setminus\fF_0$ is open.
Suppose that a function $f_0\in\CL(\fF_0)$ admits a continuous extension
$f$ to $\fF$ such that $f$ is analytic in $\Omega$ and
$|f(z)z^{-2}|\to0$ as $z\to\be$
in each unbounded component of $\O$.\footnote{The last condition holds automatically if $\Omega$ is bounded.} 
Then $f\in\CL(\fF)$ and $\|f\|_{\CL(\fF)}=\|f_0\|_{\CL(\fF_0)}$.
\end{thm}

The authors do not know an answer to the following question.
Let $f$ be a continuous function in the closed unit disk that is harmonic 
inside the disk. Suppose that $f\in\OL(\T)$. Must it be true that 
$f\in\OL(\dd)$?
A similar question can be posed for the half-plane as well as for other domains. 

Recall that if $T$ is a contraction on a Hilbert space $\h$, then by the Sz\"okefalvi-Nagy theorem (see \cite{SNF}, Chapter I, \S\:5), there exists its unitary dilation {\it unitary dilation}, i.e., a unitary operator $U$ on a Hilbert space $\K$, $\h\subset\K$, such that
$T^n=P_\h U^n\big|\h$, $n\ge0$. A dilation can always be chosen minimal.
This allows us to define the following linear and multiplicative calculus:
$\f\mapsto\f(T)\df P_\h\f(U)\big|\h$, $\f\in \CA$. The
{\it semispectral measure $\mE_T$ of the contraction} $T$ is defined by $\mE_T(\D)\df P_\h E_U(\D)\big|h$, where $E_U$ is the spectral measure of $U$ and $\D$ is a Borel subset of $\T$. It is easy to see that
$$
\f(T)=\int_\T\f(\z)\,d\mE(\z),\quad\f\in \CA.
$$

\begin{thm}
\label{teorokomszha}
Let $f\in\CL(\dd)$,  let $T_1$ and $T_2$ be contractions on a Hilbert space 
$\h$, and let $R\in\mB(\h)$. Then
\bay
\label{doifszha}
f(T_1)R-Rf(T_2)=
\int_\T\int_\T\big(\dg f)(\z,\t)\,d\mE_1(\z)(T_1R-RT_2)\,d\mE_2(\t),
\ey
where $\mE_1$ and $\mE_2$ are the semi-spectral measures of $T_1$ and $T_2$,
and the following inequality holds: 
\bay
\label{KSlip}
\|f(T_1)R-Rf(T_2)\|\le\|f\|_{\CL(\dd)}\|T_1R-RT_2\|.
\ey
\end{thm}

\Pf Let $\{\f_n\}_{n\ge1}$ and $\{\psi_n\}_{n\ge1}$ be sequences of functions in the disk-algebra that satisfy the conclusion of Theorem \ref{KSTc}. By
\rf{doipolcpm}, we have

\begin{align*}
\int_\T\int_\T\big(\dg f)(\z,\t)&
\,d\mE_1(\z)(T_1R-RT_2)\,d\mE_2(\t)
=\sum_{n\ge1}\f_n(T_1)(T_1R-RT_2)\psi_n(T_2)\\[.2cm]
&=
\sum_{n\ge1}T_1\f_n(T_1)R\psi_n(T_2)-\sum_{n\ge1}\f_n(T_1)R\psi_n(T_2)T_2\\[.2cm]
&=\int_\T\int_\T\z\big(\dg f)(\z,\t)
\,d\mE_1(\z)R\,d\mE_2(\t)
-\int_\T\int_\T\t\big(\dg f)(\z,\t)
\,d\mE_1(\z)R\,d\mE_2(\t)\\[.2cm]
&=\int_\T\int_\T(f(\z)-f(\t))
\,d\mE_1(\z)R\,d\mE_2(\t)=f(T_1)R-Rf(T_2)
\end{align*}
which proves formula \rf{doifszha}, which, in turn, immediately implies inequality 
\rf{KSlip}. $\bl$

\medskip

Inequality \rf{KSlip} was proved by Kissin and Shulman in \cite{KS3} by a different method. The proof given here is similar to the proof of Theorem 4.1 in \cite{Pe9}, see also \cite{Pe*}.
In the case $f\in B_{\be,1}^1(\T)\cap \CA$ and $R=I$, Theorem \ref{teorokomszha} was proved in \cite{Pe*}. 

A similar result can be also proved for dissipative operators, see \cite{AP*} for perturbations of functions of dissipative operators.

\medskip

\section{\bf Operator Lipschitz functions and linear-fractional transformations}
\setcounter{equation}{0}
\label{dlp}

\medskip

Let $\dlC$ denote the M\"obius group of linear-fractional transformations of the extended complex plane $\widehat \C\df\C\cup\{\be\}$. In other words,
$$
\dlC=\Big\{\f:~\f(z)=\frac{az+b}{cz+d}, \quad a,\,b,\,c,\,d\in \C,\quad ad-bc\ne0\Big\}.
$$

The set of linear-fractional transformations of the complex plane
is denoted by $\lC$, i.e.,
$$
\lC=\{\f\in\dlC:~\f(\be)=\be\}=\{\f:~\f(z)=az+b, a,b\in\C, a\ne0\}.
$$

Let $\widehat\R$ denote the one point compactification of $\R$, $\widehat\R\df\R\cup\{\be\}$.
Put 
$$
\dlR\df\{\f\in\dlC:~\f(\widehat\R)=\widehat\R\}\quad\mbox{and} \quad
\lR\df\{\f\in\lC:~\f(\R)=\R\}.
$$

With each linear-fractional transformation $\f$ and each function $f$
on a closed set  $\fF,\:\fF\subset\C$, we associate the function
${\mathcal Q}_\f f$ on the set
$\fF_\f\df\C\cap\f^{-1}(\fF\cup\{\be\})$
defined by
$$
({\mathcal Q}_\f f)(z)\df\left\{\begin{array}{ll}\dfrac{f(\f(z))}{\f^{\,\prime}(z)},&\text {if}\,\,\,
z\in\C,\,\,\f(z)\in\fF\,\,\,\text{and}\,\,\, \f(z)\ne\be,\\[.3cm]
\quad 0,&\text {if}\,\,\,z\in\C\,\,\,\text{and}\,\,\, \f(z)=\be.
\end{array}\right.
$$

It is easy to see that if $\f\in\lC$, then $\fF_\f=\f^{-1}(\fF)$, ${\mathcal Q}_\f f=(\f'(0))^{-1}(f\circ\f)$,
${\mathcal Q}_\f(\OL(\fF))=\OL(\fF_\f)$, ${\mathcal Q}_\f(\CL(\fF))=\CL(\fF_\f)$, $\|{\mathcal Q}_\f f\|_{\OL(\fF_\f)}=\|f\|_{\OL(\fF)}$
for all $f$ in $\OL(\fF)$ and $\|{\mathcal Q}_\f f\|_{\CL(\fF_\f)}=\|f\|_{\CL(\fF)}$
for all $f$ in $\CL(\fF)$. Therefore, we will be mostly interested in the case when $\f\notin\lC$. Note that if $\fF=\C$, then $\fF_\f=\C$  for every $\f$ in $\dlC$.
However, if $\fF=\R$, then $\fF_\f=\R$ for every $\f$ in $\dlR$.

Let $a\in\fF$, where $\fF$ is a closed subset of $\C$.
Put 
\bay
\label{3101}
\!\!\!\OL_a(\fF)\df\{f\in\OL(\fF):f(a)=0\}\!\!\quad\mbox{and}\quad\!\!
\CL_a(\fF)\df\{f\in\CL(\fF):f(a)=0\}.
\ey
Clearly, $\OL_a(\fF)$ and $\CL_a(\fF)$ are Banach spaces.

\begin{thm}
\label{zamfone}
Let $\fF$ be a closed subset of $\C$, $a\in\fF$ and let
$\f$ be an automorphism in $\dlC$ such that $a\df\f(\be)$.
Then ${\mathcal Q}_\f(\OL_a(\fF))\subset\OL(\fF_\f)$ and
$$
\|{\mathcal Q}_\f f\|_{\OL(\fF_\f)}\le3\|f\|_{\OL(\fF)}
$$
for every $f$ in $\OL_a(\fF)$.
\end{thm}

\Pf Consider first the case $\f(z)=\phi(z)\df z^{-1}$.
Then $a=0$ and we have to prove that ${\mathcal Q}_\phi(\OL_0(\fF))\subset\OL(\fF_\f)$
and $\|{\mathcal Q}_\phi f\|_{\OL(\fF_\f)}\le3\|f\|_{\OL(\fF)}$.
Let $f\in\OL_0(\fF)$. We may assume that $\|f\|_{\OL(\fF)}=1$. 
Then everything reduces to the inequality
$$
\|({\mathcal Q}_\f f)(N)R-R({\mathcal Q}_\f f)(N)\|\le3\max(\|NR-RN\|,\|N^*R-RN^*\|)
$$
for arbitrary bounded operators $N$ and $R$ such that $N$ is normal and $\s(N)\subset\fF_\phi$.
We define the function $h$ by $h(z)=zf(z^{-1})$ for $z\ne0$ and 
$h(0)=0$.
It is easy to see that $\sup|h|\le\|f\|_{\Li(\fF)}\le\|f\|_{\OL(\R)}=1$,  for $f(0)=0$. Note that $({\mathcal Q}_\phi f)(N)=-Nh(N)$.
Thus, we have to prove that
$$
\|Nh(N)R-RNh(N)\|\le3\max\{\|NR-RN\|,\|N^*R-RN^*\|\}.
$$
Let us use the following elementary identity:
\begin{align}
\label{42}
Nh(N)R-RNh(N)&=h(N)(NR-RN)\nonumber\\[.2cm]
&+h(N)RN-NRh(N)
+(NR-RN)h(N).
\end{align}

Note that 
$$
\|h(N)(NR-RN)\|\le\|NR-RN\|\le\max\{\|NR-RN\|,\|N^*R-RN^*\|\}.
$$
In a similar way we can estimate the norm of $(NR-RN)h(N)$.

It remains to prove that
$$
\|h(N)RN-NRh(N)\|\le\max\{\|NR-RN\|,\|N^*R-RN^*\|\}.
$$

If $N$ is invertible, then
\begin{multline*}
\|h(N)RN-NRh(N)\|=\|f(N^{-1})NRN-NRNf(N^{-1}\|)\\
\le
\max\{\|RN-NR\|,\|(N^*)^{-1}NRN-NRN(N^*)^{-1}\|\}\\
=\max\{\|NR-RN\|,\|(N^*)^{-1}N(RN^*-N^*R)N(N^*)^{-1}\|\}
\\
=\max\{\|NR-RN\|,\|N^*R-RN^*\|\}.
\end{multline*}

If $0$ is a limit point of $\fF_\phi$ (i.e., the set $\fF$ is unbounded), the proof can be concluded, for in this case each normal operator $N$ with spectrum in $\fF_\phi$
can be approximated arbitrarily accurately by a normal operator $M$ such that 
$MN=NM$ and $\s(M)\subset\fF_\f\setminus\{0\}$. This follows from Lemma \ref{apr}.

Suppose now that $0$ is an isolate point of $\fF_\phi$.
Consider a noninvertible normal operator $N$ with spectrum $\fF_\phi$.
Then $N$ can be represented as $N=\0\oplus N_0$,  where $N_0$ is an invertible normal operator.
Note that ${\mathcal Q}_\phi(N)=\0\oplus N_0^2f(N_0^{-1})$. Let $P$ be the orthogonal projection onto the subspace, on which $N_0$ is defined. It is easy to see that
\begin{multline*}
\|h(N)RN-NRh(N)\|=\|P(h(N)RN-NRh(N))P\|\\=\|h(N)PRPN-NPRPh(N)\|
=\|h(N_0)(PRP)N_0-N_0(PRP)h(N_0)\|\\
\le\max(\|N_0(PRP)-(PRP)N_0\|,\|N_0^*(PRP)-(PRP)N_0^*\|\\
\le\max(\|NR-RN\|,\|N^*R-RN^*\|).
\end{multline*}

Let us proceed to the general case. Put $b=\f^{-1}(\be)$.
Clearly, $\f(z)=a+c\phi(z-b)$, where $c\in\C\setminus\{0\}$. Thus,
everything reduces to the case $a=b=0$, i.e.,  $\f=c\phi$ because translations 
preserve the operator Lipschitz norm. Finally, the case $\f=c\phi$ reduces easily
to the case $\f=\phi$ already treated. $\bl$

\medskip

{\bf Example.} Let $\f(z)=z^{-1}$, $\fF=\C$, $f=\ov z$.
Then $f\in\OL_{0}(\C)$ and $\|f\|_{\OL(\C)}=1$. Moreover, $({\mathcal Q}_\f f)(z)=-\ov z^{-1}z^2$ and
$$
3=\|{\mathcal Q}_\f f\|_{\Li(\T)}\le\|{\mathcal Q}_\f f\|_{\Li(\C)}\le\|{\mathcal Q}_\f f\|_{\OL(\C)}\le3\|f\|_{\OL(\C)}=3.
$$
This example shows that $\|{\mathcal Q}_\f f\|_{\OL(\C)}=\|{\mathcal Q}_\f f\|_{\Li(\C)}=3$ and the constant
$3$ in Theorem \ref{zamfone} is best possible.

\medskip

Theorem \ref{zamfone} implies easily the following result:

\begin{thm}
\label{zamf}
Let $\f\in\dlC$, $a=\f(\be)$ and $b=\f^{-1}(\be)$. Suppose that $\fF$ is a closed set in $\C$, that contains $a$.
Then ${\mathcal Q}_\f(\OL_a(\fF))=\OL_b(\fF_\f)$ and
$$
\frac13\|f\|_{\OL(\fF)}\le\|{\mathcal Q}_\f f\|_{\OL(\fF_\f)}\le3\|f\|_{\OL(\fF)}
\quad\mbox{for every}\quad f~\mbox{~in~}~\OL_a(\fF).
$$
\end{thm}

\Pf Note that $({\mathcal Q}_\f(\OL_a(\fF)))(b)=0$. Thus, it follows from Theorem \ref{zamfone} that ${\mathcal Q}_\f(\OL_a(\fF))\subset\OL_b(\fF_\f)$ and 
$\|{\mathcal Q}_\f f\|_{\OL(\fF_\f)}\le3\|f\|_{\OL(\fF)}$.
To prove that ${\mathcal Q}_\f(\OL_a(\fF))\supset\OL_b(\fF_\f)$ and obtain the desired lower estimate for $\|{\mathcal Q}_\f f\|_{\OL(\fF_\f)}$, it suffices to apply Theorem \ref{zamfone}
to the closed set $\fF_\f$ and the linear-fractional transformation $\f^{-1}$. 
$\bl$

\medskip

We present one more related result.

\begin{thm}
\label{zamfcom}
Let $\f\in\dlC\setminus\lC$ and let $a=\f(\be)$. Suppose that $\fF$ is a closed
subset of $\C$ such that $a\not\in\fF$. If $z_0$ is one of the closest points of  $\fF$ to $a$, then ${\mathcal Q}_\f(\OL_{z_0}(\fF))\subset\OL(\fF_\f)$ and
$$
\|{\mathcal Q}_{\f}f\|_{\OL(\fF_\f)}
\le 5\|f\|_{\OL(\fF)}
$$
for every $f\in\OL_{z_0}(\fF)$.
\end{thm}

\Pf 
As in the proof of Theorem \ref{zamfone},  it suffices to consider the case
$\f(z)=\phi(z)\df z^{-1}$.  Let $f\in\OL_{z_0}(\fF)$ and $\|f\|_{\OL_{z_0}(\fF)}=1$. We have to prove that
$$
\|({\mathcal Q}_\f f)(N)R-R({\mathcal Q}_\f f)(N)\|\le5\max(\|NR-RN\|,\|N^*R-RN^*\|)
$$
for arbitrary normal operators $N_1$ and $N_2$ such that $\s(N_1),\;\s(N_2)\subset\fF_\phi$.
Let $h$ denote the same as in the proof of Theorem \ref{zamfone}. However, we cannot say now that $\sup|h|\le1$. We have
\begin{multline*}
\sup_{z\in\fF_\phi}|h(z)|\le\sup\{|zf(z^{-1})|:z\in\phi^{-1}(\fF)\}=\sup\{|z|^{-1}|f(z)-f(z_0)|:z\in\fF\}\\
\le\sup\{|z|^{-1}|z-z_0|:z\in\fF\}\le\sup\{1+|z|^{-1}|z_0|:z\in\fF\}=2.
\end{multline*}
Repeating now the reasoning of the proof of Theorem \ref{zamfone}, we obtain
\begin{align*}
\|({\mathcal Q}_\f f)(N)R-R({\mathcal Q}_\f f)(N)\|
&\le(1+2\sup|h(z)|)\max(\|NR-RN\|,\|N^*R-RN^*\|)\\[.2cm]
&\le5\max(\|NR-RN\|,\|N^*R-RN^*\|). \quad\bl
\end{align*}

\medskip

{\bf Example.} Let $\f(z)=z^{-1}$, $\fF=\T$, $z_0=1$, $f=1-\ov z$.
Then $f\in\OL_{z_0}(\T)$ and $\|f\|_{\OL(\T)}=1$. It is easy to verify that $({\mathcal Q}_\f f)(z)=z^3-z^2$ and $\|{\mathcal Q}_\f f\|_{\Li(\T)}\ge5$. Then
$$
5\le\|{\mathcal Q}_\f f\|_{\Li(\T)}\le\|{\mathcal Q}_\f f\|_{\OL(\T)}\le\|z^3\|_{\OL(\T)}+\|z^2\|_{\OL(\T)}=5.
$$
This example shows that the constant 5 in Theorem \ref{zamfcom} is best possible.

\medskip

{\bf Remark 1.} We can introduce the following generalization of ${\mathcal Q}_\f$, by putting
$$
({\mathcal Q}_{n,\f} f)(z)\df\left\{\begin{array}{ll}\dfrac{|\f'(z)|^nf(\f(z))}{(\f^{\,\prime}(z))^{n+1}},&\text {if}\,\,\,
z\in\C,\,\,\f(z)\in\fF\,\,\,\text{and}\,\,\, \f(z)\ne\be,\\[.3cm]
\quad 0,&\text {if}\,\,\,z\in\C\,\,\,\text{and}\,\,\, \f(z)=\be,
\end{array}\right.
$$
where $n\in\Z$. Then analogs of Theorems \ref{zamfone}, \ref{zamf} and \ref{zamfcom} for the operators 
${\mathcal Q}_{n,\f}$ hold with constants depensing on $n$. Analogs of Theorems \ref{zamfone} and \ref{zamf} can be found in \cite{A1}. An analog of Theorem  \ref{zamfcom} can be obtained in the same way.

\medskip

{\bf Remark 2.} The proofs of Theorems \ref{zamfone}, \ref{zamf} and \ref{zamfcom} also work for the spaces of commutator Lipschitz functions. The case of Theorems \ref{zamfone} and \ref{zamf} 
is treated in \cite{A1}. Clearly, in the case of $\CL(\fF)$ we can speak about
generalizations to the operators ${\mathcal Q}_{n,\f}$ (see Remark 1) only for ``sparse'' sets $\fF$.
For example, if $\fF$ has interior points, then such generalizations are impossible because the functions in $\CL(\fF)$ are analytic in the interior of $\fF$.

\medskip

Later we will be mostly interested in the case when $\fF=\R$
and $\fF=\T$. In these cases we have isometric equality $\CL(\fF)=\OL(\fF)$.

Note that Theorem \ref{zamfone} implies the following result:

\begin{thm}
\label{zamrt}
Let $\f\in\dlC$. Suppose that $\f(\widehat\R)=\T$. Then
$$
\|{\mathcal Q}_\f f\|_{\OL(\R)}\le3\|f\|_{\OL(\T)}
$$
for every $f\in\OL_a(\T)$, where $a=\f(\be)$.
\end{thm}

\Pf Let us apply Theorem \ref{zamf} to $\fF=\T$. Then $\fF_\f=\R\cup\{\f^{-1}(\be)\}$ and we have
$$
\|{\mathcal Q}_\f f\|_{\OL(\R)}\le\|{\mathcal Q}_\f f\|_{\OL(\R\cup\{\f^{-1}(\be)\})}\le3\|f\|_{\OL(\T)}
$$
for every $f\in\OL_a(\T)$. $\bl$

Put $\OL'(\R)\df\{f':f\in\OL(\R)\}$ and $\|f'\|_{\OL'(\R)}\df\|f\|_{\OL(\R)}$.
Then $\OL'(\R)$ is a Banach space of functions on $\widehat\R$.

\begin{thm}
\label{dxa}
Let $f\in\OL(\R)$. Then $(x-a)^{-1}(f(x)-f(a))\in(\OL)'(\R)$ and
$$
\big\|(x-a)^{-1}(f(x)-f(a))\big\|_{(\OL)'(\R)}\le\|f\|_{\OL(\R)}\quad
\mbox{for every}\quad a\in\R.
$$
\end{thm}

\Pf It suffices to consider the case $a=0$ and $f(0)=0$. Put
$$
F(x)=\int_0^x\frac{f(t)}t\,dt=\int_0^1\frac{f(tx)}t\,dt.
$$
We have to prove that $F\in\OL(\R)$ and $\|F\|_{\OL(\R)}\le\|f\|_{\OL(\R)}$.
Note that for every $t$ in $(0,1]$, the function $x\mapsto t^{-1}{f(tx)}$ belongs to $\OL_0(\R)$ (see \rf{3101}) and 
$\|t^{-1}{f(tx)}\|_{\OL(\R)}=\|f\|_{\OL(\R)}$
for every $t$ in $(0,1]$. Consequently,
$$
\|F\|_{\OL(\R)}\le\int_0^1\|t^{-1}f(tx)\|_{\OL(\R)}\,dt=\|f\|_{\OL(\R)}. \quad\bl
$$

\medskip

{\bf Remark.} One can prove in a similar way that for a closed nondegenerate interval $J$ and for an arbitrary function $f$ in $\OL(J)$ we have
$$
\big\|(x-a)^{-1}(f(x)-f(a))\big\|_{(\OL)'(J)}\le\|f\|_{\OL(J)}
\quad\mbox{for every}\quad a\in J,
$$
where $\OL'(J)\df\{g':g\in\OL(J)\}$ and $\|g'\|_{\OL'(J)}\df\|g\|_{\OL(J)}$.

\begin{thm}
\label{delab}
If $f\in\OL(\R)$, then $(x-a-b{\rm i})^{-1}(f(x)-f(a))\in(\OL)'(\R)$ and
$$
\big\|(x-a-b{\rm i})^{-1}(f(x)-f(a))\big\|_{(\OL)'(\R)}\le2\|f\|_{\OL(\R)} 
\quad\mbox{for all}\quad a,b\in\R.
$$
\end{thm}

\Pf It suffices to consider the case when $a=0$, $b=1$, $f(0)=0$ and $\|f\|_{\OL(\R)}=1$.  
It follows from Theorem \ref{dxa} that 
$$
\big\|(x-{\rm i})^{-1}f(x)\big\|_{(\OL)'(\R)}
\le\big\|(x-{\rm i})^{-1}{xf(x)}\big\|_{\OL(\R)}.
$$
It remains to prove that $\big\|(x-{\rm i})^{-1}{xf(x)}\big\|_{\OL(\R)}\le2$. Let
$A$ and $B$ be self-adjoint operators. We have
\begin{align*}
A(A-{\rm i}I)^{-1}f(A)-B(B-{\rm i}I)^{-1}f(B)
&=A(A-{\rm i}I)^{-1}(f(A)-f(B))\\
&+
{\rm i}(A-{\rm i}I)^{-1}(B-A)(B-{\rm i}I)^{-1}f(B),
\end{align*}
whence
$$
\|A(A-{\rm i}I)^{-1}f(A)-B(B-{\rm i}I)^{-1}f(B)\|\\
\le\|f(A)-f(B)\|+\|A-B\|\cdot\|g(B)\|\le2\|A-B\|,
$$
where $g(t)=(t-{\rm i})^{-1}f(t)$. $\bl$

\begin{cor}
\label{cordelab}
Let $f\in\OL(\R)$. Then $(x-a-b{\rm i})^{-1}f(x)\in(\OL)'(\R)$ and
$$
\big\|(x-a-b{\rm i})^{-1}f(x)\|_{(\OL)'(\R)}
\le\big(2+|f(a)|/|b|\big)\|f\|_{\OL(\R)} 
\quad\mbox{for}\quad a,b\in\R,~\;b\ne0.
$$
\end{cor}

\Pf  We can assume that $a=0$, $b=1$ and $\|f\|_{\OL(\R)}=1$. 
Using Theorem \ref{delab} and Example 2 in  \S\:\ref{prim}, we obtain
$$
\left\|\dfrac{f(x)}{x-{\rm i}}\right\|_{(\OL)'(\R)}\le\left\|\dfrac{f(x)-f(0)}{x-{\rm i}}\right\|_{(\OL)'(\R)}
+|f(0)|\cdot\|(x-{\rm i})^{-1}\|_{(\OL)'(\R)}\le2+|f(0)|.\quad\bl
$$

\begin{thm}
\label{rinv}
Let $h\in\OL'(\R)$. Then $h\circ\f\in\OL'(\R)$ for all $\f\in\dlR$ and
$
\frac19\|h\|_{\OL'(\R)}\le\|h\circ\f\|_{\OL'(\R)}\le9\|h\|_{\OL'(\R)}.
$
\end{thm}

\Pf The result is obvious if $\f\in\lR$. In this case
$\|h\|_{\OL'(\R)}=\lb\|h\circ\f\|_{\OL'(\R)}=\|h\|_{\OL'(\R)}$. Thus, everything reduces to the case
$\f(t)=\phi(t)\df t^{-1}$.
Let $h=f'$ for a function $f\in\OL(\R)$ such that $f(0)=0$ and 
$\|f\|_{\OL(\R)}=\|h\|_{\OL'(\R)}$.
It follows from Theorem \ref{zamf} that $\|x^2f(x^{-1})\|_{\OL(\R)}\le3\|h\|_{\OL'(\R)}$, whence
$$
\|(x^2f(x^{-1}))'\|_{\OL'(\R)}=\|2xf(x^{-1})-h(x^{-1})\|_{\OL'(\R)}\le3\|h\|_{\OL'(\R)}.
$$
Theorem \ref{dxa} implies the following inequality:
$$
\|xf(x^{-1})\|_{\OL'(\R)}\le\|x^2f(x^{-1})\|_{\OL(\R)}\le3\|h\|_{\OL'(\R)}.
$$
Hence,
$$
\|h(x^{-1})\|_{\OL'(\R)}\le\|(x^2f(x^{-1}))'\|_{\OL'(\R)}+2\|xf(x^{-1})\|_{\OL'(\R)}\le9\|h\|_{\OL'(\R)}.
$$
Applying this inequality to $h(x^{-1})$, we obtain
$
\frac19\|h(x^{-1})\|_{\OL'(\R)}\le\|h\|_{\OL'(\R)}.$ $\bl$

\medskip

\section{\bf The spaces  $\bs{\OL(\R)}$ and $\bs{\OL(\T)}$}
\setcounter{equation}{0}
\label{olriort}

\medskip

The main purpose of this and the next sections is to ``transplant''
Theorem \ref{rinv} from the line to the circle.

It is easy to see that if $f\in\OL(\T)$, then $f(e^{{\rm i}t})\in\OL(\R)$ and
$\|f(e^{{\rm i}t})\|_{\OL(\R)}\le\|f\|_{\OL(\T)}$. 
We show here that the converse also holds, i.e.,
each $(2\pi)$-periodic function $F$ in $\OL(\R)$ can be represented in the form $F=f(e^{{\rm i}t})$, where
$f\in\OL(\T)$ and $\|f\|_{\OL(\T)}\le\const\|F\|_{\OL(\R)}$. This can be deduced 
easily from Lemma 9.8 of \cite{AP5}, see also Lemma 5.7 of \cite{A1}.

\begin{lem}
\label{111}
Let $h(x,y)=\dfrac{x-y}{e^{{\rm i}x}- e^{{\rm i}y}}$. Then
$
\|h\|_{\fM(J_1\times J_2)}\le\frac{3\sqrt2\pi}4
$
for all intervals $J_1$ and $J_2$ such that 
$J_1-J_2\subset[-\frac32\pi,\frac32\pi]$.
\end{lem}

\Pf Consider the $3\pi$-periodic function $\xi$ such that 
$\xi(t)=t(2\sin(t/2))^{-1}$
for $t\in[-\frac32\pi,\frac32\pi]$. Then
$$
\|h\|_{\fM(J_1\times J_2)}=\big\|e^{\frac{{\rm i}x}2}h(x,y)e^{\frac{{\rm i}y}2}\big\|_{{\fM}(J_1\times J_2)}
=\|\xi(x-y)\|_{{\fM}(J_1\times J_2)}
$$
because $x-y\in[-\frac32\pi,\frac32\pi]$ for $x\in J_1$ and $y\in J_2$. Let us expand the function $\xi$ in the Fourier series:
$$
\xi(t)=\sum_{n\in\Z}a_ne^{\frac23n{\rm i}t}=a_0+2\sum_{n=1}^\be a_n\cos\frac23nt
$$
because $a_n=a_{-n}$ for every $n\in\Z$. Clearly, $a_0>0$. Note that the function $\xi$ is convex on $[-\frac32\pi,\frac32\pi]$. It follows that $(-1)^na_n\ge0$
for all natural $n$, see Theorem 35 in \cite{HR}.
It remains to observe that
\begin{align*}
\|\xi(x-y)\|_{{\fM}(J_1\times J_2)}&\le\|\xi(x-y)\|_{{\fM}(\R\times\R)}
\le\sum_{n\in\Z}|a_n|\cdot\big\|e^{\frac23n{\rm i}x}e^{-\frac23n{\rm i}y}\big\|_{{\fM}(\R\times\R)}\\
&=\sum_{n\in\Z}|a_n|=\xi\left(\frac{3\pi}2\right)=\frac{3\sqrt2\,\pi}4.\quad\bl
\end{align*}

\begin{thm}
\label{2pper}
Let $f$ be a continuous function on $\T$.
Then
$$
\|f(e^{{\rm i}x})\|_{\OL(\R)}\le\|f\|_{\OL(\T)}\le\frac{3\sqrt2\pi}2\|f(e^{{\rm i}x})\|_{\OL(\R)}.
$$
\end{thm}

\Pf As observed above, the first inequality is obvious. Let us prove the second one. Put $g(x)\df f(e^{{\rm i}x})$. We may assume that $\|g\|_{\OL(\R)}<\be$.
Then
$g$ is differentiable everywhere on $\R$. It follows that
$f$ is differentiable everywhere on $\T$. It follows from Theorems \ref{su}  and \ref{muld} that
$$
\|g\|_{\OL(\R)}=\|\dg g\|_{\fM(\R\times\R)}~{\text and}\,
~\|f\|_{\OL(\T)}=\|\dg f\|_{\fM(\T\times\T)}\!=\!\|(\dg f)(e^{{\rm i}x},e^{{\rm i}y})\|_{\fM([0,2\pi)\times[-\frac\pi2,\frac{3\pi}2))}.
$$
Therefore, we have to prove that
$$
\Big\|(\dg f)(e^{{\rm i}x},e^{{\rm i}y})\Big\|_{\fM([0,2\pi)\times[-\frac\pi2,\frac{3\pi}2))}
\le\frac{3\sqrt2\pi}2\|\dg g\|_{\fM(\R\times\R)}.
$$

We denote by $\chi_{jk}$ the characteristic function of
$\bs{J}_{j,k}\df[j\pi,(j+1)\pi)\times[k\pi-\frac\pi2,k\pi+\frac{\pi}2)$,
where $j,k\in\Z$. Note that
$$
\chi_{jk}(x,y)(\dg f)(e^{{\rm i}x},e^{{\rm i}y})
=\chi_{jk}(x,y)h(x,y)(\dg g)(x,y),
$$
where $h$ denotes the same as in Lemma \ref{111}.
This together with Lemma \ref{111} yields
\bay
\label{n111}
\Big\|(\dg f)(e^{{\rm i}x},e^{{\rm i}y})\Big\|_{\fM(\bs{J}_{j,k})}
\le\frac{3\sqrt2\pi}4\|\dg g\|_{\fM(\R\times\R)}
\ey
for $(j,k)\in\{0,1\}$, $(j,k)\ne(1,0)$.

The case when $j=1$ and $k=0$ should be considered separately, for in this case 
$J_1-J_2\not\subset[-\frac{3\pi}2,\frac{3\pi}2]$ and we cannot apply Lemma \ref{111} directly.

Let now $j=1$ and $k=0$. We have
$$
\chi_{10}(x+2\pi,y)(\dg f)(e^{{\rm i}x},e^{{\rm i}y})
=\chi_{10}(x+2\pi,y)h(x,y)(\dg g)(x,y).
$$
Now, applying Lemma \ref{111}, we obtain
$$
\|(\dg f)(e^{{\rm i}x},e^{{\rm i}y})\|_{\fM(\bs{J}_{1,0})}
=\|(\dg f)(e^{{\rm i}x},e^{{\rm i}y})\|_{\fM(\bs{J}_{-1,0})}\le\frac{3\sqrt2\pi}4\|\dg g\|_{\fM(\R\times\R)}.
$$
Put also $\bs{J}\df[0,2\pi)\times[-\frac\pi2,\frac{3\pi}2)$. Then
\begin{align*}
\big\|(\dg f&)(e^{{\rm i}x},e^{{\rm i}y})\big\|_{\fM(\bs{J})}\le
\big\|(\chi_{00}(x,y)+\chi_{11}(x,y))(\dg f)(e^{{\rm i}x},e^{{\rm i}y})\big\|_{\fM(\bs{J})}\\
&+
\big\|(\chi_{01}(x,y)+\chi_{10}(x,y))(\dg f)(e^{{\rm i}x},e^{{\rm i}y})\big\|_{\fM(\bs{J})}\\
&\le\max\big\{\big\|(\dg f)(e^{{\rm i}x},e^{{\rm i}y})\big\|_{\fM(\bs{J}_{0,0})},
\big\|(\dg f)(e^{{\rm i}x},e^{{\rm i}y})\big\|_{\fM(\bs{J}_{1,1})}\big\}\\
&+
\max\big\{\big\|(\dg f)(e^{{\rm i}x},e^{{\rm i}y})\big\|_{\fM(\bs{J}_{0,1})},
\big\|(\dg f)(e^{{\rm i}x},e^{{\rm i}y})\big\|_{\fM(\bs{J}_{1,0})}\big\}
\le\frac{3\sqrt2\pi}2\|\dg g\|_{\fM(\R\times\R)}.~\;\bl
\end{align*}

\medskip

{\bf Remark.} It follows from the proof of the theorem that 
$$
\|f(e^{{\rm i}x})\|_{\OL(\R)}\le\|f\|_{\OL(\T)}\le\frac{3\sqrt2\pi}2\|f(e^{{\rm i}x})\|_{\OL(J)}
$$
for every $f$ in $C(\T)$, where $J$ is an interval of length $3\pi$. 

\medskip

\section{\bf The spaces $\bs{(\OL)'(\R)}$ and $\bs{(\OL)_{\rm loc}'(\T)}$}
\setcounter{equation}{0}
\label{olriolt+}

\medskip

The space $(\OL)'(\R)$ has been defined in \S\:\ref{dlp}. We define the space $(\OL)_{\rm loc}'(\T)$ by
$$
(\OL)_{\rm loc}'(\T)\df\left\{f: f(e^{{\rm i}t})\in(\OL)'(\R)\right\}\quad\text{and}\quad
\|f\|_{(\OL)_{\rm loc}'(\T)}\df\left\|f(e^{{\rm i}t})\right\|_{(\OL)'(\R)}.
$$
Note that
$
\|f\|_{L^\be(\T)}=\|f(e^{{\rm i}t})\|_{L^\be(\R)}\le\|f(e^{{\rm i}t})\|_{(\OL)'(\R)}
=\|f\|_{(\OL)_{\rm loc}'(\T)}.
$

We need the following elementary lemma. 

\begin{lem} 
\label{l132}
Let $f\in\Li(\T)$. Then $f\in(\OL)_{\rm loc}'(\T)$ and
$$
\|f\|_{(\OL)_{\rm loc}'(\T)}\le|\widehat f(0)|+\frac\pi{\sqrt3}\|f\|_{\Li(\T)}.
$$
\end{lem}

\Pf Note that $\|f'\|_{L^2(\T)}\le\|f'\|_{L^\be(\T)}\le\|f\|_{\Li(\T)}$ and $\|z^n\|_{(\OL)'(\T)}=1$
for every $n\in\Z$. Consequently,
\begin{align*}
\|f\|_{(\OL)_{\rm loc}'(\T)}&\le\sum_{n\in\Z}|\widehat f(n)|\le|\widehat f(0)|+\Big(\sum_{n\ne0}n^2|\widehat f(n)|^2\Big)^{\frac12}
\Big(\sum_{n\ne0}\frac1{n^2}\Big)^{\frac12}\\
&=|\widehat f(0)|\!+\!\frac\pi{\sqrt3}\|f'\|_{L^2(\T)}\le|\widehat f(0)|\!+\!
\frac\pi{\sqrt3}\|f'\|_{L^\be(\T)}
\le|\widehat f(0)|\!+\!\frac\pi{\sqrt3}\|f\|_{\Li(\T)}. ~\bl
\end{align*}

\begin{cor}
\label{132}
The space $\OL(\T)$ is contained in  $(\OL)_{\rm loc}'(\T)$ and
$$
\|f\|_{(\OL)_{\rm loc}'(\T)}\le|\widehat f(0)|+\frac\pi{\sqrt3}\|f\|_{\OL(\T)}.
$$
\end{cor}

{\bf Remark.}  One can see from the proof of Lemma \ref{l132} that
$$
\|f\|_{(\OL)_{\rm loc}'(\T)}\le|\widehat f(0)|+\frac\pi{\sqrt3}\|f'\|_{L^2(\T)}.
$$

\begin{thm}
\label{rt}
If  $f\in\OL(\T)$, then $zf'(z)\in(\OL)_{\rm loc}'(\T)$  and
$\|zf'(z)\|_{(\OL)_{\rm loc}'(\T)}\le\|f\|_{\OL(\T)}$.
If $f\in(\OL)_{\rm loc}'(\T) $ and $\int_\T f(z)\,d\m(z)=0$, then there exists a function $F$ in $\OL(\T)$ such that $zF'(z)=f$ and 
$\|F\|_{\OL(\T)}\le\const\|f\|_{(\OL)_{\rm loc}'(\T)}$.
\end{thm}

\Pf The first statement is obvious because if $f\in\OL(\T)$, then
$$
\int_0^x e^{{\rm i}t}f'(e^{{\rm i}t})\,dt={\rm i}f(1)-{\rm i}f(e^{{\rm i}x})
$$
and $\|f'\|_{(\OL)_{\rm loc}'(\T)}=\|f(e^{{\rm i}x})\|_{\OL(\R)}\le\|f\|_{\OL(\T)}$.

Let us prove the second statement. Put
 $F(e^{{\rm i}x})\df{\rm i}\int_0^x f(e^{{\rm i}t})\,dt$. $F$ is well defined because $\int_0^{2\pi} f(e^{{\rm i}t})\,dt=2\pi\int_\T f(z)\,d\m(z)=0$.
Clearly, $zF'(z)=f(z)$. It remains to observe that $\|f\|_{(\OL)_{\rm loc}'(\T)}=\|F(e^{{\rm i}x})\|_{\OL(\R)}$
and apply Theorem \ref{2pper}. $\bl$

 \begin{cor}
 \label{rt2}
 A function $f$ on $\T$ belongs to
 $(\OL)_{\rm loc}'(\T)$ if and only if it can be represented in the form
 $f=\widehat f(0)+zF'(z)$, where $F\in\OL(\T)$. Herewith
 $$
 \|f\|_{(\OL)_{\rm loc}'(\T)}
 \le|\widehat f(0)|+\|F\|_{\OL(\T)}\le\const\|f\|_{(\OL)_{\rm loc}'(\T)}.
 $$
 \end{cor}

 \Pf It is easy to see that $\|1\|_{(\OL)_{\rm loc}'(\T)}=1$. This together with Theorem \ref{rt} implies that if
 $f=\widehat f(0)+zF'(z)$ for a function $F$ in $\OL(\T)$,
 then $f\in(\OL)_{\rm loc}'(\T)$ and
 \begin{align*}
\|f\|_{(\OL)_{\rm loc}'(\T)}&\le\|\widehat f(0)+zF'(z)\|_{(\OL)_{\rm loc}'(\T)}\le|\widehat f(0)|+\|zF'(z)\|_{(\OL)_{\rm loc}'(\T)}\\
&\le\|f\|_{(\OL)_{\rm loc}'(\T)}+
\|F\|_{\OL(\T)}\le c\|f\|_{(\OL)_{\rm loc}'(\T)}.
 \end{align*}
Let $f\in(\OL)_{\rm loc}'(\T)$. 
Then by Theorem \ref{rt}, the function $f-\widehat f(0)$
can be represented in the form $f-\widehat f(0)=zF'(z)$, where $F\in\OL(\T)$. $\bl$

\begin{cor}
\label{rt3}
If $f\in(\OL)_{\rm loc}'(\T)$, then $z^nf(z)\in(\OL)_{\rm loc}'(\T)$ for every $n$ in $\Z$.
\end{cor}

\Pf It suffices to consider the case when $f=zF'(z)$, where $F\in\OL(\T)$.
Then
$$
z^nf(z)=z^{n+1}F'(z)=z(z^nF(z))'-nz^nF(z)\in(\OL)_{\rm loc}'(\T)
$$
because $z^nF(z)\in\OL(\T)$ and $\OL(\T)\subset(\OL)_{\rm loc}'(\T)$ by Corollary \ref{132}. $\bl$

 \begin{cor}
 \label{rt4}
 A function $f$ on $\T$ belongs to
 $(\OL)_{\rm loc}'(\T)$ if and only if it can be represented in the form
 $f=\widehat f(-1)z^{-1}+F'(z)$, where $F\in\OL(\T)$; herewith
 $$
 \|f\|_{(\OL)_{\rm loc}'(\T)}
 \le|\widehat f(-1)|+\|F\|_{\OL(\T)}\le\const\|f\|_{(\OL)_{\rm loc}'(\T)}.
 $$
 \end{cor}

 \Pf Put $g(z)\df zf(z)$.  Then $\widehat f(-1)=\widehat g(0)$ and
 $g(z)=\widehat g(0)+zF'(z)$. It remains to refer to Corollaries \ref{rt2} and  \ref{rt3}. $\bl$
 
The following assertion is obvious.

\begin{lem}
\label{mol}
Let $f,g\in\OL(J)$, where $J$ is a bounded closed
interval of $\R$. Then $fg\in\OL(J)$ and 
$$
\|fg\|_{\OL(J)}\le\Big(\m(J)\|g\|_{\OL(J)}+\max_{J}|g|\Big)\|f\|_{\OL(J)}.
$$
\end{lem}

\begin{lem}
\label{zetdzeta}
Let $f\in\OL(\T)$ and $\z\in\T$. Then 
$\dfrac{f(z)-f(\z)}{z-\z}\in(\OL)_{\rm loc}'(\T)$ and
$$
\left\|\dfrac{f(z)-f(\z)}{z-\z}\right\|_{(\OL)_{\rm loc}'(\T)}
\le\const\|f\|_{\OL(\T)}.
$$
\end{lem}

\Pf It suffices to consider the case $\z=1$. We may assume that $f(1)=0$.
We have to estimate the $\OL(\R)$-seminorm of the function $\Phi$,
$$
\Phi(x)\df\int_0^x\frac{f(e^{{\rm i}t})}{e^{{\rm i}t}-1}\,dt.
$$
Clearly, $\Phi$ can be represented in the form $\Phi(x)=\l x+\Phi_0(x)$,
where $\Phi_0$ is a $2\pi$-periodic function. We have
$$
|\l|=\left|\frac1{2\pi}\int_0^{2\pi}\frac{f(e^{{\rm i}t})}{e^{{\rm i}t}-1}\,dt\right|\le
\frac1{2\pi}\int_0^{2\pi}\frac{|f(e^{{\rm i}t})-f(1)|}{|e^{{\rm i}t}-1|}\,dt\le\|f\|_{\Li(\T)}\le\|f\|_{\OL(\T)}.
$$
Therefore, it remains to estimate the $\OL(\R)$-seminorm of $\Phi_0$.

We first estimate $\|\Phi\|_{\OL([-\frac{3\pi}2,\frac{3\pi}2])}$. By the remark to Theorem \ref{dxa} and by Lemma \ref{mol},
$$
\|\Phi\|_{\OL([-\frac{3\pi}2,\frac{3\pi}2])}\le\left\|\frac{tf(e^{{\rm i}t})}{e^{{\rm i}t}-1}\right\|_{\OL([-\frac{3\pi}2,\frac{3\pi}2])}
\le\const\|f(e^{{\rm i}t})\|_{\OL([-\frac{3\pi}2,\frac{3\pi}2])}
\le\const\|f\|_{\OL(\T)}
$$
because the function $t\mapsto\frac{t}{e^{{\rm i}t}-1}$ is infinitely differentiable on $[-\frac{3\pi}2,\frac{3\pi}2]$. Hence, 
$$
\|\Phi_0\|_{\OL([-\frac{3\pi}2,\frac{3\pi}2])}\le\const\|f\|_{\OL(\T)}.
$$
Using now the remark to Theorem \ref{2pper}, we obtain
$$
\|\Phi_0\|_{\OL(\R)}\le\frac{3\sqrt2\pi}2\|\Phi_0\|_{\OL([-\frac{3\pi}2,\frac{3\pi}2])}\le\const\|f\|_{\OL(\T)}.\quad\bl
$$

\begin{thm}
\label{pere}
Let $f$ be a function on $\T$ and let $\psi$ be a linear-fractional transformation
such that $\psi(\widehat\R)=\T$.
Then $f\in(\OL)_{\rm loc}'(\T)$
if and only if $f\circ\psi\in(\OL)'(\R)$. Moreover,
\bay
\label{0310}
c_1\|f\|_{(\OL)_{\rm loc}'(\T)}
\le\|f\circ\psi\|_{(\OL)'(\R)}\le c_2\|f\|_{(\OL)_{\rm loc}'(\T)},
\ey
where $c_1$ and $c_2$ are absolute positive constants.
\end{thm}

\Pf Put $a=\psi^{-1}(0)$. It is easy to see that $a\in\C\setminus\R$ and $\psi(z)=\z(z-\ov a)^{-1}(z-a)$
for all $z\in\widehat\C$, where $|\z|=1$. Without loss of generality we may assume that
$\z=1$. Let us first prove the second inequality. Let $f\in(\OL)_{\rm loc}'(\T)$.    By Corollary \ref{rt4}, $f$ can be represented in the form $f(z)=\widehat f(-1)z^{-1}+F'(z)$, where
$F\in\OL(\T)$ and 
$|\widehat f(-1)|+\|F\|_{\OL(\T)}\le c\|f\|_{(\OL)_{\rm loc}'(\T)}$.
We have
\begin{align*}
\|f\circ\psi\|_{(\OL)'(\R)}&=
\big\|\widehat f(-1)(1/\psi)+F'\circ\psi\big\|_{(\OL)'(\R)}\\
&\le c\left\|1/\psi\right\|_{(\OL)'(\R)}\|f\|_{(\OL)_{\rm loc}'(\T)}+\|F'\circ\psi\big\|_{(\OL)'(\R)}.
\end{align*}
Note that
$$
\|1/\psi\|_{(\OL)'(\R)}=\|(t-a)^{-1}(t-\ov a)\|_{(\OL)'(\R)}\le1+2|\im a|\cdot\|(t-a)^{-1}\|_{(\OL)'(\R)}
\le3
$$
which follows easily from Example 2 in \S\:\ref{prim}.
Let us estimate $\|F'\circ\psi\big\|_{(\OL)'(\R)}$.
We select a function $F$ so that
$F(1)=F(\psi(\be))=0$. It follows from Theorem \ref{zamrt} that
$$
\big\|(F\circ\psi)/\psi'\big\|_{\OL(\R)}=\|{\mathcal Q}_\psi F\|_{\OL(\R)}
\le3\|F\|_{\OL(\T)}\le\const\|f\|_{(\OL)_{\rm loc}'(\T)}.
$$
Consequently,
$$
\big\|F'\circ\psi-(F\circ\psi)\psi''/(\psi')^2\big\|_{(\OL)'(\R)}=
\big\|\big((F\circ\psi)/\psi'\big)'\big\|_{(\OL)'(\R)}\le\const\|f\|_{(\OL)_{\rm loc}'(\T)}.
$$
It remains to estimate
$$
\big\|(F\circ\psi)\psi''/(\psi')^2\big\|_{(\OL)'(\R)}=
\big\|({\mathcal Q}_\psi F)\psi''/\psi'\big\|_{(\OL)'(\R)}
=2\big\|(z-\ov a)^{-1}{\mathcal Q}_\psi F\big\|_{(\OL)'(\R)}.
$$
Using Theorem \ref{delab}, we obtain
\begin{align*}
\big\|(z-\ov a)^{-1}&{\mathcal Q}_\psi F\big\|_{(\OL)'(\R)}
\le\left\|(z-\ov a)^{-1}({\mathcal Q}_\psi F-({\mathcal Q}_\psi F)(\re a))\right\|_{(\OL)'(\R)}\\
&+|({\mathcal Q}_\psi F)(\re a)|\cdot\|(z-\ov a)^{-1}\|_{(\OL)'(\R)}\\
&\le2\|{\mathcal Q}_\psi F\|_{\OL(\R)}+2|F(-1)|\cdot|\im a|\cdot\|(z-\ov a)^{-1}\|_{(\OL)'(\R)}\\
&\le6\|F\|_{\OL(\T)}+2|F(-1)-F(1)|\le10\|F\|_{\OL(\T)}\le\const\|f\|_{(\OL)_{\rm loc}'(\T)}.
\end{align*}

Let us now prove the first inequality.
We can select a function $g\in\OL(\R)$ such that
$g'(t)\df f(\psi(t))\in(\OL)'(\R)$ and $g(\re a)=0$.
Let $\vk$ denote the linear-fractional transformation, which is the inverse of 
$\psi$, i.e., $\vk(z)=(1-z)^{-1}(a-\ov a z)$.
It follows from Theorem \ref{zamfcom} that 
\bay
\label{1-z}
\big\|(2\im a)^{-1}(1-z)^2g(\vk(z))\big\|_{\OL(\T)}\le 5\|g\|_{\OL(\R)}.
\ey
Therefore,  
$$
\big\|(\im a)^{-1}(z-1)g(\vk(z))+f(z)\big\|_{(\OL)_{\rm loc}'(\T)}\le5\|g\|_{\OL(\R)}=5\|f\circ\psi\|_{(\OL)'(\R)}
$$
By Corollary \ref{rt4}.
It remains to prove that 
$$
\big\|(\im a)^{-1}(z-1)g(\vk(z))\big\|_{(\OL)_{\rm loc}'(\T)}\le\const\|f\circ\psi\|_{(\OL)'(\R)}.
$$
This follows immediately from 
\rf{1-z} and from Lemma \ref{zetdzeta}. $\bl$

\begin{thm}
\label{kp}
Let $f\in(\OL)_{\rm loc}'(\T)$ and let $\f$ be a linear-fractional transformation such that $\f(\T)=\T$. Then $f\circ\f\in(\OL)_{\rm loc}'(\T)$
and $c^{-1}\|f\|_{(\OL)_{\rm loc}'(\T)}\le\|f\circ\f\|_{(\OL)_{\rm loc}'(\T)}\le c\|f\|_{(\OL)_{\rm loc}'(\T)}$ for some positive number $c$.
\end{thm}

\Pf  This theorem is substantially an analog for the circle $\T$ of Theorem   \ref{rinv}, which concerns the line $\R$. Theorem \ref{pere} allows us to ``transplant'' Theorem \ref{rinv} from the line $\R$ to the circle $\T$.
$\bl$

\medskip

\section{\bf Around the sufficient condition by Arazy--Barton--Froedman}
\setcounter{equation}{0}
\label{dostol+}

\medskip

We consider in this section a sufficient condition for operator Lipschitzness
for the circle $\T$ that was found by Arazy, Barton and Friedman \cite{ABF} as well as its analog for the line $\R$. Following \cite{A2}, we show how to deduce these
sufficient conditions from Theorem \ref{coint7}. Then we introduce the notion of a Carleson measure in the strong sense and reformulate these sufficient conditions in terms of Carleson measures in the strong sense.
We also show how to deduce from them the sufficient conditions in terms of Besov classes, see \S\:\ref{Dost}. We start with the case of the line.

Put $(\CL)'(\C_+)\df\{g':g\in\CL(\C_+)\}$ and
$\|g'\|_{(\CL)'(\C_+)}=\|g\|_{\CL(\C_+)}$.
Clearly,
$(\CL)'(\C_+)$ is a Banach space.
Note that the functions of class
$(\CL)'(\C_+)$ are defined everywhere on $\clos\C_+\cup\{\be\}$.
It is easy to see that for every $g$ in $\CL(\C_+)$, the Poisson integral of $g'\big|\R$ coincides with 
$g'\big|\C_+$. Indeed, it suffices to observe that for every $t>0$,
the Poisson integral of $t^{-1}(g_t-g)\big|\R$ coincides with $t^{-1}(g_t-g)\big|\C_+$,
where $g_t(z)\df g(z+t)$ and pass to the limit as $t\to 0^+$.
We denote by $(\OL)_+(\R)$ the space of functions $f\in\OL(\R)$ having an analytic extension to the upper half-plane $\C_+$ that is continuous 
in its closure. Put
$(\OL)'_+(\R)\df\{g':g\in\OL_+(\R)\}$.

It follows from Theorem \ref{olclr} that the space $(\OL)'_+(\R)$
can be identified in a natural way with the space 
$(\CL)'(\C_+)$. Moreover, 
$$
\|f'\|_{(\OL)'(\R)}=\|f\|_{\OL(\R)}=\|f\|_{\CL(\C_+)}
=\|f'\|_{(\CL)'(\C_+)}
$$
for every $f\in(\CL)(\C_+)$.

An analog of the Arazy--Barton--Friedman theorem for the half-plane can be stated in the following way: 

\begin{thm}
\label{abfap}
Let $f$ be a function analytic in the upper half-plane. Suppose that
$$
\sup_{t\in\R}\int_{\C_+}\frac{(\im w)|f'(w)|\,d\m_2(w)}{|t-w|^2}<+\be.
$$
Then $f$ has finite angular boundary values\footnote{By $f(\be)$ we undestand $\lim f(z)$ as $|z|\to\be$ and $z$ remains in a closed angle with vertex in $\R$ and such that all its points except for the vertex are in $\C_+$.} 
everywhere on
$\widehat\R$, which will be denoted by the same letter $f$, $f\in(\CL)'(\C_+)$,
and
$$
\|f-f(\be)\|_{(\CL)'(\C_+)}\le\frac2\pi\sup_{t\in\R}\int_{\C_+}
\frac{(\im w)|f'(w)|\,d\m_2(w)}{|t-w|^2}.
$$
\end{thm}

\begin{lem}
\label{berg}
Let $f$ be an analytic function in the upper half-plane $\C_+$. Suppose that
$$
\int_{\C_+}(\im w)(1+|w|^2)^{-1}|f'(w)|\,d\m_2(w)<+\be.
$$
Then $f$ has finite angular value
$f(\be)$ at infinity and
$$
f(z)-f(\be)=\frac{2\rm i}{\pi}\int_{\C_+}\frac{(\im w)f'(w)\,d\m_2(w)}{(\ov w-z)^2}
$$
for every $z\in\C_+$.
\end{lem}

\Pf Put
$$
g(z)\df\frac{2{\rm i}}{\pi}\int_{\C_+}\frac{(\im w)f'(w)\,d\m_2(w)}{(\ov w-z)^2}
$$
for $z\in\C_+$.
Clearly, $g$ is analytic in $\C_+$ and
$$
g'(z)=\frac{4{\rm i}}{\pi}\int_{\C_+}\frac{f'(w)\,d\m_2(w)}{(\ov w-z)^3}=f'(z)
$$
for every $z\in\C_+$. The last equality follows from the fact that 
$4{\rm i}(\pi)^{-1}(\ov w-z)^{-3}$ is the reproducing kernel
for the Bergman space of functions analytic in $\C_+$ that belong to
$L^2(\C_+,y\,d\m_2(x+{\rm i}y))$. This is a well-known and easily verifiable fact. It remains to prove that the nontangential limit of $g$ at infinity is zero.
It follows from the equality
$$
g(z)=\frac{2{\rm i}}{\pi}\int_{\C_+}
\left(\frac{\ov w-{\rm i}}{\ov w-z}\right)^2
\frac{f'(w)\,d\m_2(w)}{(\ov w-{\rm i})^2}
$$
and from the Lebesgue dominated convergence theorem that
the restriction of $g(z)$ to any half-plane 
$\e{\rm i}+\C_+$ with $\e>0$ tends to zero as $|z|\to\be$. $\bl$

\medskip

{\bf Proof of Theorem \ref{abfap}.}
Put
\bey
F(z)\df
\frac{{2\rm i}}\pi\int_{\C_+}(\im w)f'(w)\left(\frac1{\ov w-z}-\frac1{\ov w}\right)\,d\m_2(w)
=\frac{2{\rm i}z}\pi\int_{\C_+}\frac{(\im w)f'(w)\,d\m_2(w)}{\ov w(\ov w-z)}
\eey
for every $z\in\C$ such that $\im z\ge0$.
The convergence of the integrals for real $z$ follows from the Cauchy--Bunyakovsky inequality if we take into account the following obvious inequality:
$$
\int_{\C_+}\frac{(\im w)|f'(w)|\,d\m_2(w)}{|z-\ov w|^2}
\le\int_{\C_+}\frac{(\im w)|f'(w)|\,d\m_2(w)}{|\re z-w|^2}.
$$
Note that
\bay
\label{FprimG}
F'(z)
=\frac{2{\rm i}}{\pi}\int_{\C_+}\frac{(\im w)f'(w)\,d\m_2(w)}{(\ov w-z)^2}
=f(z)-f(\be)
\ey
by Lemma \ref{berg}. Consider the Radon measure $\mu$ in the lower half-plane 
$\C_-$,
$$
d\mu(w)\df\frac{2{\rm i}}{\pi}(\im \ov w)f'(\ov w)\,d\m_2(w).
$$
Then $F(z)=\widehat\mu_{0}(z)$, whenever $\im z\ge0$ and
$$
\|\mu\|_{{\M}(\C_-)}=\frac2\pi\sup_{z\in\C_+}\int_{\C_+}
\frac{(\im w)|f'(w)|\,d\m_2(w)}{|\ov z-w|^2}=\frac2\pi\sup_{t\in\R}\int_{\C_+}
\frac{(\im w)|f'(w)|\,d\m_2(w)}{|t-w|^2}.
$$

It follows now from Theorem \ref{coint7} that
$$
\|f-f(\be)\|_{(\CL)'(\C_+)}=\|F\|_{\CL(\C_+)}\le\frac2\pi\sup_{t\in\R}\int_{\C_+}
\frac{(\im w)|f'(w)|\,d\m_2(w)}{|t-w|^2}.\quad\bl
$$

\medskip

We denote by $\mathcal PM(\C_+)$ the space of harmonic functions
$u$ in the upper half-plane $\C_+$ and such that
$$
\|u\|_{\mathcal PM(\C_+)}\df\sup_{y>0}\int_\R|u(x+{\rm i}y)|\,dx<+\be.
$$
It is well known (see, for example, \cite{SW}, Theorems 2.3 and 2.5 of Chapter II)
that $\mathcal PM(\C_+)$ coincides with the set of functions
$u$ that can be represented in the form
$$
u(z)=(\mathcal P\nu)(z)\df\frac1\pi\int_\R\frac{\im z\,d\nu(t)}{|z-t|^2}, \quad z\in\C_+,
$$
where $\nu$ is a complex Borel measure on $\R$ and
$\|u\|_{\mathcal PM(\C_+)}=\|\nu\|_{M(\R)}\df|\nu|(\R)$.

We denote by $\mathcal PL(\C_+)$ the subspace of
$\mathcal PM(\C_+)$, which consists of the functions $u$
that correspond to absolutely continuous measures $\nu$.

A positive measure $\mu$ on $\C_+$ is called a {\it Carleson measure in the strong sense} if $\int_{\C_+}|u(z)|\,d\mu(z)<\be$ for every $u\in\mathcal PM(\C_+)$.
Note that $\mathcal PM(\C_+)$ contains the Hardy class $H^1$ in the upper half-plane $\C_+$. It follows that a Carleson measure in the strong sense must be a Carleson measure in the usual sense. {\it We denote by 
${\rm CM}_{\rm s}(\C_+)$ the space of all Radon measures $\mu$ in $\C_+$ 
such that $|\mu|$ is a Carleson measure in the strong sense}.
Put
$$
\|\mu\|_{{\rm CM}_{\rm s}(\C_+)}\df\sup\left\{\int_{\C_+}|u(z)|\,d\mu(z):u\in\mathcal PM(\C_+), \|u\|_{\mathcal PM(\C_+)}
\le1\right\}.
$$
It is easy to see that
$$
\|\mu\|_{{\rm CM}_{\rm s}(\C_+)}=\sup\left\{\int_{\C_+}|u(z)|\,d\mu(z):u\in\mathcal PL(\C_+), \|u\|_{\mathcal PM(\C_+)}
\le1\right\}
$$
and
$$
\|\mu\|_{{\rm CM}_{\rm s}(\C_+)}=\frac1\pi\sup_{t\in\R}\int_{\C_+}\frac{(\im w)\, d\mu(w)}{|t-w|^2}
=\frac1\pi\sup_{z\in\C_+}\int_{\C_+}\frac{(\im w)\, d\mu(w)}{|\ov z-w|^2}.
$$

We can state now the analog of the Arazy--Barton--Friedman theorem for the half-plane as follows.

\begin{thm} 
Let $f$ be a function analytic in $\C_+$. Suppose that 
$|f'|\,d\m_2\in{\rm CM}_{\rm s}(\C_+)$.
Then $f$ has finite nontangential boundary values everywhere on $\widehat\R$, 
$f\in(\CL)'(\C_+)$ and
$$
\|f-f(\be)\|_{(\CL)'(\C_+)}\le 2\|f'\,d\m_2\|_{{\rm CM}_{\rm s}(\C_+)},
$$
where the same symbol $f$ is used for the corresponding boundary-value function.
\end{thm}

In a similar way we can obtain one more version of the Arazy--Barton--Friedman theorem.
In the following theorem as well as in the whole section $\zh(\n u)(a)\zh$ denotes the operator norm of the differential  
$d_a u$ of a function $u$ at a point $a$.

\begin{thm}
\label{abfhp}
Let $u$ be a (complex) harmonic function
in $\C_+$. Suppose that $\zh\n u\zh\,d\m_2\in{\rm CM}_{\rm s}(\C_+)$.
Then $u$ has nontangential boundary values everywhere on
$\widehat\R$, $u|\R\in(\OL)'(\R)$ and
$$
\|u-u(\be)\|_{(\OL)'(\R)}\le2\|\zh\n u\zh\,d\m_2\|_{{\rm CM}_{\rm s}(\C_+)}.
$$
\end{thm}

\Pf Consider analytic functions $f$ and $g$ in $\C_+$ such that $f+\ov g=u$.
Then
$$
f^\prime=\frac{\partial u}{\partial z}=\frac12\left(\frac{\partial u}{\partial x}-{\rm i}\frac{\partial u}{\partial y}\right)
\quad\text{and}\quad \ov{g^\prime}=\frac{\partial u}{\partial \ov{z}}=
\frac12\left(\frac{\partial u}{\partial x}+{\rm i}\frac{\partial u}{\partial y}\right).
$$
Put
\bey
F(z)\df
\frac{{2\rm i}}\pi\int_{\C_+}(\im w)f'(w)\left(\frac1{\ov w-z}-\frac1{\ov w}\right)\,d\m_2(w)
=\frac{2{\rm i}z}\pi\int_{\C_+}\frac{(\im w)f'(w)\,d\m_2(w)}{\ov w(\ov w-z)}
\eey
and
\bey
G(z)\df
\frac{{2\rm i}}\pi\int_{\C_+}(\im w)g'(w)\left(\frac1{\ov w-z}-\frac1{\ov w}\right)\,d\m_2(w)
=\frac{2{\rm i}z}\pi\int_{\C_+}\frac{(\im w)g'(w)\,d\m_2(w)}{\ov w(\ov w-z)}
\eey
for every $z\in\C$ such that $\im z\ge0$.
Applying identity \rf{FprimG} and the same identity for $G$, we obtain
\begin{align*}
u(x)-u(\infty)&=F'(x)+\ov G'(x)\\
&=\frac{2{\rm i}}{\pi}\int_{\C_+}\frac{(\im w)f'(w)\,d\m_2(w)}{(\ov w-x)^2}
+\frac{2{\rm i}}{\pi}\int_{\C_+}\frac{(\im w)\ov{g'(w)}\,d\m_2(w)}{(w-x)^2}.
\end{align*}
Applying Theorem \ref{coint7}, we obtain
\begin{multline*}
\|u-u(\be)\|_{(\OL)'(\R)}\le\\
\frac2\pi\sup_{x\in\R}\left(\int_{\C_+}
\frac{(\im w)|f'(w)|\,d\m_2(w)}{|x-\ov w|^2}+\int_{\C_+}
\frac{(\im w)|g'(w)|\,d\m_2(w)}{|x-w|^2}\right)\\
=\sup_{x\in\R}\int_{\C_+}\frac{(\im w)(|f'(w)|+|g'(w)|)\,d\m_2(w)}{|x-w|^2}.
\end{multline*}
It remains to observe that $|f'(w)|+|g'(w)|=\zh (\n u)(w)\zh$ for every 
$w$ in $\C_+$
because the operator norm of the linear map $h\mapsto\a h+\b\ov h$
equals $|\a|+|\b|$. $\bl$

\begin{cor} 
Let $f\in\Li(\R)$.
Suppose that $\zh{\rm Hess}\, \mP f\,\zh\,d\m_2\in{\rm CM}_{\rm s}(\C_+)$.
Then $f\in\OL(\R)$.
\end{cor}

Let us now show that the Arazy--Barton--Friedman sufficient condition  implies the sufficient condition for operator Lipschitzness obtained in \cite{Pe1} and \cite{Pe3}, see Theorem \ref{Besdost} in this survey.

To obtain this sufficient condition, we need the elementary inequality:
\bay
\label{51}
\|\f\,d\m_2\|_{{\rm CM}_{\rm s}(\C_+)}\le\int_0^\be\ess\sup\{\f(x+{\rm i}y):x\in\R\}\,dy
\ey
for an arbitrary measurable nonnegative function $\f$ on $\C_+$.

We proceed now to the alternative proof of the sufficient condition obtained in \cite{Pe1}, \cite{Pe3}.

\begin{thm} 
Let $f\in B^1_{\be 1}(\R)$. Then $f\in\OL(\R)$.
\end{thm}

\Pf Clearly, $f'\in L^\be(\R)$.
Let $u$ be the Poisson integral of $f'$. It is well known, see \S\:\ref{Prel},
that the membership $f\in B^1_{\be 1}(\R)$ is equivalent to the following condition:
$$
\int_0^\be \sup_{x\in\R}\zh\n u(x+{\rm i}y)\zh\,dy<+\be.
$$
It remains to use inequality \rf{51} and refer to Theorem \ref{abfhp}. $\bl$

\medskip

Consider now the case of the disk. Put $(\CL)'(\dd)\df\{g':~g\in\CL(\C_+)\}$ and
\lb$\|g'\|_{(\CL)'(\dd)}=\|g\|_{\CL(\dd)}$.

We denote by $\mathcal PM(\dd)$ the space of complex harmonic functions
$u$ in $\dd$ such that
$$
\|u\|_{\mathcal PM(\dd)}\df\sup_{0\le r<1}\int_\T|u(r\z)|\,|d\z|<+\be.
$$
It is well known (see, for example, \cite{H}, Chapter 3)
that the space $\mathcal PM(\dd)$ coincides with the set of functions
$u$ representable in the form
$$
u(z)=(\mathcal P\nu)(z)\df\frac1{2\pi}\int_\T\frac{(1-|z|^2)\,d\nu(\z)}{|z-\z|^2}, \quad z\in\dd,
$$
where $\nu$ is a complex Borel measure on $\T$ and
$\|u\|_{\mathcal PM(\dd)}=\|\nu\|_{M(\T)}\df|\nu|(\T)$.

We denote by $\mathcal PL(\dd)$ the subspace of
$\mathcal PM(\dd)$ that consists of the functions $u$
that correspond to absolutely continuous measures $\nu$.

A positive measure $\mu$ on $\dd$ is called a {\it a Carleson measure in the strong sense} if $\int_{\dd}|u(z)|\,d\mu(z)<\be$ for every $u\in\mathcal PM(\dd)$.
Note that the space $\mathcal PM(\dd)$ contains the Hardy class $H^1$ in the unit disk $\dd$. It follows that a Carleson measure in the strong sense
is a Carleson measure in the usual sense. {\it We denote by ${\rm CM}_{\rm s}(\dd)$ the space of Radon measures $\mu$ in $\dd$
such that $|\mu|$ a Carleson measure in the strong sense}.
Put
$$
\|\mu\|_{{\rm CM}_{\rm s}(\dd)}\df\sup\left\{\int_{\dd}|u(z)|\,d\mu(z):u\in\mathcal PM(\dd), \|u\|_{\mathcal PM(\dd)}
\le1\right\}.
$$
It is easy to see that
$$
\|\mu\|_{{\rm CM}_{\rm s}(\dd)}=\sup\left\{\int_{\dd}|u(z)|\,d\mu(z):u\in\mathcal PL(\dd), \|u\|_{\mathcal PM(\dd)}
\le1\right\}.
$$
and
$$
\|\mu\|_{{\rm CM}_{\rm s}(\dd)}=\frac1{2\pi}\sup_{\z\in\T}\int_{\dd}\frac{(1-|w|^2)\, d\mu(w)}{|\z-w|^2}.
$$
Note that
\bay
\label{pm}
\sup_{\z\in\T}\int_{\dd}\frac{(1-|w|^2)\, d\mu(w)}{|\z-w|^2}
=\sup_{z\in\dd}\int_{\dd}\frac{(1-|w|^2)\, d\mu(w)}{|1-z\ov w|^2}.
\ey
This follows from the maximum principle for $L^2$-valued analytic function in
$\dd$.

Let us now state in our notation the Arazy--Barton--Friedman sufficient condition in the case of the circle, see \cite{ABF}.

\begin{thm}
\label{abfad}
Let $f$ be an analytic function in $\dd$. Suppose that
$\z^{-1}f'(\z)\,d\m_2(\z)\in{\rm CM}_{\rm s}(\dd)$.
Then $f$ has finite nontangential boundary values everywhere on $\T$, 
$f\in(\CL)'(\dd)$ and
$$
\|f-f(0)\|_{(\CL)'(\dd)}\le 2\Big\|\z^{-1}f'(\z)\,d\m_2(\z)\Big\|_{{\rm CM}_{\rm s}(\dd)}.
$$
where the same symbol $f$ is used for the corresponding boundary-value function.
\end{thm}

We need an analog of Lemma \ref{berg}.

\begin{lem}
\label{bergd}
Let $f$ be an analytic function in $\dd$. Suppose that
\lb$\int_{\dd}(1-|w|^2)|f'(w)|\,d\m_2(w)<\be$.
Then 
$$
f(z)-f(0)=\frac{1}{\pi}\int_{\dd}\frac{(1-|w|^2)f'(w)\,d\m_2(w)}{(1-z\ov w)^2\ov w}
=\frac{1}{\pi}\int_{\C\setminus\dd}\frac{(|w|^2-1)f'(\ov w^{-1})\,d\m_2(w)}{(w-z)^2\ov w^3}
$$
for every $z\in\C_+$.
\end{lem}

\Pf We prove only the first equality because the second one can be obtained from the first by the change of variable $w\mapsto\ov w^{-1}$.
For $z=0$, the desired equality follows from the mean value theorem.
It remains to observe that
$$
f'(z)=\frac{2}{\pi}\int_{\dd}\frac{(1-|w|^2)f'(w)\,d\m_2(w)}{(1-z\ov w)^3}
$$
for every $z\in\dd$, see, for example, Corollary 1.5 of the monograph \cite{HKZ}.  $\bl$.

Note also that the first equality of the lemma can be obtained by differentiating
equality (4.3) of \cite{ABF} in $z$.

\medskip

{\bf Proof of Theorem \ref{abfad}.}
Put
\begin{align*}
F(z)&\df
\frac{1}\pi\int_{\dd}\frac{(1-|w|^2)f'(w)}{\ov w^2}\left(\frac1{1-z\ov w}-1\right)\,d\m_2(w)\\
&=\frac{z}\pi\int_{\dd}\frac{(1-|w|^2)f'(w)\,d\m_2(w)}{\ov w(1-z\ov w)}
\end{align*}
for every $z$ in $\C$ such that $|z|\le1$.
The convergence of the integrals for $z\in\T$ is a consequence of the Cauchy--Bunyakovsky inequality if we take into account identity \rf{pm}.
Note that
$$
F'(z)
=\frac{1}{\pi}\int_{\dd}\frac{(1-|w|^2)f'(w)\,d\m_2(w)}{(1-z\ov w)^2\ov w}
=f(z)-f(0)
$$
By Lemma \ref{berg}. Consider the Radon measure $\mu$ in $\C\setminus\ov\dd$,
$$
d\mu(w)\df\frac{1}{\pi}\ov w^{-3}(|w|^2-1)f'(\ov w^{-1})\,d\m_2(w).
$$
Then
\begin{align*}
\widehat\mu_{0}(z)&=\frac1\pi\int_{\C\setminus\dd}\frac{(|w|^2-1)f'(\ov w^{-1})}{\ov w^3}
\left(\frac1{w-z}-\frac1w\right)\,d\m_2(w)\\
&=\frac{1}\pi\int_{\dd}\frac{(1-|w|^2)f'(w)}{\ov w^2}\left(\frac1{1-z\ov w}-1\right)\,d\m_2(w)=F(z).
\end{align*}
Note that
\begin{align*}
\|\mu\|_{{\M}(\C\setminus\clos\dd)}&=\sup_{z\in\clos\dd}\int_{\C\setminus\clos\dd}\frac{d|\mu|(w)}{|w-z|^2}
=\frac1\pi\sup_{z\in\ov\dd}\int_{\C\setminus\clos\dd}\frac{(|w|^2-1)|f'(\ov w^{-1})|}{|w-z|^2|w|^3}\,d\m_2(w)\\
&=\frac1\pi\sup_{z\in\clos\dd}\int_{\dd}\frac{(1-|w|^2)|f'(w)|}{|1-z\ov w|^2|w|}\,d\m_2(w)\\
&=\frac1\pi\sup_{\z\in\T}\int_{\dd}\frac{(1-|w|^2)|f'(w)|}{|\z-w|^2|w|}\,d\m_2(w)
=2\big\|\z^{-1}f'(\z)\,d\m_2(\z)\big\|_{{\rm CM}_{\rm s}(\dd)}.
\end{align*}

It follows now from Theorem \ref{coint7} that
$$
\|f-f(0)\|_{(\CL)'(\dd)}=\|F\|_{\CL(\dd)}\le\|\mu\|_{{\M}(\C\setminus\clos\dd)}
=2\big\|\z^{-1}f'(\z)\,d\m_2(\z)\big\|_{{\rm CM}_{\rm s}(\dd)}.\quad\bl
$$

\begin{cor}
\label{abfadcor}
Let $f$ be an analytic function in $\dd$. Suppose that
$f'\,d\m_2\in{\rm CM}_{\rm s}(\dd)$. Then 
$f$ has finite nontangential boundary values everywhere on $\T$,
$f\in(\CL)'(\dd)$ and
$
\|f-f(0)\|_{(\CL)'(\dd)}\le\const\|f'\,d\m_2\|_{{\rm CM}_{\rm s}(\dd)},
$
where the same symbol $f$ is used for the corresponding boundary-value function.
\end{cor}

\Pf It suffices to observe that for every continuous function $h$ in $\dd$,
the condition $h\,d\m_2\in{\rm CM}_{\rm s}(\dd)$ implies that $\z^{-1}h(\z)\,d\m_2(\z)\in{\rm CM}_{\rm s}(\dd)$.
It remains to apply the closed graph theorem. $\bl$

\medskip

{\bf Remark.} One can obtain without the closed graph theorem the 
explicit estimate
$$
C_s\Big(|\z|^{-1}h(\z)\,d\m_2(\z)\Big)\le\frac83C_s\big(h\,d\m_2\big)
$$
for every subharmonic function $h$ in $\dd$, but we will not need this.

\medskip

\begin{thm}  
If  $\zh \n u\zh\,d\m_2\in{\rm CM}_{\rm s}(\dd)$
for some harmonic function $u$ in $\dd$, then $u$ has nontangential boundary values everywhere on $\T$ and $u\in(\OL)'_{{\rm loc}}(\T)$.
\end{thm}

\Pf 
The function $u$ can be represented as $u=f+\ov g$, where $f$ and $g$ are analytic functions in $\dd$. It follows from Corollary \ref{abfadcor} that $f,g\in(\OL)'(\T)$.
It remains to observe that it follows from the definition of the space $(\OL)'_{{\rm loc}}(\T)$
that it is invariant under complex conjugation. $\bl$

\begin{cor} 
If  $\zh {\rm Hess}\, u\zh\,d\m_2\in{\rm CM}_{\rm s}(\dd)$
for some harmonic function $u$ in $\dd$, then $u$ admits an extension to a continuous function
on $\dd\cup\T$ and
$u\in\OL(\T)$.
\end{cor}

This easily implies the following result of \cite{Pe1} whose proof is given in  
\S\:\ref{Dost} of this survey, see Theorem \ref{Besokr}.

\begin{thm}
\label{612}
Let $f\in B^1_{\be 1}(\T)$. Then $f\in\OL(\T)$.
\end{thm}

We have deduced the Arazy--Barton--Friedman sufficient condition from Theorem \ref{coint7}.
It can be shown that Theorem \ref{coint7} gives examples of operator Lipschitz functions that do not satisfy the analog of the Arazy--Barton--Friedman condition for $\C_+$.
In \cite{A2} an example is given of a function $f$ in $\widehat\M(\clos\C_+)$ 
such that $f''\,d\m_2\not\in{\rm CM}_{\rm s}(\C_+)$. It follows from Theorem \ref{coint7} that such a function $f$ belongs to $(\CL)'(\C_+)$, though 
the Arazy--Barton--Friedman condition fails for this function. A similar assertion is also true for functions in $\dd$.

\medskip

{\bf Remark.} In the the Arazy--Barton--Friedman paper \cite{ABF} it was mentioned that their sufficient condition for the operator Lipschitzness of a function
on the unit circle implies the sufficient condition
obtained in \cite{Pe1}, see also \S \ref{Dost} of this survey. It follows from the results of \cite{A2} that the Arazy--Barton--Friedman sufficient condition
can work even if $f'$ is not continuous. On the other hand, it is easy to see that if $f\in B^1_{\be 1}(\T)$, then $f'\in C(\T)$.
The same can be said about functions in $B^1_{\be 1}(\R)$ (see Theorem \ref{Besdiffer}). Indeed,
it is easy to verify that the function $f(z)=\exp(-{\rm i}z^{-1})$ satisfies the hypotheses of Theorem
\ref{abfap}, though its restriction to the real line is discontinuous at 0. Note also that in \cite{A2} it is proved that a subset if the real line is the set of discontinuity points of $f\big|\R$ for a function
$f$ satisfyng the hypotheses of Theorem \ref{abfap} if and only if it is an 
$F_\s$ set and has no interior points. The same can be said about functions on   $\T$ and $\dd$.

\medskip

It is interesting to compare the sufficient conditions for operator Lipschitzness given in this section with the necessary condition given in \S\:\ref{Neob}. A combination of these conditions is given in the following theorem.

\begin{thm}
If $f\in\Li(\R)$ and $\zh{\rm Hess}\mP f\zh\,d\m_2$ is a Carleson measure in the strong sense, then $f\in\OL(\R)$.
If $f\in\OL(\R)$, then $\zh{\rm Hess}\mP f\zh\,d\m_2$ is a Carleson measure.
\end{thm}

A similar assertion also holds for functions on the circle $\T$.

\medskip

\section{\bf In which case does the equality $\bs{\OL(\fF)=\Li(\fF)}$ holds?}
\setcounter{equation}{0}
\label{orav}

\medskip

\begin{thm}
\label{dal}
Suppose that $\OL(\fF)=\Li(\fF)$ for a closed subset $\fF$ of $\C$. Then  $\fF$ is finite.
\end{thm}

\Pf  Suppose that $\fF$ is infinite. Then 
$\fF$ has a limit point
$a\in\widehat\C\df\C\cup\{\be\}$. If $a\in\C$, we may assume that
$a=0$. The case $a=\be$
 can be considered in the same way. Moreover, the case $a=\be$ reduces to the case
 $a=0$ with the help of linear-fractional transformations.

 Suppose first that $\fF\subset\R$. Then it is easy to construct a function
 $f\in\Li(\fF)$ that has no derivative at 0. Clearly, $f\not\in\OL(\fF)$.

 To get rid of the assumption $\fF\subset\R$, we need the following lemma.

 \begin{lem}
 \label{ldal}
 Let $0<q<1$ and let
 $\{a_n\}_{n\ge1}$ be a sequence of positive numbers such that
 $a_{n+1}\le q a_n$ for every $n\ge1$.
 Then for each numerical sequence $b_n$ satisfying the condition
 $\sum_{n\ge1}|b_n|a^{-1}_n<+\be$, there exists $v\in\OL(\R)$
 such that $v(a_n)=b_n$ for every $n\ge1$.
 \end{lem}

 \Pf Let us fix a function $\f$ of class $C^\be(\R)$ such that $\f(0)=1$ and
 $\supp\f\subset[-\delta,\delta]$, where $\delta$ will be chosen at the end of the proof. Put $v(t)\df\sum\limits_{n\ge1}b_n\f(a_n^{-1}(t-a_n))$.
 Then
 $$
 \|v\|_{\OL(\R)}\le\sum_{n\ge1}|b_n|\cdot\|\f(a_n^{-1}(t-a_n))\|_{\OL(\R)}=\|\f\|_{\OL(\R)}
 \sum_{n\ge1}|b_n|a^{-1}_n<+\be
 $$
 and $v(a_n)=b_n$ for every $n\ge1$ provided $\delta$ is sufficiently small. $\bl$

\medskip

Let us continue the proof of Theorem \ref{dal}. 
It is well known that each Lipschitz function on a subset of $\C$ can be extended to a Lipschitz function on the whole complex plane $\C$, see, for example, \cite{St}, Chapter VI, \S\:2, Theorem 3. 
Thus, it suffices to consider the case when $\fF\setminus\{0\}$ consists 
of the terms of a sequence $\{\l_n\}_{n\ge1}$ that tend to 0 arbitrarily rapidly. Let $\l_n=a_n+{\rm i}b_n$. We may assume that
$\lim\limits_{n\to\be}\frac{\l_n}{|\l_n|}=1$ and
the real sequences $\{a_n\}_{n\ge1}$ and $\{b_n\}_{n\ge1}$
satisfy the hypotheses of Lemma \ref{ldal}. Put $h(t)\df t+{\rm i}v(t)$, where
$v$ means the same as in Lemma \ref{ldal}.
Now the case of the set $\fF$ reduces to the case of the set
$\re\fF$ that has been already treated because 
$$
\|A-B\|\le\|h(A)-h(B)\|\le(1+\|v\|_{\OL(\R)})\|A-B\|
$$
for arbitrary self-adjoint operators $A$ and $B$ such that $\s(A),\:\s(B)\subset\re\fF$. $\bl$

\medskip

\begin{center}
\bf\large Concluding remarks
\end{center}
\label{zaklzam}

\addtocontents{toc}{\vspace*{.03cm}\hspace*{-.35cm}\textbf{Concluding remarks}\hfill\pageref{zaklzam}}

\medskip

In this section we briefly discuss certain results that were not covered in the main part of the paper.

\medskip

{\bf 1. Operator moduli of continuity.} For a continuous function $f$ on $\R$, the operator modulus of continuity $\O_f$ is defined by
$$
\O_f(\d)\df\sup\{\|f(A)-f(B)\|:~A~\mbox{and}~B\;~
\mbox{are\;~self-adjoint\;~operators},~\:\|A-B\|<\d\}.
$$
Operator moduli of continuity were introduced in \cite{AP2}
and were studied in detail in \cite{AP5}. Theorem \ref{prmone} stated in this paper means that if $f\in\L_\o(\R)$, where $\o$ is a modulus of continuity, then
$$
\O_f(\d)\le\const\o_*(\d),\quad\mbox{where}\quad
\o_*(\d)\df\d\int_\d^\be\frac{\o(t)}{t^2}dt.
$$

In \cite{AP5} it was discussed to what extend such estimates are sharp, and considerably sharper estimates were obtained for continuous ``piecewise convex-concave'' functions $f$. In particular, the following best possible estimate was obtained:
$$
\big\|\,|A|-|B|\,\big\|\le C\|A-B\|\log\left(2+\log\frac{\|A\|+\|B\|}{\|A-B\|}\right)
$$
for bounded self-adjoint operators $A$ and $B$. This inequality considerably improves an estimate by Kato obtained in \cite{Kat}.

\medskip

{\bf 2. Commutator estimates for normal operators.} Lemma \ref{predosnnaformKotSh} allows us to obtain the following quasicommutator estimate:
$$
\|f(N_1)R-Rf(N_2)\|\le\const\o_*\big(\max\{\|N_1R-RN_2\|,\|N_1^*R-RN_2^*\|\}\big)
$$
for an arbitrary modulus of continuity $\o$, for every function $f$ of class $\L_\o(\R)$,
for every linear operator $R$ of norm 1, and for arbitrary normal operators $N_1$ and $N_2$ (see \cite{APPS}). In \cite{AP6} the quasicommutator norm on the left-hand side of the inequality was estimated only in terms of the norm $\|N_1R-RN_2\|$. However, on the right-hand side $\o_*$ has to be replaced with $\o_{**}\df(\o_*)_*$:
$$
\|f(N_1)R-Rf(N_2)\|\le\const\o_{**}\big(\{\|N_1R-RN_2\|\}\big).
$$
Note that in the case of H\"older classes, i.e., $\o(t)=t^\a$, $0<\a<1$, we have  $\o_{**}(t)\le\const(1-\a)^{-2}t^\a$. In other words, we obtain a commutator 
H\"older estimate.

\medskip

{\bf3. The Nikolskaya--Farforovskaya approach to operator H\"older functions.}
In \cite{FN} the authors offer an alternative approach to operator H\"older functions. It is based on the following assertion:

\medskip

{\it Let  $0<\a<1$. Then $\L_\a(\Z)\subset\OL(\Z)$.
Moreover, there exists $c_\a$ such that 
$\|f\|_{\OL(\Z)}\le c_\a\|f\|_{\L_\a(\Z)}$.}

\medskip

Theorem \ref{44} can be deduced from this result with the help of the following easily verifiable inequality: $\O_f(\d)\le2\o_f(\d/2)+2\|f(\d x)\|_{\OL(\Z)}$. The assertion itself can be proved by interpolating functions of class $\L_\a(\Z)$ by functions of class $B_{\be,1}^1(\R)$ and by applying Theorem \ref{Besdost}, though in 
\cite{FN} a completely different proof was given.

\medskip

{\bf4. Functions of collections of commuting self-adjoint opetrators.} The study of functions of normal operators is equivalent to the study of functions of pairs of commuting self-adjoint operators. In \cite{NP} the results of  \cite{APPS} (see \S\:\ref{OLnaplBes} of this survey) were generalized to the case of functions 
of an arbitrary number of commuting self-adjoint operators. However, completely new methods were used.

Note also that in \cite{A3} some results on linear-fractional substitutions obtained in \cite{A1}
(see \S\:\ref{dlp} of this survey) were generalized to operator Lipschitz functions of several variables. In the multidimensional situation the role of  linear-fractional transformations is played by M\"obius transformations, i.e., compositions of finitely many inversions. 

\medskip

{\bf5. Lipschitz functions of collections of commuting self-adjoint operators.} In  \cite{KPSS} the results of \cite{PS} were generalized to the case of functions of $n$ commuting self-adjoint operators and Lipschitz type estimates in the norm of $\bS_p$, $1<p<\be$, were obtained for Lipschitz functions on $\R^n$.

\medskip

{\bf6. Functions of pairs of noncommuting self-adjoint operators.} 
The paper \cite{ANP} is devoted to the study of functions $f(A,B)$ of not necessarily commuting self-adjoint operators $(A,B)$. Such functions are defined in terms of double operator integrals. The authors of \cite{ANP} studied the 
behavior of such functions under perturbations of the pair.
It turned out that unlike in the case of functions of commuting operators, Lipschitz type estimates in the operator norm and in the trace norm differ from each other dramatically. In particular, it was shown in \cite{ANP} that for $f\in B_{\be,1}^1(\R^2)$, the following inequality holds
$$
\|f(A_1,B_1)-f(A_2,B_2)\|_{\bS_p}\le\const\|f\|_{B_{\be,1}^1}
\max\{\|A_1-A_2\|_{\bS_p},\|B_1-B_2\|_{\bS_p}\}
$$
for $p\in[1,2]$.
Such an inequality was obtained earlier in \cite{APPS} for functions of commuting operators for $p\ge1$. However, in the case of functions of noncommuting operators, {\it this inequality is false for $p>2$ and for the operator norm}, see \cite{ANP}.

The main tool used in \cite{ANP} is triple operator integrals and certain {\it modified} Haagerup tensor products of $L^\be$ spaces that were introduced in \cite{ANP}.

\medskip

{\bf7. Operator Lipschitz functions and the Lifshits--Krein trace formula.} Let $A$ and $B$ be self-adjoint operators with trace class difference $A-B$. For each such pair, there is a unique real function $\xi$ in $L^1(\R)$ called the {\it spectral shift function} such that for sufficiently nice functions $f$ on $\R$, the following Lifshits--Krein trace formula holds
$$
\trace(f(A)-f(B))=\int_\R f'(t)\xi(t)\,dt
$$
(see \cite{Li} and \cite{Kr}). M.G. Krein showed in \cite{Kr} that this formula holds for functions $f$, whose derivative is the Fourier transform of a complex measure. In \cite{Pe3} the trace formula was extended to functions $f$ of Besov class $B_{\be,1}^1(\R)$. Theorem \ref{tsentrez} of this survey says that for the operator $f(A)-f(B)$ to be in trace class under the assumption that $A-B$ is in trace class, it is necessary and sufficient that $f$ be operator Lipschitz. Finally, in the recent paper \cite{Pe7} it was shown that for operator Lipschitz functions, the left-hand side of the Lifshits--Krein trace formula not only makes sense, but also coincides with its right-hand side. In other words, {\it the Lifshits--Krein trace formula holds for arbitrary self-adjoint operators $A$ and $B$ with trace class difference if and only if the function $f$ is operator Lipschitz}.

\medskip

To conclude the paper, we mention a recent survey \cite{Pe10}, in which applications of multiple operator integrals in different problems of perturbation theory are considered.

\medskip

\footnotesize

\noindent
\begin{tabular}{p{9cm}p{15cm}}
A.B. Aleksandrov & V.V. Peller \\
St.Petersburg Branch & Department of Mathematics \\
Steklov Institute of Mathematics  & Michigan State University \\
Fontanka 27, 191023 St.Petersburg & East Lansing, Michigan 48824\\
Russia&USA
\end{tabular}

\end{document}